\documentclass{aims}
\usepackage{amsmath}
\usepackage{paralist}
\usepackage{epsfig} 
\usepackage{epstopdf}
\usepackage[colorlinks=true]{hyperref}

\usepackage{amsfonts,amssymb,stmaryrd,amssymb}
\usepackage{graphicx,color}
\usepackage{multimedia}
\usepackage{float}
\usepackage{placeins}

\hypersetup{urlcolor=blue, citecolor=red}

\def\currentvolume{X}
 \def\currentissue{X}
  \def\currentyear{20XX}
   \def\currentmonth{XX}
    \def\ppages{X--XX}
     \def\DOI{10.3934/xx.xx.xx.xx}

\newtheorem{theorem}{Theorem}[section]

\newtheorem{lemma}[theorem]{Lemma}
\newtheorem{proposition}{Proposition}

\theoremstyle{definition}
\newtheorem{definition}[theorem]{Definition}
\newtheorem{remark}{Remark}

\newtheorem{TTheorem}{Theorem}[section]

\newcommand \p {\partial}

\newcommand \com {\mathrm{cof}}

\newcommand \R {\mathbb{R}}

\renewcommand \L {\mathrm{L}}
\newcommand \W {\mathrm{W}}

\newcommand \LL {\mathbf{L}}
\newcommand \MM {\mathbf{M}}
\newcommand \NN {\mathbf{N}}
\newcommand \GG {\mathbf{G}}
\renewcommand \H {\mathrm{H}}

\newcommand \I {\mathrm{I}}
\newcommand \Id {\mathrm{Id}}

\renewcommand \d {\mathrm{d}}

\renewcommand \div {\mathrm{div}}
\renewcommand \det {\mathrm{det}}
\newcommand \rot {\mathrm{curl}}

\newcommand \trace {\mathrm{trace}}

\title[Stabilization of a fluid-solid system: Part II]
      {Stabilization of a fluid-solid system, by the deformation of the self-propelled solid.\\ Part II: The nonlinear system.}

\author[S\'ebastien Court]{S\'ebastien Court}



\email{sebastien.court@math.univ-toulouse.fr}

\thanks{This work is partially supported by the ANR-project CISIFS, 09-BLAN-0213-03.}

\begin{document}
\maketitle

\centerline{\scshape S\'ebastien Court }
\medskip
{\footnotesize
 \centerline{Institut de Math\'ematiques de Toulouse}
   \centerline{Universit\'e Paul Sabatier}
   \centerline{118 route de Narbonne}
   \centerline{31062 Toulouse Cedex 9, FRANCE}
}

\bigskip

 \centerline{(Communicated by the associate editor name)}

\begin{abstract}
In this second part we prove that the full nonlinear fluid-solid system introduced in Part I is stabilizable by deformations of the solid that have to satisfy nonlinear constraints. Some of these constraints are physical and guarantee the {\it self-propelled} nature of the solid. The proof is based on the boundary feedback stabilization of the linearized system. From this boundary feedback operator we construct a deformation of the solid which satisfies the aforementioned constraints and which stabilizes the nonlinear system. The proof is made by a fixed point method.
\end{abstract}

\section{Introduction}
\textcolor{black}{The fluid-solid system aforementioned is completely described in Part I (section 1.1). For more clarity, let us remind it briefly:
\begin{eqnarray}
\frac{\p u}{\p t} - \nu \Delta u + (u\cdot \nabla )u + \nabla p  =  0 , & \quad &
x \in \mathcal{F} (t),\quad t\in (0,\infty), \label{prems} \\
\div \ u   =  0 ,  & \quad &  x\in \mathcal{F} (t), \quad t\in (0,\infty), \label{deus}
\end{eqnarray}
\begin{eqnarray}
u  =  0 , & \quad & x \in \p \mathcal{O} ,\quad t\in (0,\infty),  \label{trois} \\
u  =  h'(t) +  \omega (t)\wedge(x-h(t)) + w(x,t) , & \quad &  x\in \p \mathcal{S}(t),\quad t\in
(0,\infty), \label{quatre}
\end{eqnarray}
\begin{eqnarray}
M h''(t)  =  - \int_{\p \mathcal{S}(t)} \sigma(u,p) n \d \Gamma  , & \quad & t\in (0,\infty), \label{cinq} \\
\left(I\omega\right)' (t)  = - \int_{\p \mathcal{S}(t)} (x-h(t))\wedge  \sigma(u,p) n \d \Gamma, & \quad & t\in(0,\infty),
\label{six}
\end{eqnarray}
\begin{eqnarray}
u(y,0)  =  u_0 (y), \  y\in \mathcal{F}(0) ,\quad h'(0)=h_1 \in \R^3 ,\quad \omega(0) = \omega_0 \in \R^3 ,
\end{eqnarray}
where
\begin{eqnarray}
\mathcal{S}(t) = h(t) + \mathbf{R}(t)X^{\ast}(\mathcal{S}(0),t), & \quad & \mathcal{F}(t) = \mathcal{O} \setminus \overline{\mathcal{S}(t)},
\end{eqnarray}
and where the velocity $w$ is defined by the following change of frame
\begin{eqnarray}
w(x,t) & = & \mathbf{R}(t)\ w^{\ast}\left(\mathbf{R}(t)^{T}(x-h(t)), t\right), \quad x\in \mathcal{S}(t). \label{wwwstar} \label{ders}
\end{eqnarray}
The unknowns of system \eqref{prems}--\eqref{ders} are the quadruplet $(u,p,h',\omega)$, for which the function $\omega$ is related to the rotation $\mathbf{R}$ by:
\begin{eqnarray*}
\left\{\begin{array} {ccccc}
\displaystyle \frac{\d \mathbf{R}}{\d t} & = & \mathbb{S}\left( \omega\right) \mathbf{R} \\
\mathbf{R}(0) & = & \I_{\R^3}
\end{array} \right.,
\quad \text{with }
\mathbb{S}(\omega) = \left(
\begin{matrix}
0 & -\omega_3 & \omega_2 \\
\omega_3 & 0 & -\omega_1 \\
-\omega_2 & \omega_1 & 0
\end{matrix} \right). \label{rotation}
\end{eqnarray*}
The control is the mapping $X^{\ast}$, which is related to the function $w^{\ast}$ by:
\begin{eqnarray*}
\frac{\p X^{\ast}}{\p t}(y,t) = w^{\ast}(X^{\ast}(y,t),t), \quad X^{\ast}(y,0) = y-h(0)=y, \quad y\in \mathcal{S}(0). \label{cauchystar}
\end{eqnarray*}
Let us also remind the hypotheses that the mapping $X^{\ast}$ has to satisfy:
\begin{description}
\item[H1] For all $t\geq 0$, $X^{\ast}(\cdot,t)$ is a $C^{1}$-diffeomorphism from $\mathcal{S}(0)$ onto $\mathcal{S}^{\ast}(t)$. \\
\item[H2] The volume of the solid has to be constant with respect to the time. That is equivalent to say that
\begin{eqnarray*}
 \int_{\p \mathcal{S}(0)} \frac{\p X^{\ast}}{\p t} \cdot \left( \com \nabla X^{\ast} \right)n\d \Gamma =  0.
\end{eqnarray*}
\item[H3] The deformation of the solid does not modify its linear momentum:
\begin{eqnarray*}
\int_{\mathcal{S}(0)} \rho_{\mathcal{S}}(y,0) X^{\ast}(y,t) \d y & = & 0.
\end{eqnarray*}
\item[H4] The deformation of the solid does not modify its angular momentum:
\begin{eqnarray*}
\int_{\mathcal{S}(0)} \rho_{\mathcal{S}}(y,0) X^{\ast}(y,t)\wedge \frac{\p X^{\ast}}{\p t}(y,t) \d y & = & 0.
\end{eqnarray*}
\end{description}
For more explanations on the model we invite the reader to refer to the introduction of Part I. Let us introduce the main result of this two-part work, for which the first part is only a step.
}

\subsection{The main result and the methods}
The existence of global strong solutions for system \eqref{prems}--\eqref{ders} has been proven in \cite{SMSTT} in dimension 2, and more recently in dimension 3 \textcolor{black}{in \cite{Court},} for small data in a framework where the deformations of the solid are limited in regularity. The main result of this second part is Theorem \ref{maintheorem}, that we can state as follows:
\begin{TTheorem}
Let be $\lambda >0$. Then for $(u_0,h_1,\omega_0)$ small enough in $\mathbf{H}^1(\mathcal{F}) \times \R^3 \times \R^3$, there exists a deformation $X^{\ast}$ satisfying the hypotheses {\bf H1}--{\bf H4} given above and also
\begin{eqnarray*}
e^{\lambda t}\frac{\p X^{\ast}}{\p t} & \in & \L^2(0,\infty;\mathbf{H}^3(\mathcal{S}(0))) \cap \H^1(0,\infty;\mathbf{H}^1(\mathcal{S}(0)))
\end{eqnarray*}
such that the solution $(u,p,h',\omega)$ of system \eqref{prems}--\eqref{ders} satisfies
\begin{eqnarray*}
\|e^{\lambda t} (u,p,h',\omega)\|
_{\H^{2,1}(Q_{\infty})\times \L^2(0,\infty;\mathbf{H}^1(\mathcal{F}(t))) \times \H^1(0,\infty;\R^3) \times \H^1(0,\infty;\R^3)}
& \leq & C.
\end{eqnarray*}
\end{TTheorem}
The notation are explained in section \ref{secdef}, below. The proof of this theorem is based on the preliminary stabilization of the linearized system in Part I. The idea is the following: If the perturbations of the system (which are represented by the initial conditions) are small enough, the behavior of the nonlinear system is close to the evolution of the linearized system. The strategy we follow in this second part consists in rewriting system \eqref{prems}--\eqref{ders} in cylindrical domains (that means in domains whose the space component does not depend on time), by defining a change of variables and a change of unknowns. Then we focus on the nonlinear system thus obtained. The feedback boundary control obtained in Part I enables us to stabilize the nonhomogeneous linear part of this system, whereas we have now to consider as a control function a deformation of the solid which satisfies the nonlinear constraints {\bf H1}--{\bf H4} given above. So the difficulty is to define properly from this boundary feedback control a deformation of the solid which stabilizes the full nonlinear system. For that the boundary feedback control is extended inside the solid as it is done in the section 6 of Part I; The deformation thus obtained satisfies the linearized version of the constraints stated in the hypotheses {\bf H1}--{\bf H4}. Then we project the displacement associated with this mapping on a set representing the displacements satisfying the nonlinear constraints. The deformation obtained is said to be {\it admissible}, that is to say it lies in a well-chosen functional space and it obeys the hypotheses {\bf H1}--{\bf H4} (see Definition \ref{defcontrol}). This projection method enables us to decompose the deformation velocity - on the fluid-solid interface - into two parts: The first part of this decomposition stabilizes the linear part of the nonlinear system, whereas the remaining part satisfies good Lipschitz properties with respect to the boundary feedback control. This point is essential if we want to prove by a fixed point method that such a deformation stabilizes the nonlinear system.\\
Besides the technical aspects of this work, a particular difficulty is the consideration of a control that has to satisfy nonlinear physical constraints. The originality of our contribution can be read in a perspective which concerns the study and the control of the swim of a deformable structure at an intermediate Reynolds number. By controlling the velocity of the environing fluid in a bounded domain, the solid stabilizes to zero the full system which is already dissipative (because of the viscosity), but the strength lies in the fact that this stabilization is obtained for all prescribed exponential decay rate.\\
Other papers treat of this issue: Let us quote for instance the work of Khapalov \cite{Khapalov1, Khapalov2} where the incompressible Navier-Stokes equations are considered. More recently let us mention the work of Chambrion \& Munnier \cite{Munnier1, Munnier2} dealing with perfect fluids, and where geometric methods have been used. Concerning the swim at a low Reynolds number, let us quote the recent paper \cite{Loheac} where prescribed types of deformations are considered. At a high Reynolds number, the work of Glass \& Rosier \cite{Glass} consists in applying the Coron return method in order to prove the local controllability of the position and the velocity of a boat, which is able to impose a velocity on a part of its boundary in order to move itself in a fluid satisfying the incompressible Euler equations. The result we give concerns only the velocities of the system.

\subsection{Plan}
The functional framework is given in section \ref{secdef} for the unknowns as well as for the control function and the changes of variables. Then the change of variables and the change of unknowns are introduced in section \ref{secchange}. Technical results in relation with the changes of variables are stated and proven in Appendix~A and Appendix~B. In section \ref{secchange} we also give the nonlinear system written in cylindrical domains, in a form which enables us to make the connection with the linearized system. The main result of Part I is used in section \ref{linearsec} where the feedback stabilization of the nonhomogeneous linear system is treated. Then we construct in section \ref{secdecompcontrol} an {\it admissible} deformation which is supposed to stabilize the full nonlinear system in section \ref{secnonlinear}. The link with the actual unknowns is made in section \ref{secconclusion} \textcolor{black}{that we can consider} as a conclusion.

\section{Definitions and notation} \label{secdef}
\textcolor{black}{Like in Part I, we denote by
\begin{eqnarray*}
\mathcal{F} = \mathcal{F}(0) & \text{and} & \mathcal{S}= \mathcal{S}(0)
\end{eqnarray*}
the domains occupied at time $t=0$ by the fluid and the solid respectively. The domain $\mathcal{S}$ is assumed to be smooth enough. The inertia matrix of the solid at $t=0$ is still denoted by $I_0$. We set $\mathcal{S}^{\ast}(t) = X^{\ast}(\mathcal{S},t)$, $\mathcal{F}^{\ast}(t) = \mathcal{O} \setminus \overline{\mathcal{S}^{\ast}(t)}$,
\begin{eqnarray*}
S^0_{\infty} = \mathcal{S} \times (0,\infty), & \quad & Q_{\infty}^0 = \mathcal{F} \times (0,\infty)
\end{eqnarray*}
and
\begin{eqnarray*}
Q_{\infty}  =  \displaystyle \bigcup_{t\geq 0} \mathcal{F}(t) \times \{ t \}.
\end{eqnarray*}
Let us keep in mind that the density $\rho_{\mathcal{S}}$ at time $t=0$ is assumed to be constant with respect to the space:
\begin{eqnarray*}
\rho_{\mathcal{S}}(y,0) = \rho_{\mathcal{S}} > 0.
\end{eqnarray*}
For more simplicity we have assumed - without loss of generality -  that $h_0 = 0$. This enables us to write
\begin{eqnarray*}
\int_{\mathcal{S}}y \d y & = & 0.
\end{eqnarray*}
}
The cofactor matrix associated with some matrix field $\mathbf{A}$ is denoted by $\com \mathbf{A}$. Let us keep in mind that when this matrix is invertible we have the property
\begin{eqnarray*}
\com \mathbf{A}^T & = & (\det \mathbf{A}) \mathbf{A}^{-1}.
\end{eqnarray*}
Let us introduce some functional spaces. \textcolor{black}{For the multidimensional Sobolev spaces} we use the notation
\begin{eqnarray*}
 \mathbf{H}^s(\Omega) = [\H^s(\Omega)]^d \ \text{or } [\H^s(\Omega)]^{d^k},
\end{eqnarray*}
for some positive integer $k$, for all bounded domain $\Omega$ of $\R^2$ or $\R^3$.

\subsection{Functional framework for the unknowns}
The velocity $u$ will be considered in the following functional space
\begin{eqnarray*}
\H^{2,1}(Q_{\infty}) & = & \L^2(0,\infty;\mathbf{H}^{2}(\mathcal{F}(t))) \cap \H^{1}(0,\infty;\mathbf{L}^2(\mathcal{F}(t)))
\end{eqnarray*}
that we endow and define with the norm given by
\begin{eqnarray*}
 \| u\|^2_{\H^{2,1}(Q_{\infty})} & = & \int_0^{\infty} \left\|u(\cdot, t) \right\|_{\mathbf{H}^2(\mathcal{F}(t))}^2 \d t + \int_0^{\infty} \left\|\frac{\p u}{\p t}(\cdot, t) \right\|_{\mathbf{L}^{2}(\mathcal{F}(t))}^2 \d t \ < \ + \infty .
\end{eqnarray*}
Analogously we can define spaces of type $\H^{s_1}(0,T;\H^{s_2}(\Omega(t)))$ and $\H^{s_1}(0,T;\mathbf{H}^{s_2}(\Omega(t)))$ for all time-depending domain $\Omega(t)$, where $s_1$ and $s_2$ are non-negative integers. The pressure $p$ will be considered in $\L^2(0,T;\H^1(\mathcal{F}(t)))$; At each time $t$ it is determined up to a constant that we fix such that $\int_{\mathcal{F}(t)} p = 0$. Thus in particular from the Poincar\'e-Wirtinger inequality the pressures $P$ defined in $\mathcal{F}$ can be estimated in $\H^1(\mathcal{F})$ as follows\footnote{In the following the symbols $C$ or $C_0$ will denote more generally a generic positive constant which does not depend on time or on the unknowns.}
\begin{eqnarray*}
\| P \|_{\H^1(\mathcal{F})} & \leq & C \| \nabla P \|_{\L^2(\mathcal{F})}.
\end{eqnarray*}
The same estimate will be considered for other functions which play the role of a pressure in $\mathcal{F}(0)$.\\
More classically in the cylindrical domain $Q_{\infty}^0$ we set
\begin{eqnarray*}
\H^{2,1}(Q_{\infty}^0) & = & \L^{2}(0,\infty; \mathbf{H}^2(\mathcal{F})) \cap \H^{1}(0,\infty; \mathbf{L}^2(\mathcal{F})),
\end{eqnarray*}
and we keep in mind the continuous embedding
\begin{eqnarray*}
\H^{2,1}(Q_{\infty}^0) & \hookrightarrow & \L^{\infty}(0,\infty; \mathbf{H}^1(\mathcal{F})).
\end{eqnarray*}

\subsection{Functional framework for the control and the changes of variables}
Let us introduce some functional spaces for the solid displacements - defined in the domain $\mathcal{S}$ - and the changes of variables which will be defined in the domain $\mathcal{F}$. We will mainly consider mappings $X^{\ast}$ satisfying $X^{\ast}(\cdot,0) = \Id_{\mathcal{F}}$, and thus we will consider the displacements
\begin{eqnarray*}
Z^{\ast} = X^{\ast} - \Id_{\mathcal{S}} & \in & \mathcal{W}_{\lambda}(S_{\infty}^0),
\end{eqnarray*}
where the space $\mathcal{W}_{\lambda}(S_{\infty}^0)$ is defined as follows
\begin{eqnarray*}
Z^{\ast} \in \mathcal{W}_{\lambda}(S_{\infty}^0) \Leftrightarrow \left\{ \begin{array} {lll}
\displaystyle e^{\lambda t}\frac{\p Z^{\ast}}{\p t} \in \L^2(0,\infty;\mathbf{H}^3(\mathcal{S})) \cap \H^{1}(0,\infty;\mathbf{H}^1(\mathcal{S})), \\
\hfill \\
Z^{\ast}(y,0) = 0, \ \displaystyle \frac{\p Z^{\ast}}{\p t}(y , 0) = 0 \quad \forall y\in \mathcal{S}.
\end{array} \right.
\end{eqnarray*}
We endow it with the scalar product
\begin{eqnarray*}
\langle Z^{\ast}_1 ; Z^{\ast}_2 \rangle_{\mathcal{W}_{\lambda}(S_{\infty}^0)} & := &
\left\langle e^{\lambda t}\frac{\p Z^{\ast}_1}{\p t} ; e^{\lambda t}\frac{\p Z^{\ast}_2}{\p t} \right\rangle_{\L^2(0,\infty;\mathbf{H}^3(\mathcal{S})) \cap \H^{1}(0,\infty;\mathbf{H}^1(\mathcal{S}))}
\end{eqnarray*}
which makes it be a Hilbert space, because of the continuous embedding
\begin{eqnarray}
\mathcal{W}_{\lambda}(S_{\infty}^0) & \hookrightarrow & \L^{\infty}(0,\infty;\mathbf{H}^3( \mathcal{S})) \cap \W^{1,\infty}(0,\infty;\mathbf{H}^1( \mathcal{S})). \label{ineqstar}
\end{eqnarray}
Indeed, for $X^{\ast}-\Id_{\mathcal{S}} \in \mathcal{W}_{\lambda}(S_{\infty}^0)$, we have the following estimates
\begin{eqnarray*}
\| X^{\ast} - \Id_{\mathcal{S}} \|_{\L^{\infty}(0,\infty;\mathbf{H}^3(\mathcal{S}))}
& \leq & \int_0^{\infty}e^{-\lambda s}\left\|e^{\lambda s} \frac{\p X^{\ast}}{\p t}(\cdot,s) \right\|_{\mathbf{H}^3(\mathcal{S})} \d s, \\
& \leq & \frac{1}{\sqrt{2\lambda}}\left\|e^{\lambda t} \frac{\p X^{\ast}}{\p t} \right\|_{\L^2(0,\infty;\mathbf{H}^3(\mathcal{S}))}, \\
\left\| \frac{\p X^{\ast}}{\p t} \right\|_{\L^{\infty}(0,\infty;\mathbf{H}^1(\mathcal{S}))}
& \leq & \int_0^{\infty}e^{-\lambda s}\left\|e^{\lambda s} \frac{\p^2 X^{\ast}}{\p t^2}(\cdot,s) \right\|_{\mathbf{H}^1(\mathcal{S})} \d s, \\
& \leq & \frac{1+\lambda}{\sqrt{2\lambda}}\left\|e^{\lambda t} \frac{\p X^{\ast}}{\p t} \right\|_{\H^1(0,\infty;\mathbf{H}^1(\mathcal{S}))}.
\end{eqnarray*}
In the rewriting of system \eqref{prems}--\eqref{ders} in cylindrical domains, and in the final fixed point method, the mapping which has the most important role is denoted by $\tilde{X}$. We could consider the same type of functional space for this mapping in $Q_{\infty}^0$, but we will only need to consider - and we will only have - estimates of $\tilde{X} - \Id_{\mathcal{F}}$ in the space
\begin{eqnarray*}
\L^{\infty}(0,\infty;\mathbf{H}^3( \mathcal{F})) \cap \W^{1,\infty}(0,\infty;\mathbf{H}^1( \mathcal{F})).
\end{eqnarray*}
Thus for more clarity we set
\begin{eqnarray*}
\mathcal{H}(S_{\infty}^0) & = & \L^{2}(0,\infty;\mathbf{H}^3( \mathcal{S})) \cap \H^{1}(0,\infty;\mathbf{H}^1( \mathcal{S})), \\
\mathcal{H}(Q_{\infty}^0) & = & \L^{2}(0,\infty;\mathbf{H}^3( \mathcal{F})) \cap \H^{1}(0,\infty;\mathbf{H}^1( \mathcal{F})), \\
\tilde{\mathcal{W}}(Q_{\infty}^0) & = & \L^{\infty}(0,\infty;\mathbf{H}^3( \mathcal{F})) \cap \W^{1,\infty}(0,\infty;\mathbf{H}^1( \mathcal{F})).
\end{eqnarray*}
The main reason for which we choose to consider the aforementioned Lagrangian mappings in such functional spaces is the following: The changes of variables $\tilde{X}$ will be - indirectly - obtained through extensions of the deformations of the solid $X^{\ast}$; We will need to consider displacements $\tilde{X}-\Id_{\mathcal{F}}$ which have the regularity indicated in the spaces given above, and for that we will have to consider displacements of the solid which lie in a Hilbert space (for the projection method) and which satisfy at least
\begin{eqnarray*}
e^{\lambda t}\frac{ \p X^{\ast}}{\p t}_{\left| \p \mathcal{S} \right.} & \in &
\L^{2}(0,\infty;\mathbf{H}^{5/2}(\p \mathcal{S})) \cap \H^{1}(0,\infty;\mathbf{H}^{1/2}(\p \mathcal{S})).
\end{eqnarray*}

\subsection{Definitions of stabilization for the nonlinear system}
Let us specify the notion of {\it admissible} control and the definition we give to the stabilization of the main system.

\begin{definition} \label{defcontrol}
Let be $\lambda >0$. A deformation $X^{\ast}$ is said {\it admissible} for the nonlinear system \eqref{prems}--\eqref{ders} if the displacement obeys
\begin{eqnarray*}
X^{\ast} - \Id_{\mathcal{S}} & \in & \mathcal{W}_{\lambda}(S_{\infty}^0),
\end{eqnarray*}
if $X^{\ast}(\cdot ,t)$ is a $C^1$-diffeomorphism from $\mathcal{S}$ onto $\mathcal{S}^{\ast}(t)$ for all $t \geq 0$, and if for all $t\geq 0$ it satisfies the following hypotheses
\begin{eqnarray}
\int_{\p \mathcal{S}} \frac{\p X^{\ast}}{\p t}(y,t) \cdot \left(\com \nabla X^{\ast}(y,t)\right) n  \d \Gamma(y)  =  0, & & \label{const1} \\
\int_{\mathcal{S}} X^{\ast}(y,t) \d y = 0, & & \label{const2} \\
\int_{\mathcal{S}} X^{\ast}(y,t) \wedge \frac{\p X^{\ast}}{\p t}(y,t) \d y = 0. & & \label{const3}
\end{eqnarray}
\end{definition}

\begin{remark}
First, the constraint which forces $X^{\ast}$ to be a $C^1$-diffeomorphism can be relaxed if the control $X^{\ast}$ stays close to the identity $\Id_{\mathcal{S}}$. Indeed, in this work the data are assumed to be small enough, so that we will lead to consider the displacement $X^{\ast} - \Id_{\mathcal{S}}$ small enough in the space $\mathcal{W}_{\lambda}(S_{\infty}^0)$ and so in $\L^{\infty}(0,\infty;\mathbf{H}^3(\mathcal{S}))$; Thus we can assume that this constraint is always satisfied.
\end{remark}

\begin{definition} \label{defstab}
We say that system \eqref{prems}--\eqref{ders} is stabilizable with an arbitrary exponential decay rate if for each $\lambda >0$ there exists an {it admissible} deformation $X^{\ast}$ (in the sense of Definition \ref{defcontrol}) and a positive constant $C$ - depending only on $u_0$, $h_1$ and $\omega_0$ - such that the solution $(u,p,h', \omega)$ of system \eqref{prems}--\eqref{ders} satisfies
\begin{eqnarray*}
\|e^{\lambda t} (u,p,h',\omega)\|
_{\H^{2,1}(Q_{\infty})\times \L^2(0,\infty;\mathbf{H}^1(\mathcal{F}(t))) \times \H^1(0,\infty;\R^3) \times \H^1(0,\infty;\R^3)}
& \leq & C.
\end{eqnarray*}
\end{definition}

\section{Change of variables and change of unknowns} \label{secchange}
In order to use a change of unknowns which will enable us to rewrite the main system \eqref{prems}--\eqref{ders}, we first extend to the whole domain $\overline{\mathcal{O}}$ the mappings $X_{\mathcal{S}}(\cdot,t) = h(t) + \mathbf{R}(t)X^{\ast}(\cdot,t)$, initially defined on $\mathcal{S}$. Then by this extension - denoted by $X$ - we will define new unknowns whose interest lies in the fact that they are supposed to satisfy a system written in cylindrical domains, that is to say domains whose the space component does not depend on time.

\subsection{The change of variables}
For a vector field $h \in \H^2(0,\infty;\R^3)$ and a rotation $\mathbf{R} \in \H^2(0,\infty;\R^9)$ which provides an angular velocity $\omega \in \H^1(0,\infty;\R^3)$, and for an {\it admissible} deformation $X^{\ast}$ - in the sense of Definition \ref{defcontrol} - the aim of this subsection is to construct a mapping $X$ which has to satisfy the properties
\begin{eqnarray*}
\left\{ \begin{array} {lll}
\det \nabla X = 1, & \quad & \text{in } \mathcal{F} \times (0,\infty), \\
X = h+\mathbf{R} X^{\ast} , & \quad & \text{on } \p \mathcal{S} \times (0,\infty), \\
X = \Id_{\p \mathcal{O}}, & \quad & \text{on } \p \mathcal{O} \times (0,\infty),
\end{array} \right.
\end{eqnarray*}
such that for all $t\geq 0$ the function $X(\cdot,t)$ maps $\mathcal{F}$ onto $\mathcal{F}(t)$, $\p \mathcal{S}$ onto $\p \mathcal{S}(t)$, and leaves invariant the boundary $\p \mathcal{O}$. For that, let us first define an intermediate extension. We can extend the mapping $X^{\ast}$ to $\mathcal{F}$ as follows:\\
If we assume that the mapping $X^{\ast}$ satisfies the hypothesis {\bf H2}, that is to say for all $t \geq 0$ we have the condition
\begin{eqnarray*}
\int_{\p \mathcal{S}} \frac{\p X^{\ast}}{\p t} \cdot \left( \com \nabla X^{\ast} \right)n\d \Gamma & = &  0,
\end{eqnarray*}
and if we also assume that the function $\displaystyle (y,t) \mapsto e^{\lambda t} \frac{\p X^{\ast}}{\p t}$ is small enough in $\L^{2}(0,\infty;\mathbf{H}^3(\mathcal{S})) \cap \H^{1}(0,\infty;\mathbf{H}^1(\mathcal{S}))$, then there exists a mapping $\overline{X}^{\ast} \in \L^{\infty}(0,\infty;\mathbf{H}^3(\mathcal{F})) \cap \W^{1,\infty}(0,\infty;\mathbf{H}^1(\mathcal{F}))$ satisfying
\begin{eqnarray*}
\left\{ \begin{array} {lcl}
\det \nabla \overline{X}^{\ast}  = 1 & \quad & \text{in } \mathcal{F} \times (0,\infty), \\
\overline{X}^{\ast} = X^{\ast} & \quad & \text{on } \p \mathcal{S} \times (0,\infty), \\
\overline{X}^{\ast} = \Id_{\p \mathcal{O}} & \quad & \text{on } \p \mathcal{O} \times (0,\infty). 
\end{array} \right.
\end{eqnarray*}

The existence and other properties of such an extension $\overline{X}^{\ast}$ are summed up in the statement of Proposition \ref{lemmaxtension}, in Appendix~A. From this intermediate extension, the purpose in now to define the extension $X$ aforementioned.\\
The mapping $X_{\mathcal{S}}(\cdot,t)$ is obtained from $X^{\ast}(\cdot,t)$ by composing it to the left by the rigid displacement
\begin{eqnarray*}
X^{R}(\cdot,t) \ : \  x^{\ast} \mapsto h(t) + \mathbf{R}(t)x^{\ast}.
\end{eqnarray*}
We cannot do the same thing for obtaining the mapping $X$ from $\overline{X}^{\ast}$ because of the boundary condition on $\p \mathcal{O}$ that has to be preserved. Thus we define an extension of $X^{R}(\cdot,t)$ to the whole domain. For that we can use the same process which has been introduced in \cite{TT}, and thus construct $\overline{X}^R(\cdot,t)$ which satisfies for all $t\geq 0$ the following properties:
\begin{eqnarray*}
\left\{ \begin{array} {ll}
\det \nabla \overline{X}^{R}(\cdot,t) = 1 & \text{ in } \R^2 \text{ or } \R^3, \label{ch4_ext1}\\
\overline{X}^{R}(\cdot,t) = X^{R}(\cdot,t) &  \text{ on } \p \left(X^{\ast}(\mathcal{S},t)\right), \\
\overline{X}^{R}(\cdot,t) = \Id_{\p \mathcal{O}} & \text{ on } \p \mathcal{O}. \label{ch4_ext3}
\end{array} \right.
\end{eqnarray*}
So we define $X$ as follows
\begin{eqnarray*}
X(\cdot,t) & = & \overline{X}^R(\cdot,t) \circ \overline{X}^{\ast}(\cdot,t),
\end{eqnarray*}
and we denote by $Y(\cdot,t)$ its inverse for all $t \geq 0$.\\

The constructions of the mappings aforementioned are quite technical, that is why we develop the details of these constructions only in Appendix~A, in the same time as the regularity deduced on $X$. However, let us note an important point:

\begin{remark}
In Appendix~A, the definition of the mapping $X$ is conditioned by a smallness assumption on the solid deformation $X^{\ast}$. Since the deformation which will be chosen in sections \ref{secnonlinear} and \ref{secconclusion} for stabilizing the full nonlinear system will actually depend only on - and be controlled by - the initial data $(u_0,h_1,\omega_0)$ that we will assume small enough, it is possible to proceed like this.\\
Note that the definition of such a change of variables is also done in \cite{SMSTT} or in \cite{Court}. The change of variables we use in this paper has the same properties as the ones constructed in these articles. But the way we proceed here is not the same: In \cite{SMSTT} the mapping $X$ is constructed by extending the Eulerian velocity $w$ to the fluid part, but this means is not suitable in a framework where the role of the Lagrangian mapping representing the deformation of the solid is central and where its regularity is limited; Indeed the velocity $w(\cdot,t)$ is defined through $w^{\ast}(\cdot,t)$ which is itself defined on the domain $X^{\ast}(\mathcal{S},t)$. Concerning the means used in \cite{Court}, the problem solved for constructing the mapping $X$ in this paper is similar to ours, but it requires a smallness assumption on the time existence. Here the hypothesis we make is the smallness of $e^{\lambda t}\frac{\p X^{\ast}}{\p t}$ in an infinite time horizon space.
\end{remark}

\subsection{Rewriting the main system in cylindrical domains}
We use the change of variables given above in order to transform system \eqref{prems}--\eqref{ders} into a system which deals with non-depending time domains. For that we set the change of unknowns
\begin{eqnarray}
\begin{array} {lll}
\tilde{u} (y,t) =  \mathbf{R}(t)^T u(X(y,t),t), & \quad & u(x,t)  =  \mathbf{R}(t) \tilde{u} (Y(x,t),t), \label{tildeu} \\
\tilde{p} (y,t)  =  p(X(y,t),t), & \quad & p(x,t)  =  \tilde{p}(Y(x,t),t), \label{tildep}
\end{array}
\end{eqnarray}
for $x\in \overline{\mathcal{F}(t)}$ and $y\in \overline{\mathcal{F}}$, and
\begin{eqnarray}
\tilde{h}'(t)  =  \mathbf{R}(t)^T h'(t), & \quad & \tilde{\omega}(t)  =  \mathbf{R}(t)^T \omega(t). \label{tildeh}
\end{eqnarray}

\begin{remark} \label{remarkc}
Let us notice that if $\tilde{h}'$ and $\tilde{\omega}$ are given, then by using the second equality of \eqref{tildeh} we see that $\mathbf{R}$ satisfies the Cauchy problem
\begin{eqnarray*}
\begin{array} {ccccc}
\displaystyle \frac{\d}{\d t}(\mathbf{R}) & = & \mathbb{S}\left(
\mathbf{R}\tilde{\omega}\right) \mathbf{R} & = & \mathbf{R} \mathbb{S}\left(\tilde{\omega}\right) \\
\mathbf{R}(t=0) & = & \I_{\R^3}, & &
\end{array}
\quad \text{with }
\mathbb{S}(\tilde{\omega}) = \left(
\begin{matrix}
0 & -\tilde{\omega}_3 & \tilde{\omega}_2 \\
\tilde{\omega}_3 & 0 & -\tilde{\omega}_1 \\
-\tilde{\omega}_2 & \tilde{\omega}_1 & 0
\end{matrix} \right). \label{rotationpb}
\end{eqnarray*}
So $\mathbf{R}$ is determined in a unique way. Thus it is obvious to see that, in \eqref{tildeh}, $h'$ and $\omega$ are also
determined in a unique way. Moreover, since we have
\begin{eqnarray*}
u(x,t)  =  \mathbf{R}(t) \tilde{u}(Y(x,t) ,t), & \quad & p(x,t)  =  \tilde{p}(Y(x,t) ,t),
\end{eqnarray*}
and since the mapping $Y$ depends only on $h$, $\omega$ and the control $X^{\ast}$, we finally see that if $(\tilde{u}, \tilde{p}, \tilde{h}', \tilde{\omega})$ is given, then $(u,p,h',\omega)$ is determined in a unique way.
\end{remark}

For what follows, it is convenient to define the mappings
\begin{eqnarray}
\begin{array} {ll}
\tilde{X}(y,t) =  \mathbf{R}(t)^T(X(y,t) - h(t)), & \quad (y,t) \in \mathcal{F}\times (0,\infty),  \label{definitilde} \\
&  \\
\tilde{Y}(\tilde{x},t) = Y(h(t) + \mathbf{R}(t)\tilde{x},t), & \quad (\tilde{x},t) \in
\displaystyle \bigcup_{t \geq 0} \mathbf{R}(t)^T\left(\mathcal{F}(t)-h(t)\right)\times \{ t \}.
\end{array}
\end{eqnarray}
The regularity, dependence with respect to the unknowns, and estimates for the mappings $\tilde{X}$ and $\tilde{Y}$ are given in Appendix~B.\\

Then, like in \cite{Court}, system \eqref{prems}--\eqref{ders} whose the unknowns are $(u,p,h',\omega)$ is rewritten in the cylindrical domain $\mathcal{F}\times (0,\infty)$ as the following system whose the unknowns are $(\tilde{u}, \tilde{p}, \tilde{h}', \tilde{\omega})$:
\begin{eqnarray*}
\frac{\p \tilde{u}}{\p t} - \nu \LL \tilde{u} + \MM
(\tilde{u}, \tilde{h}', \tilde{\omega})+ \NN \tilde{u} + \tilde{\omega}(t)\wedge\tilde{u}+ \GG \tilde{p}
 =  0, & \quad & \text{in } \mathcal{F}\times (0,\infty), \label{premsfix} \\
\div \ \tilde{u}  =  g, & \quad & \text{in } \mathcal{F}\times (0,\infty),
\label{deusfix}
\end{eqnarray*}
\begin{eqnarray*}
\tilde{u} = 0 , & \quad & \text{in } \p \mathcal{O} \times (0,\infty), \label{troisfix} \\
\tilde{u}  =  \tilde{h}'(t) + \tilde{\omega} (t) \wedge X^{\ast}(y,t)+ \frac{\p X^{\ast}}{\p t}(y,t) ,
& \quad & (y,t)\in \p \mathcal{S}\times (0,\infty), \label{quatrefix}
\end{eqnarray*}
\begin{eqnarray*}
M \tilde{h}''(t) = - \int_{\p \mathcal{S}} \tilde{\sigma}(\tilde{u},
\tilde{p}) \nabla \tilde{Y}(\tilde{X})^T n \d \Gamma - M\tilde{\omega}(t)\wedge \tilde{h}'(t) , & \quad & t \in (0,\infty)   \label{cinqfix} \\
I^{\ast}(t)\tilde{\omega}' (t) =  -  \int_{\p \mathcal{S}} X^{\ast}(y,t) \wedge \left(
\tilde{\sigma}(\tilde{u},\tilde{p})\nabla \tilde{Y}(\tilde{X})^T n\right) \d \Gamma & \quad &  \nonumber \\
 - {I^{\ast}}'(t)\tilde{\omega}(t)+ I^{\ast}(t)\tilde{\omega}(t)\wedge\tilde{\omega}(t), & \quad & t \in (0,\infty) \label{sixfix}
\end{eqnarray*}
\begin{eqnarray*}
\tilde{u}(y,0)  =  u_0 (y), \  y\in \mathcal{F} , \quad \tilde{h}'(0)=h_1 \in \R^3 ,\quad \tilde{\omega}(0) = \omega_0 \in \R^3 , \label{dersfix}
\end{eqnarray*}
where
\begin{eqnarray*}
I^{\ast}(t) &  = & \rho_{\mathcal{S}}\left( \int_{\mathcal{S}}|X^{\ast}(y,t)|^2\d y\right)\I_{\R^3} \quad \text{in dimension 2}, \\
I^{\ast}(t) &  = & \rho_{\mathcal{S}} \int_{\mathcal{S}}\left(|X^{\ast}(y,t)|^2\I_{\R^3}-X^{\ast}(y,t)\otimes X^{\ast}(y,t) \right)\d y \quad \text{in dimension 3},
\end{eqnarray*}
where $[ \cdot ]_i$ specifies the i-th component of a vector
\begin{eqnarray*}
& & [ \LL  \tilde{u} ]_i(y,t)  =  [ \nabla \tilde{u}(y,t) \Delta \tilde{Y}(\tilde{X}(y,t),t)]_i  +  \nabla^2 \tilde{u}_i(y,t) : \left(\nabla \tilde{Y} \nabla \tilde{Y}^T \right)(\tilde{X}(y,t),t),
 \label{LL} \\
& & \MM (\tilde{u}, \tilde{h}', \tilde{\omega})(y,t) = -\nabla \tilde{u} (y,t) \nabla \tilde{Y}(\tilde{X}(y,t),t)
 \left(\tilde{h}'(t) + \tilde{\omega}(t) \wedge \tilde{X}(y,t) + \frac{\p \tilde{X}}{\p t}(y,t)\right) ,\nonumber \\  \label{MM} \\
& & \NN  \tilde{u}(y,t) = \nabla \tilde{u}(y,t) \nabla \tilde{Y}(\tilde{X}(y,t),t) \tilde{u}(y,t),  \label{NN} \\
& & \GG  \tilde{p} (y,t)  =  \nabla \tilde{Y}(\tilde{X}(y,t),t)^T \nabla \tilde{p}(y,t), \label{GG} \\
& & \tilde{\sigma}(\tilde{u},\tilde{p})(y,t) =
\nu\left(\nabla \tilde{u}(y,t) \nabla \tilde{Y}(\tilde{X}(y,t),t) + \nabla \tilde{Y}(\tilde{X}(y,t),t)^T \nabla \tilde{u}(y,t)^T \right) - \tilde{p}(y,t) \I_{\R^3} \nonumber
\end{eqnarray*}
and
\begin{eqnarray*}
g(y,t)  & = &  \trace\left( \nabla \tilde{u}(y,t) \left(\I_{\R^3} -
\nabla \tilde{Y}\left(\tilde{X}(y,t),t\right) \right) \right) \nonumber
\\
& = & \nabla \tilde{u}(y,t) : \left(\I_{\R^3} - \nabla \tilde{Y}\left(\tilde{X}(y,t),t\right)^T\right)
. \label{expg}
\end{eqnarray*}
\textcolor{black}{
Notice that we can actually express this nonhomogeneous divergence term as $g = \div \ G$, where
\begin{eqnarray*}
G(y,t) & = &  \left(\I_{\R^3} - \nabla \tilde{Y}(\tilde{X}(y,t),t)\right)\tilde{u}(y,t) .
\end{eqnarray*}
Indeed, we can calculate
\begin{eqnarray*}
\div \ G & = &  \left(\I_{\R^3} - \nabla \tilde{Y}(\tilde{X})\right)^T : \nabla \tilde{u} - \tilde{u} \cdot \div \left( \nabla \tilde{Y}(\tilde{X})^T\right),
\end{eqnarray*}
and the second term of this expression vanishes, because we have by construction
\begin{eqnarray*}
 \nabla \tilde{Y}(\tilde{X})^T  = \det(\nabla \tilde{X}) \nabla \tilde{Y}(\tilde{X})^T = \com(\nabla \tilde{X}) ,
\end{eqnarray*}
and the Piola identity (see \cite{Ciarlet}, the first part of the proof of Theorem 1.7-1 p.39) is nothing else than
\begin{eqnarray*}
\div \left(\com(\nabla \tilde{X})\right) & = & 0.
\end{eqnarray*}
}
For $\lambda > 0$, we now set the following change of unknowns:
\begin{eqnarray}
\hat{u} = e^{\lambda t}\tilde{u}, \quad \hat{p} = e^{\lambda t}\tilde{p}, \quad \hat{h}' = e^{\lambda t}\tilde{h}', \quad \hat{\omega} = e^{\lambda t}\tilde{\omega}
.  \label{chhat}
\end{eqnarray}
The idea of this second change of unknowns is the following: If we find a control $X^{\ast}$ such that the quadruplet $(\hat{u}, \hat{p}, \hat{h}', \hat{\omega})$ is bounded in some infinite-time horizon space, then the intermediate unknowns $(\tilde{u}, \tilde{p}, \tilde{h}', \tilde{\omega})$ will be stabilized with an exponential decay rate, and it will be sufficient to deduce from that the same property for the actual unknowns $(u,p,h',\omega)$ (see section \ref{secconclusion}).\\

The system satisfied by $(\tilde{u}, \tilde{p}, \tilde{h}', \tilde{\omega})$ is then transformed into
\begin{eqnarray}
\frac{\p \hat{u}}{\p t} - \lambda \hat{u} -\nu \Delta \hat{u} + \nabla \hat{p}
 =  F(\hat{u},\hat{p},\hat{h}',\hat{\omega}), & \quad &  \text{in } \mathcal{F}\times (0,\infty), \label{hpremsfix} \\
\div \ \hat{u}  =  \div \ G(\hat{u}), & \quad & \text{in } \mathcal{F}\times (0,\infty), \label{hdeusfix}
\end{eqnarray}
\begin{eqnarray}
\hat{u} = 0 , & \quad & \text{on } \p \mathcal{O}\times (0,\infty),  \label{htroisfix} \\
\hat{u}  =  \hat{h}'(t) + \hat{\omega} (t) \wedge y+ e^{\lambda t}\frac{\p X^{\ast}}{\p t} + W(\hat{\omega}) , & \quad & (y,t)\in \p \mathcal{S}\times (0,\infty), \label{hquatrefix}
\end{eqnarray}
\begin{eqnarray}
M \hat{h}'' - \lambda M \hat{h}'  =  - \int_{\p \mathcal{S}} \sigma(\hat{u},\hat{p}) n  \d \Gamma + F_M(\hat{u},\hat{p},\hat{h}',\hat{\omega}),
& \quad & \text{in } (0,\infty),  \label{hcinqfix} \\
I_0\hat{\omega}' (t) -\lambda I_0\hat{\omega}  =   -  \int_{\p \mathcal{S}} y \wedge \sigma(\hat{u},\hat{p}) n  \d \Gamma + F_I(\hat{u},\hat{p},\hat{\omega}), & \quad  & \text{in } (0,\infty),  \label{hsixfix}
\end{eqnarray}
\begin{eqnarray}
\hat{u}(y,0)  =  u_0 (y), \  y\in \mathcal{F} , \quad \hat{h}'(0)=h_1 \in \R^3 ,\quad \hat{\omega}(0) =\omega_0 \in \R^3 , \label{hdersfix}
\end{eqnarray}
with
\begin{eqnarray*}
F(\hat{u},\hat{p},\hat{h}',\hat{\omega}) & = & \nu (\LL - \Delta) \hat{u} - e^{-\lambda t}\MM (\hat{u}, \hat{h}', \hat{\omega}) - e^{-\lambda t}\NN \hat{u} - (\GG-\nabla) \hat{p} - e^{-\lambda t}\hat{\omega} \wedge \hat{u},  \\
G(\hat{u},\hat{h}',\hat{\omega}) & = & \left(\I_{\R^3}- \nabla \tilde{Y}(\tilde{X}(y,t),t)\right)\hat{u}, \label{0rhsG} \\
W(\hat{\omega}) & = & \hat{\omega} \wedge \left(X^{\ast} - \Id\right), \label{0rhsW} \\
F_M(\hat{u},\hat{p},\hat{h}',\hat{\omega}) & = & -Me^{-\lambda t} \hat{\omega} \wedge \hat{h}'(t) \nonumber \\
 & & - \nu\int_{\p \mathcal{S}}\left(\nabla \hat{u} \left(\nabla \tilde{Y}(\tilde{X}) - \I_{\R^3}\right) + \left({\nabla \tilde{Y}(\tilde{X})} - \I_{\R^3}\right)^T\nabla \hat{u}^T \right)\nabla \tilde{Y}(\tilde{X})^T n \d \Gamma \nonumber \\
 & & - \int_{\p \mathcal{S}}\sigma(\hat{u},\hat{p})\left(\nabla \tilde{Y}(\tilde{X})-\I_{\R^3}\right)^Tn\d \Gamma, \label{0rhsFM} \\
F_I(\hat{u},\hat{p},\hat{\omega}) & = & -\left(I^{\ast} - I_0\right) \hat{\omega}' + \lambda\left(I^{\ast}-I_0\right)\hat{\omega} - {I^{\ast}}'\hat{\omega} +e^{-\lambda t}I^{\ast}\hat{\omega} \wedge \hat{\omega} \nonumber \\
& &  - \nu\int_{\p \mathcal{S}}y\wedge \left(\nabla \hat{u} \left(\nabla \tilde{Y}(\tilde{X}) - \I_{\R^3}\right)
+ ({\nabla \tilde{Y}(\tilde{X})} - \I_{\R^3})^T\nabla \hat{u}^T\right)\nabla \tilde{Y}(\tilde{X})^T n \d \Gamma \nonumber \\
& & - \int_{\p \mathcal{S}}y\wedge \left(\sigma(\hat{u},\hat{p})(\nabla \tilde{Y}(\tilde{X})-\I_{\R^3})^Tn\right)\d \Gamma \nonumber \\
& & + \int_{\p \mathcal{S}}\left(X^{\ast}-\Id\right)\wedge \left(\tilde{\sigma}(\hat{u},\hat{p})
\nabla \tilde{Y}(\tilde{X})^T n\right)\d \Gamma. \label{0rhsFI}
\end{eqnarray*}

\begin{remark} \label{remarkcc}
An important remark is the following: Since system \eqref{prems}--\eqref{ders} and the system satisfied by $(\tilde{u}, \tilde{p}, \tilde{h}', \tilde{\omega})$ above  are equivalent, and since for an {\it admissible} control $X^{\ast}$ satisfying the constraint \eqref{const1} the compatibility condition is satisfied for system \eqref{prems}--\eqref{ders}, in system \eqref{hpremsfix}--\eqref{hdersfix} the underlying compatibility condition enables us to have automatically the following equality
\begin{eqnarray*}
\int_{\p \mathcal{S}} G(\hat{u},\hat{h}',\hat{\omega})\cdot n \d \Gamma & = &
\int_{\p \mathcal{S}} \left(e^{\lambda t}\frac{\p X^{\ast}}{\p t} + W(\hat{\omega})\right) \cdot n \d \Gamma
\end{eqnarray*}
as soon as $\hat{u} = 0$ on $\p \mathcal{O}$.
\end{remark}

\section{The stabilized nonhomogeneous linearized system} \label{linearsec}
\subsection{A nonhomogeneous linear system}
In this section, let us consider the nonhomogeneous linear system suggested by the writing of system \eqref{hpremsfix}--\eqref{hdersfix}:
\begin{eqnarray}
\frac{\p \hat{U}}{\p t} - \lambda \hat{U} - \nu \Delta \hat{U} + \nabla \hat{P}  =  \mathbb{F} , &  & \textrm{in $\mathcal{F}\times (0,\infty)$},  \label{linnhbb1}\\
\div \  \hat{U}  =  \div \ \mathbb{G}, &  & \textrm{in $\mathcal{F}\times (0,\infty)$}, \label{linnhbb2}
\end{eqnarray}
\begin{eqnarray}
\hat{U}  =  0 , & & \textrm{in $\p \mathcal{O}\times (0,\infty)$}, \\
\hat{U}=  \hat{H}'(t) + \hat{\Omega}(t) \wedge y + \zeta(y,t) + \mathbb{W}(y,t), &  & \textrm{in } \p \mathcal{S}\times(0,\infty), \label{linnhbb4}
\end{eqnarray}
\begin{eqnarray}
M \hat{H}''(t)-\lambda M \hat{H}'(t) = - \int_{\p \mathcal{S}} \sigma(\hat{U},\hat{P}) n \d \Gamma + \mathbb{F}_M, & &  t\in (0,\infty), \\
I_0 \hat{\Omega}'(t) - \lambda I_0\hat{\Omega} (t) = -  \int_{\p \mathcal{S}} y \wedge \sigma(\hat{U},\hat{P}) n \d \Gamma + \mathbb{F}_I , & &  t\in(0,\infty),
\end{eqnarray}
\begin{eqnarray}
\hat{U}(y,0)  =  u_0 (y), \  y \in \mathcal{F}, \quad \hat{H}'(0)=h_1 \in \R^3 ,\quad \hat{\Omega}(0) = \omega_0 \in \R^3, \label{linnhbb9}
\end{eqnarray}
We assume that the data satisfy
\begin{eqnarray*}
\begin{array} {ll}
\mathbb{F} \in \L^2(0,\infty;\mathbf{L}^2(\mathcal{F})), & \\
\mathbb{G} \in \H^{2,1}(Q_{\infty}^0), &
\mathbb{W} \in \L^2(0,\infty;\mathbf{H}^{3/2}(\p \mathcal{S})) \cap \H^1(0,\infty;\mathbf{H}^{1/2}(\p \mathcal{S})), \\
\mathbb{F}_M \in \L^2(0,\infty;\R^3), & \mathbb{F}_I \in \L^2(0,\infty;\R^3),
\end{array}
\end{eqnarray*}
and the following compatibility conditions
\begin{eqnarray*}
\int_{\p \mathcal{S}} \zeta \cdot n \d \Gamma = 0,  \quad
\int_{\p \mathcal{S}} \mathbb{G} \cdot n \d \Gamma = \int_{\p \mathcal{S}} \mathbb{W} \cdot n \d \Gamma, \quad
\mathbb{G}(\cdot,0) = 0, \quad \mathbb{G}_{| \p \mathcal{O}} = 0.
\end{eqnarray*}

\begin{remark}
Note that in this system the control $\zeta$ does not stand for $ e^{\lambda t}\frac{\p X^{\ast}}{\p t}$, but it represents a boundary control which satisfies the linearized version of the constraints \eqref{const1}--\eqref{const3}. That is why this notation $\zeta$ is only used for a control for a linear system, like in Part I. Thus the right-hand-side term $\mathbb{W}$ in \eqref{linnhbb4} stands for
\begin{eqnarray*}
W(\hat{\omega}) + e^{\lambda t}\frac{\p X^{\ast}}{\p t} - \zeta.
\end{eqnarray*}
\end{remark}

\subsection{Operator formulation}
The nonhomogeneous divergence condition \eqref{linnhbb2} can be lifted by setting
\begin{eqnarray*}
\overline{U} & = & \hat{U} - \mathbb{G}.
\end{eqnarray*}
Then the quadruplet $(\overline{U}, \hat{P}, \hat{H}', \hat{\Omega})$ is supposed to satisfy the following system
\begin{eqnarray*}
\frac{\p \overline{U}}{\p t} - \lambda \hat{U} - \nu \Delta \hat{U} + \nabla \hat{P}  =  \mathbb{F} + \mathbb{F}_{\mathbb{G}}
, &  & \textrm{in $\mathcal{F}\times (0,\infty)$},  \label{linnhbbb1}\\
\div \  \overline{U}  =  0, &  & \textrm{in $\mathcal{F}\times (0,\infty)$}, \label{linnhbbb2}
\end{eqnarray*}
\begin{eqnarray}
\overline{U}  =  0 , & & \textrm{in $\p \mathcal{O}\times (0,\infty)$}, \nonumber \\
\overline{U}=  \hat{H}'(t) + \hat{\Omega}(t) \wedge y + \zeta(y,t) + \mathbb{W}(y,t)-\mathbb{G}(y,t), & & \textrm{in } \p \mathcal{S}\times(0,\infty), \label{linnhbbb4}
\end{eqnarray}
\begin{eqnarray*}
M \hat{H}''(t)-\lambda M \hat{H}'(t) = - \int_{\p \mathcal{S}} \sigma(\overline{U},\hat{P}) n \d \Gamma + \mathbb{F}_M + \mathbb{F}_{M,\mathbb{G}},
& &  t\in (0,\infty), \\
I_0 \hat{\Omega}'(t) - \lambda I_0\hat{\Omega} (t) = -  \int_{\p \mathcal{S}} y \wedge \sigma(\overline{U},\hat{P}) n \d \Gamma + \mathbb{F}_I + \mathbb{F}_{I, \mathbb{G}}  , & &   t\in(0,\infty),
\end{eqnarray*}
\begin{eqnarray*}
\overline{U}(y,0)  =  u_0 (y), \  y \in \mathcal{F}, \quad \hat{H}'(0)=h_1 \in \R^3 ,\quad \hat{\Omega}(0) = \omega_0 \in \R^3, \label{linnhbbb9}
\end{eqnarray*}
where
\begin{eqnarray*}
\mathbb{F}_{\mathbb{G}} & = & -\frac{\p \mathbb{G}}{\p t} + \lambda \mathbb{G} + \nu \Delta \mathbb{G}, \\
\mathbb{F}_{M, \mathbb{G}} & = & -2\nu \int_{\p \mathcal{S}} D(\mathbb{G})n\d \Gamma, \\
\mathbb{F}_{I, \mathbb{G}} & = & -2\nu \int_{\p \mathcal{S}} y\wedge D(\mathbb{G})n\d \Gamma.
\end{eqnarray*}
By following the steps of the operator formulation used in Part I for the homogeneous linear system, the right-hand-side $\mathbb{W} - \mathbb{G}$ in \eqref{linnhbbb4} can be lifted and the pressure of the resulting system can be eliminated, so that the latter can be formally rewritten as follows
\begin{eqnarray}
\frac{\d \overline{\mathbf{U}}}{\d t} & = & \mathcal{A}_{\lambda}\overline{\mathbf{U}} + \mathcal{B}_{\lambda}\zeta
+ \mathbf{F} + \mathcal{B}_{\lambda} \left(\mathbb{W} - \mathbb{G}\right) \label{eqevol1} \\
\left( \Id -\mathbb{P} \right) \overline{U} & = & \left( \Id -\mathbb{P} \right)\left(L_0(\hat{H}') + \hat{L}_0(\hat{\Omega})\right) \nonumber
\end{eqnarray}
with $\overline{\mathbf{U}} = \left(\mathbb{P}\overline{U}, \hat{H}', \hat{\Omega} \right)^T$ and $\mathbf{F} = \displaystyle \left(\mathbb{F} + \mathbb{F}_{\mathbb{G}} , \mathbb{F}_M + \mathbb{F}_{M, \mathbb{G}} , \mathbb{F}_I + \mathbb{F}_{I, \mathbb{G}} \right)$.
In this formulation the operators $\mathbb{P}$, $L_0$, $\mathcal{A}_{\lambda}$ and $\mathcal{B}_{\lambda}$ are given in Part I (see section 3).\\

Let us now remind the most important result established in Part I of this work.

\subsection{Boundary feedback stabilization}
We replace the control $\zeta$ by the feedback operator denoted by $\mathcal{K}_{\lambda}$ and defined in section 5 of Part I. Then in the evolution equation \eqref{eqevol1} the operator $\mathcal{A}_{\lambda}$ becomes $\mathcal{A}_{\lambda} + \mathcal{B}_{\lambda} \mathcal{K}_{\lambda}$. This latter is stable, so that, \textcolor{black}{by following \cite{Bensoussan} (see Theorem 3.1 of Chapter 1) for instance}, we can estimate
\begin{eqnarray*}
& & \|\overline{U} \|_{\H^{2,1}(Q_{\infty}^0)} + \| \hat{P} \|_{\L^2(0,\infty;\H^1(\mathcal{F}))} +
\| H' \|_{\H^1(0,\infty;\R^3)} + \| \hat{\Omega} \|_{\H^1(0,\infty;\R^3)}\\
& & \leq C\left( \|u_0\|_{\mathbf{H}^1(\mathcal{F})} + | h_1 |_{\R^3} + |\omega_0 |_{\R^3}
 + \| \mathbb{F}\|_{\L^2(0,\infty;\mathbf{L}^2(\mathcal{F}))} + \| \mathbb{G} \|_{\H^{2,1}(Q_{\infty}^0)} \right. \\
& & \left. +
\| \mathbb{W}\|_{\L^2(0,\infty;\mathbf{H}^{3/2}(\p \mathcal{F}))\cap \H^1(0,\infty;\mathbf{H}^{1/2}(\p \mathcal{F}))}
+ \| \mathbb{F}_M\|_{\L^2(0,\infty;\R^3)} + \| \mathbb{F}_I\|_{\L^2(0,\infty;\R^3)} \right).
\end{eqnarray*}
Let us keep in mind that we have $\hat{U} = \overline{U} + \mathbb{G}$, and so we have also the following estimate
\begin{eqnarray}
& & \|\hat{U} \|_{\H^{2,1}(Q_{\infty}^0)} + \| \hat{P} \|_{\L^2(0,\infty;\H^1(\mathcal{F}))} +
\| H' \|_{\H^1(0,\infty;\R^3)} + \| \hat{\Omega} \|_{\H^1(0,\infty;\R^3)} \nonumber \\
& & \leq C\left( \|u_0\|_{\mathbf{H}^1(\mathcal{F})} + | h_1 |_{\R^3} + |\omega_0 |_{\R^3}
 + \| \mathbb{F}\|_{\L^2(0,\infty;\mathbf{L}^2(\mathcal{F}))} + \| \mathbb{G} \|_{\H^{2,1}(Q_{\infty}^0)} \right. \nonumber \\
& & \left. +
\| \mathbb{W}\|_{\L^2(0,\infty;\mathbf{H}^{3/2}(\p \mathcal{F}))\cap \H^1(0,\infty;\mathbf{H}^{1/2}(\p \mathcal{F}))}
+ \| \mathbb{F}_M\|_{\L^2(0,\infty;\R^3)} + \| \mathbb{F}_I\|_{\L^2(0,\infty;\R^3)} \right). \nonumber \\ \label{superestK}
\end{eqnarray}
for some independent constant $C >0$.

\section{Construction of a stabilizing admissible deformation} \label{secdecompcontrol}
Let us consider a boundary stabilizing control $\zeta \in \L^2(0,\infty;\mathbf{H}^{5/2}(\p \mathcal{S})) \cap \H^1(0,\infty;\mathbf{H}^{1/2}(\p \mathcal{S}))$ which can be chosen in a feedback form, as described above. The purpose of this section consists in defining from this boundary function a deformation $X^{\ast}$ which is {\it admissible} in the sense of Definition \ref{defcontrol}, and which has to have a satisfying Lipschitz behavior with respect to the function $\zeta$.

\subsection{statement}
The main result of this section is the following:

\begin{theorem} \label{thdecompsuper}
Let be $\zeta \in \L^2(0,\infty;\mathbf{H}^{5/2}(\p \mathcal{S})) \cap \H^1(0,\infty;\mathbf{H}^{1/2}(\p \mathcal{S}))$ satisfying
\begin{eqnarray*}
\int_{\p \mathcal{S}} \zeta \cdot n \d \Gamma & = & 0.
\end{eqnarray*}
If $\zeta$ is small enough in $\L^2(0,\infty;\mathbf{H}^{5/2}(\p \mathcal{S})) \cap \H^1(0,\infty;\mathbf{H}^{1/2}(\p \mathcal{S}))$, then there exists a mapping $X^{\ast}$ which is {\it admissible} in the sense of Definition \ref{defcontrol}, and which satisfies
\begin{eqnarray}
\| X^{\ast} - \Id_{\mathcal{S}} \|_{\mathcal{W}_{\lambda}(S_{\infty}^0)} & \leq &
\| \zeta \|_{\L^2(\mathbf{H}^{5/2}(\p \mathcal{S})) \cap \H^1(\mathbf{H}^{1/2}(\p \mathcal{S}))} , \label{estsuper} \label{petitobis2} \label{estdd12} \\
\left\| e^{\lambda t}\frac{\p X^{\ast}}{\p t} - \zeta \right\|_{\L^2(\mathbf{H}^{5/2}(\p \mathcal{S}))\cap \H^1(\mathbf{H}^{1/2}(\p \mathcal{S}))}
& = & o\left(\| \zeta \|_{\L^2(\mathbf{H}^{5/2}(\p \mathcal{S})) \cap \H^1(\mathbf{H}^{1/2}(\p \mathcal{S}))} \right). \nonumber \\ \label{estsuper0}
\end{eqnarray}
Moreover, if two functions $\zeta_1$ and $\zeta_2$ are close enough to $0$ in $\L^2(0,\infty;\mathbf{H}^{5/2}(\p \mathcal{S})) \cap \H^1(0,\infty;\mathbf{H}^{1/2}(\p \mathcal{S}))$, then the {\it admissible} deformations $X^{\ast}_1$ and $X^{\ast}_2$ that they define respectively satisfy
\begin{eqnarray}
\| X^{\ast}_2 - X^{\ast}_1 \|_{\mathcal{W}_{\lambda}(S_{\infty}^0)} & \leq &
\| \zeta_2 - \zeta_1 \|_{\L^2(\mathbf{H}^{5/2}(\p \mathcal{S})) \cap \H^1(\mathbf{H}^{1/2}(\p \mathcal{S}))},
\label{petitoter2} \label{estsuper12}
\end{eqnarray}
\begin{eqnarray}
& & \left\| \left(e^{\lambda t}\frac{\p X^{\ast}_2}{\p t} - \zeta_2 \right) - \left(e^{\lambda t}\frac{\p X^{\ast}_1}{\p t} - \zeta_1 \right) \right\|_{\L^2(\mathbf{H}^{5/2}(\p \mathcal{S}))\cap \H^1(\mathbf{H}^{1/2}(\p \mathcal{S}))} \nonumber \\
& & \leq \| \zeta_2 - \zeta_1 \|_{\L^2(\mathbf{H}^{5/2}(\p \mathcal{S})) \cap \H^1(\mathbf{H}^{1/2}(\p \mathcal{S}))} \nonumber \\
& & \quad \times K\left(\| \zeta_1 \|_{\L^2(\mathbf{H}^{5/2}(\p \mathcal{S})) \cap \H^1(\mathbf{H}^{1/2}(\p \mathcal{S}))}
+ \| \zeta_2 \|_{\L^2(\mathbf{H}^{5/2}(\p \mathcal{S})) \cap \H^1(\mathbf{H}^{1/2}(\p \mathcal{S}))} \right),
 \label{estsuper120}
\end{eqnarray}
with $K(r) \rightarrow 0$ when $r$ goes to $0$.
\end{theorem}

The proof of this theorem is divided into two steps. In the first one, \textcolor{black}{we use the results of section 6 of Part I in order to} construct from a boundary function $\zeta$ an internal solid deformation $X^{\ast}_{\zeta}$ satisfying the linearized versions of constraints \eqref{const1}--\eqref{const3}. Then in a second time we project the displacement $X^{\ast}_{\zeta} - \Id_{\mathcal{S}}$ on a space of displacements which define {\it admissible} deformations.

\subsection{First step of the proof}
The first step of the proof of Theorem \ref{thdecompsuper} consists in defining from a boundary control $\zeta$ an internal deformation $X^{\ast}_{\zeta}$ which satisfies the linearized versions of constraints \eqref{const1}--\eqref{const3}. For that, let us remind a result of Part I, which is a consequence of the addition of Proposition 5 and Corollary 1 of section 6:

\begin{proposition} \label{lemmaxistence2}
For $\zeta \in \L^2(0,\infty;\mathbf{H}^{5/2}(\p \mathcal{S})) \cap \H^1(0,\infty;\mathbf{H}^{1/2}(\p \mathcal{S}))$ satisfying
\begin{eqnarray*}
\int_{\p \mathcal{S}} \zeta \cdot n \d \Gamma & = & 0,
\end{eqnarray*}
the following system
\begin{eqnarray}
\mu \varphi - 2\div \ D(\varphi) = F(\varphi) & & \text{in } \mathcal{S} \label{Lame012} \\
\varphi =  \zeta & & \text{on } \p \mathcal{S} \label{Lame042}
\end{eqnarray}
admits a unique solution $\varphi$ in $\L^2(0,\infty;\mathbf{H}^3(\mathcal{S}) \cap \H^1(0,\infty;\mathbf{H}^1(\mathcal{S}))$, for $\mu>0$ large enough. This function $\varphi$ satisfies the conditions
\begin{eqnarray*}
\int_{\p \mathcal{S}} \varphi \cdot n \d \Gamma = 0, \quad \int_{\mathcal{S}} \varphi \d y = 0 , \quad
\int_{\mathcal{S}} y \wedge \varphi \d y = 0,
\end{eqnarray*}
and there exists a positive constant $C > 0$ such that
\begin{eqnarray}
\| \varphi \|_{\L^2(0,\infty;\mathbf{H}^3(\mathcal{S}) \cap \H^1(0,\infty;\mathbf{H}^1(\mathcal{S}))}
& \leq & C \| \zeta \|_{\L^2(0,\infty;\mathbf{H}^{5/2}(\p \mathcal{S})) \cap \H^1(0,\infty;\mathbf{H}^{1/2}(\p \mathcal{S}))}.
\label{estinftylame2}
\end{eqnarray}
Besides, if $\zeta_1, \zeta_2 \in \L^2(0,\infty;\mathbf{H}^{5/2}(\p \mathcal{S})) \cap \H^1(0,\infty;\mathbf{H}^{1/2}(\p \mathcal{S}))$, then the solutions $\varphi_1$ and $\varphi_2$ associated with $\zeta_1$ and $\zeta_2$ respectively satisfy
\begin{eqnarray}
\| \varphi_2 - \varphi_1 \|_{\L^2(0,\infty;\mathbf{H}^3(\mathcal{S}) \cap \H^1(0,\infty;\mathbf{H}^1(\mathcal{S}))} & \leq &
C \| \zeta_2 - \zeta_1 \|_{\L^2(0,\infty;\mathbf{H}^{5/2}(\p \mathcal{S})) \cap \H^1(0,\infty;\mathbf{H}^{1/2}(\p \mathcal{S}))}. \nonumber \\
\label{estinftylame122}
\end{eqnarray}
\end{proposition}

For $\varphi \in \L^2(0,\infty;\mathbf{H}^3(\mathcal{S}) \cap \H^1(0,\infty;\mathbf{H}^1(\mathcal{S}))$ obtained from $\zeta$, we now define $X^{\ast}_{\zeta}- \Id_{\mathcal{S}} \in \mathcal{W}_{\lambda}(S_{\infty}^0)$ as follows
\begin{eqnarray*}
X^{\ast}_{\zeta}(\cdot,t) = y + \int_0^t e^{-\lambda s}\varphi(\cdot,s) \d s & \Leftrightarrow &
\left\{ \begin{array} {l}
\displaystyle e^{\lambda t}\frac{\p X^{\ast}_{\zeta}}{\p t} = \varphi \\
\displaystyle X^{\ast}_{\zeta}(\cdot,0) = \Id_{\mathcal{S}}.
\end{array} \right.
\end{eqnarray*}
Thus the estimates \eqref{estinftylame2} and \eqref{estinftylame122} become
\begin{eqnarray}
\| X^{\ast}_{\zeta}- \Id_{\mathcal{S}} \| & \leq & C \| \zeta \|_{\L^2(0,\infty;\mathbf{H}^{5/2}(\p \mathcal{S})) \cap \H^1(0,\infty;\mathbf{H}^{1/2}(\p \mathcal{S}))}, \label{estLAME} \\
\| X^{\ast}_{\zeta_2}- X^{\ast}_{\zeta_1} \| & \leq & C \| \zeta_2 - \zeta_1 \|_{\L^2(0,\infty;\mathbf{H}^{5/2}(\p \mathcal{S})) \cap \H^1(0,\infty;\mathbf{H}^{1/2}(\p \mathcal{S}))}. \label{estLAME12}
\end{eqnarray}
The result of this proposition could be actually reduced to saying that there exists a linear continuous operator from $\L^2(0,\infty;\mathbf{H}^{5/2}(\p \mathcal{S})) \cap \H^1(0,\infty;\mathbf{H}^{1/2}(\p \mathcal{S}))$ to $\mathcal{W}_{\lambda}(S_{\infty}^0)$, but the means that we use for obtaining it (that is to say by considering system \eqref{Lame012}-\eqref{Lame042}) is actually useful for a key point in the second part of the proof.

\subsection{Second step of the proof}
Let us consider a control $X^{\ast}_{\zeta} \in \mathcal{W}_{\lambda}(S_{\infty}^0)$ which has been obtained in the previous subsection from a boundary velocity $\zeta$. Instead of projecting $X^{\ast}_{\zeta}$ on a set of controls satisfying the nonlinear constraints required by Definition \ref{defcontrol}, we prefer projecting the displacements $X^{\ast}_{\zeta} - \Id_{\mathcal{S}}$, because we choose the space $\mathcal{W}_{\lambda}(S_{\infty}^0)$ as an Hilbertian framework. We denote the displacement by
\begin{eqnarray*}
Z^{\ast}_{\zeta} & = & X^{\ast}_{\zeta} - \Id_{\mathcal{S}}.
\end{eqnarray*}
The goal of this subsection is to define (in a suitable way) a mapping $X^{\ast} \in \mathcal{W}_{\lambda}(S_{\infty}^0)$ which satisfies the nonlinear constraints of Definition \ref{defcontrol}. We associate with it the displacement
\begin{eqnarray*}
Z^{\ast} & = & X^{\ast} - \Id_{\mathcal{S}},
\end{eqnarray*}
so that the wanted mapping is now $Z^{\ast}$. We can decompose such a mapping as follows
\begin{eqnarray*}
Z^{\ast}  =  Z^{\ast}_{\zeta} + \left(Z^{\ast} - Z^{\ast}_{\zeta} \right) =  Z^{\ast}_{\zeta} + \left(X^{\ast} - X^{\ast}_{\zeta} \right).
\end{eqnarray*}
Let us define the differentiable mapping
\begin{eqnarray*}
\begin{array} {crcl}
 \mathfrak{F} \ : & \mathcal{W}_{\lambda}(S_{\infty}^0) & \rightarrow &
 \H^1_0(0,\infty;\R^{3}) \times \H^1_0(0,\infty;\R^{3}) \times \H^1_0(0,\infty;\R)  \\
& Z^{\ast} & \mapsto & \displaystyle \mathfrak{F}(Z^{\ast})
\end{array} \\
\mathfrak{F}(Z^{\ast})(y,t)  =  \left( \displaystyle \int_{\mathcal{S}} \frac{\p Z^{\ast}}{\p t};
\ \displaystyle \int_{\mathcal{S}} \left(Z^{\ast}+\Id_{\mathcal{S}} \right)\wedge \frac{\p Z^{\ast}}{\p t};
\ \displaystyle \ \int_{\p \mathcal{S}} (\com (\nabla Z^{\ast}+ \I_{\R^3}))^T\frac{\p Z^{\ast}}{\p t}\cdot n \right),
\end{eqnarray*}
and the spaces
\begin{eqnarray*}
\mathcal{E}_{nl} & = & \left\{Z^{\ast} \in \mathcal{W}_{\lambda}(S_{\infty}^0) \mid \mathfrak{F}(Z^{\ast})=0 \right\},\\
\mathcal{E}_{l} & = & \left\{Z^{\ast}_{\zeta} \in \mathcal{W}_{\lambda}(S_{\infty}^0) \mid D_{0}\mathfrak{F}(Z^{\ast}_{\zeta})=0 \right\},
\end{eqnarray*}
where
\begin{eqnarray*}
D_{0}\mathfrak{F}(Z^{\ast}_{\zeta}) & = & t\mapsto \left(\int_{\mathcal{S}} \frac{\p Z^{\ast}_{\zeta}}{\p t};
 \ \int_{\mathcal{S}} \Id_{\mathcal{S}}\wedge \frac{\p Z^{\ast}_{\zeta}}{\p t};
 \ \int_{\p \mathcal{S}} \frac{\p Z^{\ast}_{\zeta}}{\p t}\cdot n \right).
\end{eqnarray*}
Note that $\mathcal{E}_{l}$ is a space where lie $\Id_{\mathcal{S}}$ and $X^{\ast}_{\zeta}$. That is why the constraints satisfied by $Z^{\ast}_{\zeta}$ and $X^{\ast}_{\zeta}$ are the same. Note that the space $\mathcal{E}_{nl}$ takes into account the nonlinear constraints adapted to the displacements $Z^{\ast}$. The purpose of this paragraph is to project any displacement $Z^{\ast}_{\zeta} \in \mathcal{E}_{l}$ on the set $\mathcal{E}_{nl}$, provided that the displacement $Z^{\ast}_{\zeta}$ is close enough to $0$.
\begin{figure}[!h]
\begin{center}
\hspace*{-1.2cm}\scalebox{0.8}
{ \input{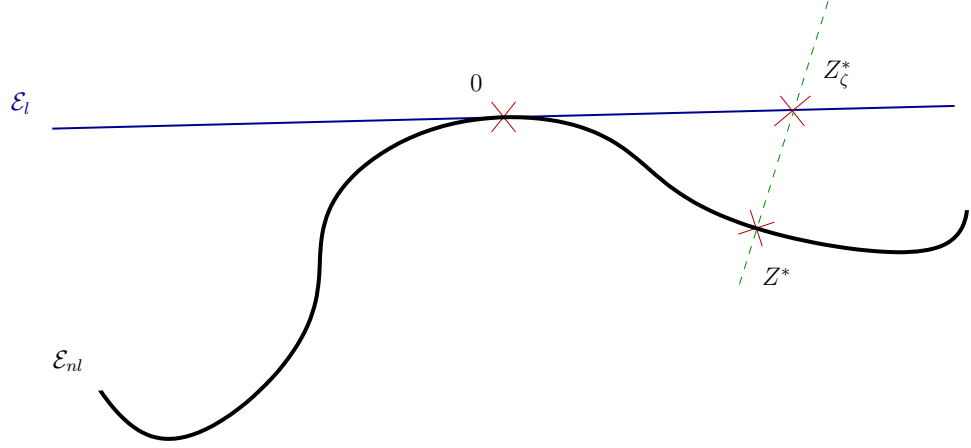} }
\caption{Projection of a given mapping on the set of admissible displacements}\label{projection}
\end{center}
\end{figure}
\FloatBarrier

The definition of such a projection is given by:
\begin{theorem} \label{thdecompsuperbis}
Let be $Z^{\ast}_{\zeta} \in \mathcal{W}_{\lambda}(S_{\infty}^0)$. If $Z^{\ast}_{\zeta}$ is small enough in $\mathcal{W}_{\lambda}(S_{\infty}^0)$, then there exists a unique mapping $Z^{\ast} \in \mathcal{E}_{nl}$ such that
\begin{eqnarray*}
\| Z^{\ast} - Z^{\ast}_{\zeta} \|^2_{\mathcal{W}_{\lambda}(S_{\infty}^0)} & = &
\min_{\mathcal{Z}^{\ast} \in \mathcal{E}_{nl}} \| \mathcal{Z}^{\ast} - X^{\ast}_{\zeta} \|^2_{\mathcal{W}_{\lambda}(S_{\infty}^0)}.
\end{eqnarray*}
Moreover, we have that
\begin{eqnarray}
\| Z^{\ast} - Z^{\ast}_{\zeta} \|_{\mathcal{W}_{\lambda}(S_{\infty}^0)} & = &
o \left(\| Z^{\ast}_{\zeta} \|_{\mathcal{W}_{\lambda}(S_{\infty}^0)}\right). \label{petitobis} \label{estdd1}
\end{eqnarray}
Thus we denote by $\mathcal{P} : Z^{\ast}_{\zeta} \mapsto Z^{\ast}$ the projection so obtained.\\
If the displacements $Z^{\ast}_{\zeta_1}$ and $Z^{\ast}_{\zeta_2}$ are close enough to $0$ in $\mathcal{W}_{\lambda}(S_{\infty}^0)$, then
\begin{eqnarray}
& & \|(Z^{\ast}_2 - Z^{\ast}_{\zeta_2}) - (Z^{\ast}_1 - Z^{\ast}_{\zeta_1})\|_{\mathcal{W}_{\lambda}(S_{\infty}^0)}  \leq
\|Z^{\ast}_{\zeta_2} - Z^{\ast}_{\zeta_1} \|_{\mathcal{W}_{\lambda}(S_{\infty}^0)} \nonumber \\
& & \qquad \times K\left(\| Z^{\ast}_{\zeta_1} \|_{\mathcal{W}_{\lambda}(S_{\infty}^0)}
+ \| Z^{\ast}_{\zeta_2} \|_{\mathcal{W}_{\lambda}(S_{\infty}^0)} \right), \label{petitoter}
\end{eqnarray}
with $Z^{\ast}_1 = \mathcal{P}Z^{\ast}_{\zeta_1}$ and $Z^{\ast}_2 = \mathcal{P}Z^{\ast}_{\zeta_2}$, and $K(r) \rightarrow 0$ when $r$ goes to $0$.
\end{theorem}

\begin{remark}
Note that in this statement we do not need to assume that $Z^{\ast}_{\zeta} \in \mathcal{E}_l$. But the way we have constructed $X^{\ast}_{\zeta}$ such that $Z^{\ast}_{\zeta} = X^{\ast}_{\zeta}-\Id_{\mathcal{S}} \in \mathcal{E}_l$ in the previous subsection will be useful in the proof below.
\end{remark}

\begin{proof}
The proof of this theorem is an application of Theorem 3.33 of \cite{Bonnans} (page 74), that we state as follows:
\begin{theorem}
Let $\mathcal{W}$ be a Hilbert space, $\mathcal{R}$ a Banach space, and $g$ a mapping of class $C^2$ from $\mathcal{W}$ to $\mathcal{R}$, such that $g^{-1}(\{0\}) \neq \emptyset$. Let be $Z^{\ast}_{\zeta} \in \mathcal{W}$, and $Z_0^{\ast} \in g^{-1}(\{0\})$ such that $D_{X_0^{\ast}}g$ is surjective. Then there exists  $\varepsilon >0$ such that if $\|Z^{\ast}_{\zeta} - Z^{\ast}_0 \|_{\mathcal{W}} \leq \varepsilon $, then the following optimization problem under equality constraints
\begin{eqnarray*}
\min_{\mathcal{Z}^{\ast} \in \mathcal{W}, \ g(\mathcal{Z}^{\ast}) = 0} \| \mathcal{Z}^{\ast} - Z^{\ast}_{\zeta} \|^2_{\mathcal{W}}
\end{eqnarray*}
admits a unique solution $Z^{\ast}$.
Moreover, the mapping $Z^{\ast}_{\zeta} \mapsto Z^{\ast}$ so obtained is $C^1$.
\end{theorem}
In order to apply this theorem with
\begin{eqnarray*}
\mathcal{W} = \mathcal{W}_{\lambda}(S_{\infty}^0), & \quad &
\mathcal{R} = \H^1_0(0,\infty;\R^{3}) \times \H^1_0(0,\infty;\R^{3}) \times \H^1_0(0,\infty;\R), \\
g = \mathfrak{F}, & \quad & Z_0^{\ast} = 0,
\end{eqnarray*}
the only nontrivial assumption to be verified is that the mapping $D_{0}\mathfrak{F}$ is surjective. For that, let us consider
\begin{eqnarray*}
(a,b,c) & \in  & \H^1_0(0,\infty;\R^{3}) \times \H^1_0(0,\infty;\R^{3}) \times \H^1_0(0,\infty;\R).
\end{eqnarray*}
An antecedent $\mathcal{Z}^{\ast}$ of this triplet $(a,b,c)$ can be obtained as
\begin{eqnarray*}
\mathcal{Z}^{\ast}(y,t) & = & \int_0^t e^{-\lambda s} \varphi(y,s) \d s,
\end{eqnarray*}
where $\varphi$ is the solution of the following system
\begin{eqnarray*}
\mu \varphi - 2\div \ D(\varphi) = F(\varphi)+F_{a,b} & \quad & \text{in } \mathcal{S}, \\
\varphi =  \frac{c}{|\p \mathcal{S} |}n & \quad & \text{on } \p \mathcal{S},
\end{eqnarray*}
for $\mu >0$ large enough, with
\begin{eqnarray*}
F_{a,b}(y,t) & = & \frac{\rho_{\mathcal{S}}}{M}a(t) + \rho_{\mathcal{S}} \left(I_0^{-1} b(t)\right) \wedge y, \\
F(\mathcal{Z}^{\ast})(y,t) & = & \frac{\rho_{\mathcal{S}}}{M}\left(\int_{\p \mathcal{S}}2D(\mathcal{Z}^{\ast})n\d \Gamma\right)
 + \left(I_0^{-1}\int_{\p \mathcal{S}} \Id_{\mathcal{S}} \wedge 2D(\mathcal{Z}^{\ast})n\d \Gamma \right) \wedge y.
\end{eqnarray*}
The previous study of system \eqref{Lame012}-\eqref{Lame042} (see section 6 of Part I) can be straightforwardly adapted to get the existence of a solution in $\varphi \in \L^2(0,\infty;\mathbf{H}^3(\mathcal{S}))\cap \H^1(0,\infty;\mathbf{H}^1(\mathcal{S}))$ for such a system, and thus a displacement $\mathcal{Z}^{\ast} \in \mathcal{W}_{\lambda}(S_{\infty}^0)$.\\
Since the projection $\mathcal{P} : Z^{\ast}_{\zeta} \mapsto Z^{\ast}$ so obtained is $C^1$, we can notice that its differential at $0$ is the identity $\Id_{\mathcal{W}_{\lambda}(S_{\infty}^0)}$, and thus a Taylor development shows that
\begin{eqnarray*}
& & Z^{\ast}  =   Z^{\ast}_{\zeta} + 0 + o\left(\|Z^{\ast}_{\zeta}\|_{\mathcal{W}_{\lambda}(S_{\infty}^0)}\right), \\
& & \|Z^{\ast} - Z^{\ast}_{\zeta} \|_{\mathcal{W}_{\lambda}(S_{\infty}^0)} = o\left(\|Z^{\ast}_{\zeta}\|_{\mathcal{W}_{\lambda}(S_{\infty}^0)}\right).
\end{eqnarray*}
For $Z^{\ast}_{\zeta_1}$ and $Z^{\ast}_{\zeta_2}$ close to $0$, the estimate \eqref{petitoter} is obtained in considering a Taylor development around $Z^{\ast}_{\zeta_1}$ for the mapping $\mathcal{P} - \Id_{\mathcal{W}_{\lambda}(S_{\infty}^0)}$:
\begin{eqnarray*}
\left(Z^{\ast}_2 - Z^{\ast}_{\zeta_2}\right) - \left(Z^{\ast}_1 - Z^{\ast}_{\zeta_1}\right) & = &
\left[D_{Z^{\ast}_{\zeta_1}}\mathcal{P} - \Id_{\mathcal{W}_{\lambda}(S_{\infty}^0)} \right]\left(Z^{\ast}_{\zeta_2} - Z^{\ast}_{\zeta_1} \right) \\
&  & + o\left(\|Z^{\ast}_{\zeta_2} - Z^{\ast}_{\zeta_1} \|_{\mathcal{W}_{\lambda}(S_{\infty}^0)} \right).
\end{eqnarray*}
Since $\mathcal{Z}^{\ast} \mapsto D_{\mathcal{Z}^{\ast}}\mathcal{P}$ is continuous at $0$, we have
\begin{eqnarray*}
D_{Z^{\ast}_{\zeta_1}}\mathcal{P} - \Id_{\mathcal{W}_{\lambda}(S_{\infty}^0)} \rightarrow 0 & \quad & \text{when }
\| Z^{\ast}_{\zeta_1} \|_{\mathcal{W}_{\lambda}(S_{\infty}^0)} \text{ goes to $0$,}
\end{eqnarray*}
and thus we obtain the announced estimate.
\end{proof}

Then, from the displacement $Z^{\ast}_{\zeta} = X^{\ast}_{\zeta} - \Id_{\mathcal{S}}$ we can define a deformation $X^{\ast}$ as follows
\begin{eqnarray*}
X^{\ast} & = & \mathcal{P}(X^{\ast}_{\zeta}-\Id_{\mathcal{S}}) + \Id_{\mathcal{S}}.
\end{eqnarray*}
This deformation is {\it admissible} in the sense of Definition \ref{defcontrol}.\\
The interest of such a decomposition (namely $X^{\ast} = X^{\ast}_{\zeta} + (X^{\ast}-X^{\ast}_{\zeta})$ with $X^{\ast}$ given by Theorem \ref{thdecompsuper}) lies in the fact that the {\it admissible} control $X^{\ast}$ so decomposed will enable us to stabilize the nonhomogeneous linear part of system \eqref{prems}--\eqref{ders} thanks to the term $X^{\ast}_{\zeta}$ (see the previous section), whereas the residual term $(X^{\ast}-X^{\ast}_{\zeta})$ satisfies the property \eqref{petitobis}, which leads to
\begin{eqnarray*}
\|X^{\ast} - X^{\ast}_{\zeta} \|_{\mathcal{W}_{\lambda}(S_{\infty}^0)} & = &
o \left( \| X^{\ast}_{\zeta}- \Id_{\mathcal{S}} \|_{\mathcal{W}_{\lambda}(S_{\infty}^0)} \right) , \label{petitoo} \\
\|X^{\ast} - \Id_{\mathcal{S}} \|_{\mathcal{W}_{\lambda}(S_{\infty}^0)} & \leq &
C\| X^{\ast}_{\zeta}- \Id_{\mathcal{S}} \|_{\mathcal{W}_{\lambda}(S_{\infty}^0)} ,  \\
\left\| e^{\lambda t} \frac{\p X^{\ast}}{\p t} - e^{\lambda t} \frac{\p X^{\ast}_{\zeta}}{\p t} \right\|
_{\L^2(\mathbf{H}^{5/2}(\p \mathcal{S}))\cap \H^1(\mathbf{H}^{1/2}(\p \mathcal{S}))} & = &
o \left( \| X^{\ast}_{\zeta}- \Id_{\mathcal{S}} \|_{\mathcal{W}_{\lambda}(S_{\infty}^0)} \right). \label{petito}
\end{eqnarray*}
By combining the second and the third inequality to the estimate \eqref{estLAME}, we then obtain respectively \eqref{estsuper} and \eqref{estsuper0}.\\
Lastly, the estimate \eqref{petitoter} is reformulated as follows
\begin{eqnarray}
& & \|(X^{\ast}_2 - X^{\ast}_{\zeta_2}) - ( X^{\ast}_1 - X^{\ast}_{\zeta_1})\|_{\mathcal{W}_{\lambda}(S_{\infty}^0)}  \leq \nonumber \\
& & K\left(\| X^{\ast}_{\zeta_1}- \Id_{\mathcal{S}} \|_{\mathcal{W}_{\lambda}(S_{\infty}^0)}
+ \| X^{\ast}_{\zeta_2}- \Id_{\mathcal{S}} \|_{\mathcal{W}_{\lambda}(S_{\infty}^0)} \right)
\times \| X^{\ast}_{\zeta_2}- X^{\ast}_{\zeta_1} \|_{\mathcal{W}_{\lambda}(S_{\infty}^0)},
 \label{petitooo12} \\
& & \left\| \left(e^{\lambda t}\frac{\p X^{\ast}_2}{\p t} - \zeta_2\right) - \left(e^{\lambda t}\frac{\p X^{\ast}_1}{\p t} - \zeta_1\right) \right\|_{\L^2(0,\infty;\mathbf{H}^{5/2}(\p \mathcal{S}))\cap \H^1(0,\infty;\mathbf{H}^{1/2}(\p \mathcal{S}))}  \leq \nonumber \\
& & K\left(\| X^{\ast}_{\zeta_1}- \Id_{\mathcal{S}} \|_{\mathcal{W}_{\lambda}(S_{\infty}^0)}
+ \| X^{\ast}_{\zeta_2}- \Id_{\mathcal{S}} \|_{\mathcal{W}_{\lambda}(S_{\infty}^0)} \right)
\times \| X^{\ast}_{\zeta_2}- X^{\ast}_{\zeta_1} \|_{\mathcal{W}_{\lambda}(S_{\infty}^0)}, \nonumber \\ \label{petito12}
\end{eqnarray}
where $K(r) \rightarrow 0$ when $r$ goes to $0$, and where
\begin{eqnarray*}
X^{\ast}_1 & = & \mathcal{P}(X^{\ast}_{\zeta_1}-\Id_{\mathcal{S}}) + \Id_{\mathcal{S}}, \\
X^{\ast}_2 & = & \mathcal{P}(X^{\ast}_{\zeta_2}-\Id_{\mathcal{S}}) + \Id_{\mathcal{S}}.
\end{eqnarray*}
The inequality \eqref{petitooo12} combined to \eqref{estLAME} and \eqref{estLAME12} leads to the estimate \eqref{estsuper12}, and the inequality \eqref{petito12} combined to \eqref{estLAME} and \eqref{estLAME12} leads to the estimate \eqref{estsuper120}.

\section{A fixed point method} \label{secnonlinear}
System \eqref{prems}--\eqref{ders} is transformed into system \eqref{hpremsfix}--\eqref{hdersfix}. Before proving the stabilization to zero, with an arbitrary exponential decay rate $\lambda >0$, of system \eqref{prems}--\eqref{ders}, let us first prove the stability of system \eqref{hpremsfix}--\eqref{hdersfix} for all $\lambda >0$, for some well-chosen deformation $X^{\ast}$.

\subsection{Back to the nonlinear system written in a cylindrical domain}
In system \eqref{hpremsfix}--\eqref{hdersfix}, the mapping $X^{\ast}$ has to be {\it admissible} (in the sense of Definition \ref{defcontrol}). It has to be chosen also in order to stabilize the linear part of this system. For that, we decompose formally the function $e^{\lambda t}\frac{\p X^{\ast}}{\p t}$ on $\p \mathcal{S}$ as follows
\begin{eqnarray*}
e^{\lambda t} \frac{\p X^{\ast}}{\p t} & = & \zeta + \left(e^{\lambda t} \frac{\p X^{\ast}}{\p t} - \zeta \right)
\end{eqnarray*}
Let us choose in system \eqref{hpremsfix}--\eqref{hdersfix} the function $\zeta$ in the following feedback form:
\begin{eqnarray*}
\zeta & = & \mathcal{K}_{\lambda}(\mathbb{P}(\hat{u}-G(\hat{u},\hat{h}',\hat{\omega})),\hat{h}',\hat{\omega}).
\end{eqnarray*}
Provided that $\zeta$ is small enough, by Theorem \ref{thdecompsuper} we can now define an {\it admissible} deformation $X^{\ast}$. From this deformation we can define a change of variables $X$ and the corresponding mappings $\tilde{X}$ and $\tilde{Y}$ (see \eqref{definitilde}) which enables us to define the right-hand-sides of system \eqref{hpremsfix}--\eqref{hdersfix}; More precisely, we rewrite this system as follows
\begin{eqnarray}
\frac{\p \hat{u}}{\p t} - \lambda \hat{u} -\nu \Delta \hat{u} + \nabla \hat{p}
 =  F(\hat{u},\hat{p},\hat{h}',\hat{\omega}), & \quad &  \text{in } \mathcal{F}\times (0,\infty), \label{hhpremsfix} \\
\div \ \hat{u}  =  \div \ G(\hat{u},\hat{h}',\hat{\omega}), & \quad & \text{in } \mathcal{F}\times (0,\infty),
\end{eqnarray}
\begin{eqnarray}
\hat{u} = 0 , & \quad & \text{on } \p \mathcal{O}\times (0,\infty),  \\
\hat{u}  =  \hat{h}'(t) + \hat{\omega} (t) \wedge y+ \mathcal{K}_{\lambda}(\mathbb{P}(\hat{u}-G(\hat{u},\hat{h}',\hat{\omega})),\hat{h}',\hat{\omega}) & & \\
+ \left(e^{\lambda t}\frac{\p X^{\ast}}{\p t}-\zeta \right) + W(\hat{u},\hat{h}',\hat{\omega}) , & \quad &
(y,t)\in \p \mathcal{S}\times (0,\infty),
\end{eqnarray}
\begin{eqnarray}
M \hat{h}'' - \lambda M \hat{h}'  =  - \int_{\p \mathcal{S}} \sigma(\hat{u},\hat{p}) n  \d \Gamma + F_M(\hat{u},\hat{p},\hat{h}',\hat{\omega}),
& \quad & \text{in } (0,\infty),   \\
I_0\hat{\omega}' (t) -\lambda I_0\hat{\omega}  =   -  \int_{\p \mathcal{S}} y \wedge \sigma(\hat{u},\hat{p}) n  \d \Gamma + F_I(\hat{u},\hat{p},\hat{\omega}), & \quad  & \text{in } (0,\infty),
\end{eqnarray}
\begin{eqnarray}
\hat{u}(y,0)  =  u_0 (y), \  y\in \mathcal{F} , \quad \hat{h}'(0)=h_1 \in \R^3 ,\quad \hat{\omega}(0) =\omega_0 \in \R^3 , \label{hhdersfix}
\end{eqnarray}
with
\begin{eqnarray}
F(\hat{u},\hat{p},\hat{h}',\hat{\omega}) & = & \nu (\LL - \Delta) \hat{u} - e^{-\lambda t}\MM (\hat{u}, \hat{h}', \hat{\omega}) - e^{-\lambda t}\NN \hat{u} - (\GG-\nabla) \hat{p} - e^{-\lambda t}\hat{\omega} \wedge \hat{u}, \label{rhsF} \\
G(\hat{u},\hat{h}',\hat{\omega}) & = & \left(\I_{\R^3}- \nabla \tilde{Y}(\tilde{X}(y,t),t)\right)\hat{u}, \label{rhsG} \\
W(\hat{u},\hat{h}',\hat{\omega}) & = & \hat{\omega} \wedge \left(X^{\ast} - \Id\right), \label{rhsW} \\
\zeta & = & \mathcal{K}_{\lambda}(\mathbb{P}(\hat{u}-G(\hat{u},\hat{h}',\hat{\omega})),\hat{h}',\hat{\omega}), \\
F_M(\hat{u},\hat{p},\hat{h}',\hat{\omega}) & = & -Me^{-\lambda t} \hat{\omega} \wedge \hat{h}'(t) \nonumber \\
 & & - \nu\int_{\p \mathcal{S}}\left(\nabla \hat{u} \left(\nabla \tilde{Y}(\tilde{X}) - \I_{\R^3}\right) + \left({\nabla \tilde{Y}(\tilde{X})} - \I_{\R^3}\right)^T\nabla \hat{u}^T \right)\nabla \tilde{Y}(\tilde{X})^T n \d \Gamma \nonumber \\
 & & - \int_{\p \mathcal{S}}\sigma(\hat{u},\hat{p})\left(\nabla \tilde{Y}(\tilde{X})-\I_{\R^3}\right)^Tn\d \Gamma, \label{rhsFM} \\
F_I(\hat{u},\hat{p},\hat{\omega}) & = & -\left(I^{\ast} - I_0\right) \hat{\omega}' + \lambda\left(I^{\ast}-I_0\right)\hat{\omega} - {I^{\ast}}'\hat{\omega} +e^{-\lambda t}I^{\ast}\hat{\omega} \wedge \hat{\omega} \nonumber \\
& &  - \nu\int_{\p \mathcal{S}}y\wedge \left(\nabla \hat{u} \left(\nabla \tilde{Y}(\tilde{X}) - \I_{\R^3}\right)
+ ({\nabla \tilde{Y}(\tilde{X})} - \I_{\R^3})^T\nabla \hat{u}^T\right)\nabla \tilde{Y}(\tilde{X})^T n \d \Gamma \nonumber \\
& & - \int_{\p \mathcal{S}}y\wedge \left(\sigma(\hat{u},\hat{p})(\nabla \tilde{Y}(\tilde{X})-\I_{\R^3})n\right)^T\d \Gamma \nonumber \\
& & + \int_{\p \mathcal{S}}\left(X^{\ast}-\Id\right)\wedge \left(\tilde{\sigma}(\hat{u},\hat{p})
\nabla \tilde{Y}(\tilde{X})^T n\right)\d \Gamma. \label{rhsFI}
\end{eqnarray}
Note that the mapping $X^{\ast}$ is entirely defined by $\zeta$ and so by the unknowns $(\hat{u},\hat{h}',\hat{\omega})$, while the mappings $\tilde{X}$ and $\tilde{Y}$ are entirely defined by $X^{\ast}$ and the unknowns $(\hat{h}',\hat{\omega})$, and thus by $(\hat{u},\hat{h}',\hat{\omega})$.

\begin{remark}
Note that the projection method used in order to define an {\it admissible} deformation $X^{\ast}$ from a boundary control $\zeta$ has been made in infinite time horizon, with regards to the functional spaces considered. It implies that the nonlinear system above is noncausal, that means that the control $\zeta$ chosen in a feedback form anticipate {\it a priori} at some time $t$ the behavior of the unknowns for later times. In practice we could define a projection method for increasing times, but the corresponding Lipschitz estimates - in infinite time horizon - obtained above do not hold anymore. Anyway, let us show that this system admits a unique solution, and thus that it makes sense for all time.
\end{remark}

\subsection{Statement}
\begin{theorem} \label{thstabnonlinX}
For $(u_0,h_1,\omega_0)$ small enough in $\mathbf{H}^1_{cc}$, system \eqref{hhpremsfix}--\eqref{hhdersfix} admits a unique solution $(\hat{u},\hat{p},\hat{h}',\hat{\omega})$ in $\H^{2,1}(Q_{\infty}^0) \times \L^2(0,\infty;\H^1(\mathcal{F})) \times \H^1(0,\infty;\R^3) \times \H^1(0,\infty;\R^3)$, and there exists a positive constant $C$ such that
\begin{eqnarray*}
\| \hat{u} \|_{\H^{2,1}(Q_{\infty}^0)} + \| \hat{p} \|_{\L^2(0,\infty;\H^1(\mathcal{F}))} + \| \hat{h}'\|_{\H^1(0,\infty;\R^3)}
+ \| \hat{\omega}\|_{\H^1(0,\infty;\R^3)} & \leq & C.
\end{eqnarray*}
\end{theorem}

\subsection{Proof of Theorem \ref{thstabnonlinX}}
Let us set
\begin{eqnarray*}
\mathbb{H} & = & \H^{2,1}(Q_{\infty}^0) \times \L^2(0,\infty;\H^1(\mathcal{F})) \times \H^1(0,\infty;\R^3) \times \H^1(0,\infty;\R^3).
\end{eqnarray*}
A solution of system \eqref{hhpremsfix}--\eqref{hhdersfix} can be seen as a fixed point of the mapping
\begin{eqnarray*}
\begin{array} {cccc}
\mathcal{N} \ : \ & \mathbb{H} & \rightarrow & \mathbb{H} \\
& (\hat{v},\hat{q},\hat{k}',\hat{\varpi}) & \mapsto & (\hat{u}, \hat{p}, \hat{h}', \hat{\omega})
\end{array}
\end{eqnarray*}
where $(\hat{u}, \hat{p}, \hat{h}', \hat{\omega})$ satisfies
\begin{eqnarray*}
\frac{\p \hat{u}}{\p t} - \lambda \hat{u} -\nu \Delta \hat{u} + \nabla \hat{p}
 =  F(\hat{v},\hat{q},\hat{k}',\hat{\varpi}) & \quad & \text{in } \mathcal{F}\times (0,\infty), \\
\div \ \hat{u}  =  \div \ G(\hat{v},\hat{k}',\hat{\varpi}) & \quad & \text{in } \mathcal{F}\times (0,\infty),
\end{eqnarray*}
\begin{eqnarray*}
\hat{u} = 0  & \quad & \text{on } \p \mathcal{O}\times (0,\infty),  \\
\hat{u} =  \hat{h}'(t) + \hat{\omega} (t) \wedge y+ \mathcal{K}_{\lambda}(\mathbb{P}(\hat{u}-G(\hat{u},\hat{h}',\hat{\omega})),\hat{h}',\hat{\omega}) & & \\
+ \left(e^{\lambda t}\frac{\p X^{\ast}}{\p t}-\zeta \right) + W(\hat{v},\hat{k}',\hat{\varpi})  & \quad & \text{on } \p \mathcal{S}\times (0,\infty),
\end{eqnarray*}
\begin{eqnarray*}
M \hat{h}'' - \lambda M \hat{h}'  =  - \int_{\p \mathcal{S}} \sigma(\hat{u},\hat{p}) n  \d \Gamma + F_M(\hat{v},\hat{q},\hat{k}',\hat{\varpi}) & \quad & \text{in } (0,\infty),   \\
I_0\hat{\omega}' (t) -\lambda I_0\hat{\omega}  =   -  \int_{\p \mathcal{S}} y \wedge \sigma(\hat{u},\hat{p}) n  \d \Gamma + F_I(\hat{v},\hat{k}',\hat{\varpi}) & \quad & \text{in } (0,\infty),
\end{eqnarray*}
\begin{eqnarray*}
\hat{u}(y,0)  =  u_0 (y), \  y\in \mathcal{F} , \quad \hat{h}'(0)=h_1 \in \R^3 ,\quad \hat{\omega}(0) =\omega_0 \in \R^3,
\end{eqnarray*}
with
\begin{eqnarray*}
\zeta & = & \mathcal{K}_{\lambda}(\mathbb{P}(\hat{v}-G(\hat{v},\hat{k}',\hat{\varpi})),\hat{k}',\hat{\varpi}),
\end{eqnarray*}
and $X^{\ast}$ which is obtained from $\zeta$ by Theorem \ref{thdecompsuper}. The system above is actually the nonhomogeneous linear system introduced in section \ref{linearsec}. In particular the estimate \eqref{superestK} gives
\begin{eqnarray}
& & \| \hat{u} \|_{\H^{2,1}(Q_{\infty}^0)} + \| \nabla \hat{p} \|_{\L^2(0,\infty;\L^2(\mathcal{F}))}+
\| \hat{h}' \|_{\H^1(0,\infty;\R^3)} + \| \hat{\omega} \|_{\H^1(0,\infty;\R^3)} \leq \nonumber \\
& &  C_0\left( \| u_0 \|_{\mathbf{H}^1(\mathcal{F})} +|h_1|_{\R^3} + |\omega_0|_{\R^3}
+ \| F(\hat{v},\hat{q},\hat{k}',\hat{\varpi}) \|_{\L^2(0,\infty;\mathbf{L}^2(\mathcal{F}))} \right. \nonumber \\
& & + \|  G(\hat{v},\hat{k}',\hat{\varpi}) \|_{\H^{2,1}(Q_\infty^0)}
+  \| W(\hat{v},\hat{k}',\hat{\varpi}) \|_{\L^2(0,\infty;\mathbf{H}^{3/2}(\p \mathcal{S}))\cap \H^1(0,\infty;\mathbf{H}^{1/2}(\p \mathcal{S}))}
 \nonumber \\
& & + \left\| e^{\lambda t}\frac{\p X^{\ast}}{\p t} - \zeta \right\|_{\L^2(0,\infty;\mathbf{H}^{3/2}(\p \mathcal{S}))\cap \H^1(0,\infty;\mathbf{H}^{1/2}(\p \mathcal{S}))}
\nonumber \\
& & \left.
+ \| F_M(\hat{v},\hat{q},\hat{k}',\hat{\varpi}) \|_{\L^2(0,\infty;\R^3)}
+ \| F_I(\hat{v},\hat{k}',\hat{\varpi}) \|_{\L^2(0,\infty;\R^3)}  \right). \nonumber \\  \label{fixenh}
\end{eqnarray}

\subsubsection{Preliminary estimates}
Some estimates given below are obtained by using a result stated in the Appendix~B of \cite{Grubb} (Proposition B.1), and that we remind in Lemma \ref{lemmaGrubb}. Let us also keep in mind the regularities provided by Proposition \ref{lemmaKtilde} for mappings $\tilde{X}$ and $\tilde{Y}(\tilde{X})$.

\begin{lemma} \label{lemmaH301}
There exists a positive constant $C$ such that for all $(\hat{u},\hat{p},\hat{h}',\hat{\omega})$ in $\mathbb{H}$ we have
\begin{eqnarray}
& & \left\| \left(\Delta - \LL\right)\hat{u} \right\|_{L^2(0,\infty;\mathbf{L}^2(\mathcal{F}))}  \leq   C
\left\| \hat{u} \right\|_{\L^2(0,\infty;\mathbf{H}^2(\mathcal{F}))} \times \\
& & \left( \| \nabla \tilde{Y}(\tilde{X})\nabla \tilde{Y}(\tilde{X})^T - \I_{\R^3} \|_{\L^{\infty}(0,\infty;\mathbf{H}^{2}(\mathcal{F}))}
+ \| \Delta \tilde{Y}(\tilde{X}) \|_{\L^{\infty}(0,\infty;\mathbf{H}^{1}(\mathcal{F}))} \right),  \label{estLL} \\
& & \left\| (\nabla - \GG)\hat{p} \right\|_{L^2(0,\infty;\mathbf{L}^2(\mathcal{F}))} \leq C\| \nabla \tilde{Y}(\tilde{X}) - \I_{\R^3} \|_{\L^{\infty}(0,\infty;\H^{2}(\mathcal{F}))} \left\| \nabla \hat{p} \right\|_{\L^2(0,\infty;\mathbf{L}^2(\mathcal{F}))}.
\nonumber \\ \label{estGG}
\end{eqnarray}
\end{lemma}

\begin{proof}
The only delicate point that consists in verifying that $\Delta \tilde{Y}(\tilde{X})$ in $\L^{\infty}(0,\infty;\mathbf{H}^{1}(\mathcal{F}))$. For that, let us consider the i-th component of $\Delta \tilde{Y}(\tilde{X})$; We write
\begin{eqnarray*}
\Delta \tilde{Y}_i(\tilde{X}(\cdot,t),t) & = & \trace \left( \nabla^2 \tilde{Y}_i(\tilde{X}(\cdot,t),t) \right)
\end{eqnarray*}
with
\begin{eqnarray*}
\nabla^2 \tilde{Y}_i(\tilde{X}(\cdot,t),t) & = & \left( \nabla \left( \nabla \tilde{Y}_i(\tilde{X}(\cdot,t),t) \right)\right) \nabla \tilde{Y}(\tilde{X}(\cdot,t),t)\\
& = & \left( \nabla \left( \nabla \tilde{Y}_i(\tilde{X}(\cdot,t),t) - \I_{\R^3} \right)\right) \nabla \tilde{Y}(\tilde{X}(\cdot,t),t),
\end{eqnarray*}
and we apply Lemma \ref{lemmaGrubb} with $s=1$, $\mu = 0$ and $\kappa = 1$ to obtain
\begin{eqnarray*}
\| \Delta \tilde{Y}_i(\tilde{X}(\cdot,t),t) \|_{\mathbf{H}^{1}(\mathcal{F})} & \leq & C \| \nabla \tilde{Y}(\tilde{X}(\cdot,t),t) -\I_{\R^3} \|_{\mathbf{H}^{2}(\mathcal{F})} \| \nabla \tilde{Y}(\tilde{X}(\cdot,t),t) \|_{\mathbf{H}^{2}(\mathcal{F})}.
\end{eqnarray*}
\end{proof}

\begin{lemma} \label{lemmaH302}
There exists a positive constant $C$ such that for all $(\hat{u},\hat{p},\hat{h}',\hat{\omega})$ in $\mathbb{H}$ we have
\begin{eqnarray}
& & \left\| \MM \hat{u} \right\|_{L^2(0,\infty;\mathbf{L}^2(\mathcal{F}))}  \leq C \| \hat{u} \|_{\L^{2}(0,\infty;\mathbf{H}^2(\mathcal{F}))}\left(\| \nabla  \tilde{Y}(\tilde{X})-\I_{\R^3} \|_{\L^{\infty}(0,\infty;\mathbf{H}^{2}(\mathcal{F}))}+1 \right) \times \nonumber \\
 & & \left( \|\hat{h}'\|_{\L^{\infty}(0,\infty;\R^3)} + \|\hat{\omega}\|_{\L^{\infty}(0,\infty;\R^3)}
 \| \tilde{X} \|_{\L^{\infty}(0,\infty;\mathbf{H}^{1}(\mathcal{F}))}
 + \left\|\frac{\p \tilde{X}}{\p t} \right\|_{\L^{\infty}(0,\infty;\mathbf{H}^{1}(\mathcal{F}))}\right) , \nonumber \\ \label{estiMM} \\
& &  \left\| \NN \hat{u} \right\|_{L^2(0,\infty;\mathbf{L}^2(\mathcal{F}))}  \leq
 C \left\| \hat{u} \right\|_{\L^{\infty}(0,\infty;\mathbf{H}^1(\mathcal{F}))} \left\| \hat{u} \right\|_{\L^{2}(0,\infty;\mathbf{H}^{2}(\mathcal{F}))} \times  \nonumber \\
& &  \qquad \qquad \left(\| \nabla  \tilde{Y}(\tilde{X})-\I_{\R^3}\|_{\L^{\infty}(0,\infty;\mathbf{H}^{2}(\mathcal{F}))}+1\right), \label{estiNN}  \\
& &  \| \hat{\omega} \wedge \hat{u} \|_{L^2(0,\infty;\mathbf{L}^2(\mathcal{F}))}  \leq  C \| \hat{\omega}\|_{\H^{1}(0,\infty;\R^3)} \| \hat{u} \|_{L^{2}(0,\infty;\mathbf{L}^2(\mathcal{F}))}. \label{estiOU}
\end{eqnarray}
\end{lemma}

\begin{proof}
There is no particular difficulty for obtaining these estimates.
\end{proof}

\begin{lemma} \label{lemmaH3}
There exists a positive constant $C$ such that for all $\hat{u} \in \H^{2,1}(Q_{\infty}^0)$ we have
\begin{eqnarray}
\left\| G(\hat{u}) \right\|_{\L^2(0,\infty;\mathbf{H}^2(\mathcal{F}))} & \leq & C\|\hat{u} \|_{\L^2(0,\infty;\mathbf{H}^2(\mathcal{F}))} \| \nabla \tilde{Y}(\tilde{X}) - \I_{\R^3} \|_{\L^{\infty}(0,\infty;\mathbf{H}^2(\mathcal{F}))}, \nonumber \\
 \left\| G(\hat{u}) \right\|_{\H^1(0,\infty;\mathbf{L}^2(\mathcal{F}))} & \leq & C \left( \|\hat{u} \|_{\H^1(0,\infty;\mathbf{L}^{2}(\mathcal{F}))} \|  \nabla \tilde{Y}(\tilde{X}) - \I_{\R^3} \|_{\L^{\infty}(0,\infty;\mathbf{H}^2(\mathcal{F}))} \right. \nonumber \\
  & &  \qquad + \left. \|\hat{u} \|_{\L^2(0,\infty;\mathbf{H}^2(\mathcal{F}))} \| \nabla \tilde{Y}(\tilde{X}) - \I_{\R^3} \|_{\W^{1,\infty}(0,\infty;\mathbf{L}^2(\mathcal{F}))} \right). \nonumber
\end{eqnarray}
\end{lemma}

\begin{proof}
The quantity $\nabla \tilde{Y}(\tilde{X})$ lies in $\L^2(0,\infty;\mathbf{H}^{2}(\mathcal{F}))$. We apply Lemma \ref{lemmaGrubb} with $s =2$, $\mu = 0$ and $\kappa = 0$ in order to get
\begin{eqnarray*}
\left\|G(\hat{u})(\cdot,t) \right\|_{\mathbf{H}^{2}(\mathcal{F})} & \leq & C \left\| \hat{u} \right\|_{\mathbf{H}^{2}(\mathcal{F})}
\|\nabla \tilde{Y}(\tilde{X}(\cdot,t),t) - I \|_{\mathbf{H}^{2}(\mathcal{F})}, \\
\left\|G(\hat{u}) \right\|_{\L^2(0,\infty;\mathbf{H}^{2}(\mathcal{F}))} & \leq & C \left\| \hat{u} \right\|_{\L^2(0,\infty;\mathbf{H}^{2}(\mathcal{F}))}
\| \nabla \tilde{Y}(\tilde{X}) - \I_{\R^3} \|_{\L^{\infty}(0,\infty;\mathbf{H}^{2}(\mathcal{F}))}.
\end{eqnarray*}
For proving the regularity $\H^1(0,\infty;\mathbf{L}^2(\mathcal{F}))$, we first write
\begin{eqnarray*}
\frac{\p G(\hat{u})}{\p t}  & = & \left( \nabla \tilde{Y}(\tilde{X}) - \I_{\R^3} \right)\frac{\p \hat{u}}{\p t} + \left(\frac{\p }{\p t}\left( \nabla \tilde{Y}(\tilde{X})-I \right)\right) \hat{u}.
\end{eqnarray*}
The quantity $\nabla \tilde{Y}(\tilde{X})$ lies in $\H^2(0,\infty;\mathbf{L}^2(\mathcal{F})) \hookrightarrow \W^{1,\infty}(0,\infty;\mathbf{L}^2(\mathcal{F}))$, so that we have the estimates
\begin{eqnarray*}
\left\| \frac{\p G(\hat{u})}{\p t} \right\|_{\L^2(0,\infty;\mathbf{L}^2(\mathcal{F}))} & \leq & C \left( \| \nabla \tilde{Y}(\tilde{X}) - \I_{\R^3} \|_{\L^{\infty}(0,\infty;\mathbf{L}^{\infty}(\mathcal{F}))} \left\| \frac{\p \hat{u}}{\p t} \right\|_{\L^2(0,\infty;\mathbf{L}^2(\mathcal{F}))} \right. \\
& & + \left. \| \nabla \tilde{Y}(\tilde{X}) - \I_{\R^3} \|_{W^{1,\infty}(0,\infty;\mathbf{L}^2(\mathcal{F}))} \|\hat{u} \|_{\L^2(0,\infty;\mathbf{L}^{\infty}(\mathcal{F}))}  \right).
\end{eqnarray*}
\end{proof}

\begin{lemma} \label{estFMFI}
There exists a positive constant $C$ such that for all $(\hat{u},\hat{p},\hat{h}',\hat{\omega})$ in $\mathbb{H}$ we have
\begin{eqnarray}
& & \left\| W(\hat{\omega}) \right\|_{\H^1(0,\infty;\mathbf{H}^{3/2}(\p \mathcal{S}))} \leq C\|\hat{\omega}\|_{\H^1(0,\infty;\R^3)}
\| X^{\ast} - \Id \|_{\tilde{\mathcal{W}}(Q_{\infty}^0)}, \label{estW} \\
& & \left\| F_M(\hat{u},\hat{p}, \hat{h}', \hat{\omega}) \right\|_{\L^2(0,\infty;\R^3)}  \leq \nonumber \\
& & C\left(  \|\hat{h}'\|_{\L^{2}(0,\infty;\R^3)} \|\hat{\omega}\|_{\H^{1}(0,\infty;\R^3)}
+  \left( \|\hat{u} \|_{\L^2(0,\infty;\mathbf{H}^2(\mathcal{F}))} + \|\hat{p} \|_{\L^2(0,\infty;\H^1(\mathcal{F}))} \right) \times \right. \nonumber \\
& &  \left. \left( \| \nabla \tilde{Y}( \tilde{X}) \|_{L^{\infty}(\p \mathcal{S}\times (0,\infty))}
\left(1+\| \nabla \tilde{Y}(\tilde{X}) - \I_{\R^3} \|_{L^{\infty}(\p \mathcal{S}\times (0,\infty))}\right)
 \right) \right), \label{estFM} \\
& & \left\| F_I(\hat{u},\hat{p},\hat{\omega}) \right\|_{\L^2(0,\infty;\R^3)}  \leq \nonumber \\
& & C \left( (1+\lambda)\| I^{\ast} - I_0 \|_{L^{\infty}(0,\infty;\R^{9})} \| \hat{\omega} \|_{\H^1(0,\infty;\R^3)} \right. \nonumber \\
& & + \|{I^{\ast}}' \|_{L^2(0,\infty;\R^9)} \| \hat{\omega}\|_{\L^{\infty}(0,\infty;\R^3)} + \| I^{\ast}\|_{\L^{\infty}(0,\infty;\R^9)} \|\hat{\omega}\|_{\L^{\infty}(0,\infty;\R^3)} \|\hat{\omega}\|_{\L^{2}(0,\infty;\R^3)} \nonumber \\
& & +  \left( \|\hat{u} \|_{\L^2(0,\infty;\mathbf{H}^2(\mathcal{F}))} +
\|\hat{p} \|_{\L^2(0,\infty;\H^1(\mathcal{F}))} \right) \times \nonumber \\
& & \left( \| \nabla \tilde{Y}( \tilde{X}) - \I_{\R^3} \|_{L^{\infty}(\p \mathcal{S}\times (0,\infty))}
+ \| \nabla \tilde{Y}( \tilde{X}) \|_{L^{\infty}(\p \mathcal{S}\times (0,\infty))} \times \right. \nonumber \\
& & \left. \left( \| \nabla \tilde{Y}(\tilde{X}) - I \|_{L^{\infty}(\p \mathcal{S}\times (0,\infty))}
+ \|  X^{\ast} - \Id \|_{L^{\infty}(\p \mathcal{S}\times (0,\infty))} \right)  \right),
\end{eqnarray}
with
\begin{eqnarray*}
\|{I^{\ast}}' \|_{L^2(0,\infty;\R^9)} & \leq & C \| X^{\ast} \|_{\L^{\infty}(0,\infty;\mathbf{L}^2(\mathcal{S}))}
\left\| \frac{\p X^{\ast}}{\p t} \right\|_{\L^2(0,\infty;\mathbf{L}^2(\mathcal{S}))} , \\
\| I^{\ast} - I_0 \|_{L^{\infty}(0,\infty;\R^{9})} & \leq & C \| X^{\ast} \|_{\L^{\infty}(0,\infty;\mathbf{L}^2(\mathcal{S}))}
\left\| e^{\lambda t}\frac{\p X^{\ast}}{\p t} \right\|_{\L^2(0,\infty;\mathbf{L}^2(\mathcal{S}))} , \\
\| I^{\ast} \|_{L^{\infty}(0,\infty;\R^{9})} & \leq &  C\| X^{\ast} \|_{\L^{\infty}(0,\infty;\mathbf{L}^2(\mathcal{S}))}.
\end{eqnarray*}
\end{lemma}

\begin{proof}
There is no particular difficulty for proving the other two estimates, if we refer to the respective expressions of $W$, $F_M$ and $F_I$ given by \eqref{rhsW}, \eqref{rhsFM} and \eqref{rhsFI}.
\end{proof}


For some radius $R>0$, let us define the ball
\begin{eqnarray*}
& & B_R = \left\{ (\hat{u}, \hat{p}, \hat{h}', \hat{\omega}) \in \mathbb{H} \mid \hat{u} = 0 \text{ on } \p \mathcal{O}, \text{ and} \right.  \\
& & \left. \|\hat{u} \|_{\H^{2,1}(Q_{\infty}^0)} + \| \hat{p}\|_{\L^2(0,\infty;\H^1(\mathcal{F}))} + \|\hat{h}'\|_{\H^1(0,\infty;\R^3)} + \|\hat{\omega}\|_{\H^1(0,\infty;\R^3)} \leq 2RC_0  \right\}
\end{eqnarray*}
with
\begin{eqnarray*}
R & = & \left(  \|u_0\|_{\mathbf{H}^1(\mathcal{F})}+ |h_1|_{\R^3} + |\omega_0|_{\R^3} \right),
\end{eqnarray*}
and where the constant $C_0$ appears in the estimate \eqref{fixenh}. Note that $B_R$ is a closed subset of $\mathbb{H}$.\\

First, for $R$ small enough, we claim that for $(\hat{v}, \hat{q}, \hat{k}', \hat{\varpi}) \in B_R$ the function
\begin{eqnarray*}
\zeta & = & \mathcal{K}_{\lambda}(\mathbb{P}(\hat{v}-G(\hat{v},\hat{k}',\hat{\varpi})),\hat{k}',\hat{\varpi})
\end{eqnarray*}
is small enough in $\L^2(0,\infty;\mathbf{H}^{5/2}(\p \mathcal{S}))\cap \H^1(0,\infty;\mathbf{H}^{1/2}(\p \mathcal{S}))$, because of the estimates of Lemma \ref{lemmaH3} and the continuity of the operator $\mathcal{K}_{\lambda}$ (see Part I, section 5). So by Theorem \ref{thdecompsuper} we can define the corresponding {\it admissible} deformation $X^{\ast}$, the mappings $\tilde{X}$ and $\tilde{Y}$ which stem, and thus - in virtue of Remark \ref{remarkcc} - we can define properly the mapping $\mathcal{N}$ in $B_R$.

\subsubsection{Stability of the set $B_R$ by the mapping $\mathcal{N}$}
Let be $(\hat{v}, \hat{q}, \hat{k}', \hat{\varpi}) \in B_R$, for $R$ small enough. In the estimates provided by the previous lemmas \ref{lemmaH301}, \ref{lemmaH302}, \ref{lemmaH3}, \ref{estFMFI}, note that from the estimate \eqref{app11}, \eqref{app12}, \eqref{app13} of Proposition \ref{lemmaKtilde} combined to \eqref{estsuper} and \eqref{estsuper0} of Theorem \ref{thdecompsuper} we can deduce
\begin{eqnarray*}
\|F(\hat{v}, \hat{q}, \hat{k}', \hat{\varpi})\|_{\L^2(0,\infty;\mathbf{L}^2(\mathcal{F}))}
& = & o(\|u_0\|_{\mathbf{H}^1(\mathcal{F})}+ |h_1|_{\R^3} + |\omega_0|_{\R^3}), \\
\|G(\hat{v}, \hat{k}', \hat{\varpi})\|_{\H^{2,1}(Q_{\infty}^0)}
& = & o(\|u_0\|_{\mathbf{H}^1(\mathcal{F})}+ |h_1|_{\R^3} + |\omega_0|_{\R^3}), \\
\| W(\hat{v}, \hat{k}', \hat{\varpi})\|_{\L^2(\mathbf{H}^{3/2}(\p \mathcal{S}))\cap \H^1(\mathbf{H}^{1/2}(\p \mathcal{S}))}
& = & o(\|u_0\|_{\mathbf{H}^1(\mathcal{F})}+ |h_1|_{\R^3} + |\omega_0|_{\R^3}), \\
\left\| e^{\lambda t}\frac{\p X^{\ast}}{\p t} -\zeta
\right\|_{\L^2(\mathbf{H}^{3/2}(\p \mathcal{S}))\cap \H^1(\mathbf{H}^{1/2}(\p \mathcal{S}))}
& = & o(\|u_0\|_{\mathbf{H}^1(\mathcal{F})}+ |h_1|_{\R^3} + |\omega_0|_{\R^3}), \\
\|F_M(\hat{v}, \hat{q}, \hat{k}', \hat{\varpi})\|_{\L^2(0,\infty;\R^3)}
& = & o(\|u_0\|_{\mathbf{H}^1(\mathcal{F})}+ |h_1|_{\R^3} + |\omega_0|_{\R^3}), \\
\|F_I(\hat{v}, \hat{k}', \hat{\varpi})\|_{\L^2(0,\infty;\R^3)}
& = & o(\|u_0\|_{\mathbf{H}^1(\mathcal{F})}+ |h_1|_{\R^3} + |\omega_0|_{\R^3}).
\end{eqnarray*}
In particular, for $R$ small enough the mapping $\mathcal{N}$ is well-defined, and so we can write the estimate \eqref{fixenh} that we remind
\begin{eqnarray}
& & \left\| \mathcal{N}(\hat{v}, \hat{q}, \hat{k}', \hat{\varpi}) \right\|_{\mathbb{H}} \leq  \nonumber \\
& &  C_0\left( \| u_0 \|_{\mathbf{H}^1(\mathcal{F})} +|h_1|_{\R^3} + |\omega_0|_{\R^3}
+ \| F(\hat{v},\hat{q},\hat{k}',\hat{\varpi}) \|_{\L^2(0,\infty;\mathbf{L}^2(\mathcal{F}))} \right. \nonumber \\
& & + \|  G(\hat{v},\hat{k}',\hat{\varpi}) \|_{\H^{2,1}(Q_\infty^0)}
+  \| W(\hat{v},\hat{k}',\hat{\varpi}) \|_{\L^2(0,\infty;\mathbf{H}^{3/2}(\p \mathcal{S}))\cap \H^1(0,\infty;\mathbf{H}^{1/2}(\p \mathcal{S}))}
 \nonumber \\
& & + \left\| e^{\lambda t}\frac{\p X^{\ast}}{\p t} - \zeta \right\|_{\L^2(0,\infty;\mathbf{H}^{3/2}(\p \mathcal{S}))\cap \H^1(0,\infty;\mathbf{H}^{1/2}(\p \mathcal{S}))}
\nonumber \\
& & \left.
+ \| F_M(\hat{v},\hat{q},\hat{k}',\hat{\varpi}) \|_{\L^2(0,\infty;\R^3)}
+ \| F_I(\hat{v},\hat{k}',\hat{\varpi}) \|_{\L^2(0,\infty;\R^3)}  \right). \label{estNNN}
\end{eqnarray}
Thus we have
\begin{eqnarray*}
& & \left\| \mathcal{N}(\hat{v}, \hat{q}, \hat{k}', \hat{\varpi}) \right\|_{\mathbb{H}} \leq  C_0\left( R + o(R)\right).
\end{eqnarray*}
This shows that for $R = \left( \|u_0\|_{\mathbf{H}^1(\mathcal{F})}+ |h_1|_{\R^3} + |\omega_0|_{\R^3} \right)$ small enough the ball $B_R$ is stable by the mapping $\mathcal{N}$.

\subsubsection{Lipschitz stability for the mapping $\mathcal{N}$}
Let $(\hat{v}_1, \hat{q}_1, \hat{k}'_1, \hat{\varpi}_1)$ and $(\hat{v}_2, \hat{q}_2, \hat{k}'_2, \hat{\varpi}_2)$ be in $B_R$. We set
\begin{eqnarray*}
(\hat{u}_1, \hat{p}_1, \hat{h}'_1, \hat{\omega}_1)  =  \mathcal{N}(\hat{v}_1, \hat{q}_1, \hat{k}'_1, \hat{\varpi}_1), & \quad &
(\hat{u}_2, \hat{p}_2, \hat{h}'_2, \hat{\omega}_2)  =  \mathcal{N}(\hat{v}_2, \hat{q}_2, \hat{k}'_2, \hat{\varpi}_2),
\end{eqnarray*}
and
\begin{eqnarray*}
\begin{array} {ccccccc}
\hat{u} = \hat{u}_2 - \hat{u}_1 ,& \quad & \hat{p} = \hat{p}_2 - \hat{p}_1,
 & \quad & \hat{h}' = \hat{h}_2' - \hat{h}_1', & \quad & \hat{\omega} = \hat{\omega}_2 - \hat{\omega}_1, \\
\hat{v} = \hat{v}_2 - \hat{v}_1 ,& \quad & \hat{q} = \hat{q}_2 - \hat{q}_1,
 & \quad & \hat{k}' = \hat{k}_2' - \hat{k}_1', & \quad & \hat{\varpi} =  \hat{\varpi} _2 -  \hat{\varpi} _1.
\end{array}
\end{eqnarray*}
We also denote
\begin{eqnarray*}
\zeta_1 & = & \mathcal{K}_{\lambda}(\mathbb{P}(\hat{v}_1-G(\hat{v}_1,\hat{k}'_1,\hat{\varpi}_1)),\hat{k}'_1,\hat{\varpi}_1), \\
\zeta_2 & = & \mathcal{K}_{\lambda}(\mathbb{P}(\hat{v}_2-G(\hat{v}_2,\hat{k}'_2,\hat{\varpi}_2)),\hat{k}'_2,\hat{\varpi}_2),
\end{eqnarray*}
and next $X^{\ast}_1$ and $X^{\ast}_2$ the {\it admissible} deformation given by Theorem \ref{thdecompsuper} which are defined by $\zeta_1$ and $\zeta_2$ respectively.\\
For $R$ small enough the quadruplet $(\hat{u}, \hat{p}, \hat{h}', \hat{\omega})$ satisfies the system
\begin{eqnarray*}
\frac{\p \hat{u}}{\p t} - \lambda \hat{u} -\nu \Delta \hat{u} + \nabla \hat{p}
 =  \overline{F}, & \quad & \text{in } \mathcal{F}\times (0,\infty), \\
\div \ \hat{u}  =  \div \ \overline{G}, & \quad & \text{in } \mathcal{F}\times (0,\infty),
\end{eqnarray*}
\begin{eqnarray*}
\hat{u} = 0 , & \quad & \text{on } \p \mathcal{O}\times (0,\infty),  \\
\hat{u}  =  \hat{h}'(t) + \hat{\omega} (t) \wedge y + \mathcal{K}_{\lambda}(\mathbb{P}(\hat{u}-\overline{G}),\hat{h}', \hat{\omega}) + \overline{W}  , & \quad & (y,t)\in \p \mathcal{S}\times (0,\infty),
\end{eqnarray*}
\begin{eqnarray*}
M \hat{h}'' - \lambda M \hat{h}'  =
- \int_{\p \mathcal{S}} \sigma(\hat{u},\hat{p}) n  \d \Gamma + \overline{F}_M, & \quad & \text{in } (0,\infty)   \\
I_0\hat{\omega}' (t) -\lambda I_0\hat{\omega}  =
-  \int_{\p \mathcal{S}} y \wedge \sigma(\hat{u},\hat{p}) n  \d \Gamma + \overline{F}_I, & \quad & \text{in } (0,\infty)
\end{eqnarray*}
\begin{eqnarray*}
\hat{u}(y,0)  = 0, \  \text{in } \mathcal{F} , \quad \hat{h}'(0)=0 \in \R^3 ,\quad \hat{\omega}(0) = 0 \in \R^3 ,
\end{eqnarray*}
with
\begin{eqnarray*}
\overline{F} & = & F(\hat{v}_2, \hat{q}_2, \hat{k}'_2, \hat{\varpi}_2) - F(\hat{v}_1, \hat{q}_1, \hat{k}'_1, \hat{\varpi}_1), \\
\overline{G} & = & G(\hat{v}_2, \hat{k}'_2, \hat{\varpi}_2) - G(\hat{v}_1, \hat{k}'_1, \hat{\varpi}_1), \\
\overline{W} & = & W(\hat{v}_2, \hat{k}'_2, \hat{\varpi}_2) - W(\hat{v}_1, \hat{k}'_1, \hat{\varpi}_1) \\
& & + \left(e^{\lambda t}\frac{\p X^{\ast}_2}{\p t}-\zeta_2 \right)
- \left(e^{\lambda t}\frac{\p X^{\ast}_1}{\p t}-\zeta_1 \right), \\
\overline{F}_M & = & F_M(\hat{v}_2, \hat{q}_2, \hat{k}'_2, \hat{\varpi}_2) - F_M(\hat{v}_1, \hat{q}_1, \hat{k}'_1, \hat{\varpi}_1) , \\
\overline{F}_I & = & F_I(\hat{v}_2, \hat{q}_2, \hat{k}'_2, \hat{\varpi}_2) - F_I(\hat{v}_1, \hat{q}_1, \hat{k}'_1, \hat{\varpi}_1).
\end{eqnarray*}
The deformations $X^{\ast}_1$ and $X^{\ast}_2$ induce respectively the mappings $(\tilde{X}_1,\tilde{Y}_1)$ and $(\tilde{X}_2,\tilde{Y}_2)$ which - partially - define the right-hand-sides above.\\
The right-hand-sides $\overline{F}$, $\overline{G}$, $\overline{W}$, $\overline{F}_M$ and $\overline{F}_I$ can be expressed as quantities which are multiplicative of the differences
\begin{eqnarray*}
& & \hat{v}, \quad \hat{q}, \quad \hat{k}', \quad \hat{\omega}, \quad (X^{\ast}_2- X^{\ast}_1),
\quad (\tilde{X}_2- \tilde{X}_1), \\
& &  \left( \nabla \tilde{Y}_2(\tilde{X}_2) - \nabla \tilde{Y}_1(\tilde{X}_1)\right), \quad
 \left( \Delta \tilde{Y}_2(\tilde{X}_2) - \Delta \tilde{Y}_1(\tilde{X}_1)\right) .
\end{eqnarray*}
For instance, the nonhomogeneous divergence condition $\overline{G}$ can be written as
\begin{eqnarray*}
\overline{G} & = & \left( \nabla \tilde{Y}_2(\tilde{X}_2) - \nabla \tilde{Y}_1(\tilde{X}_1) \right) \hat{v}_2 + (\nabla \tilde{Y}_1(\tilde{X}_1)- \I_{\R^3})\hat{v}.
\end{eqnarray*}
Then the estimates of Lemmas \ref{lemmaH301}, \ref{lemmaH302}, \ref{lemmaH3}, \ref{estFMFI} can be adapted for this right-hand-sides, so that the estimates \eqref{app11}, \eqref{app12}, \eqref{app13} of Proposition \ref{lemmaKtilde} combined to the estimates \eqref{estsuper12} and \eqref{estsuper120} of Theorem \ref{thdecompsuper} enable us to prove that for $R$ small enough we have
\begin{eqnarray*}
\| \overline{F} \|_{\L^2(0,\infty;\mathbf{L}^2(\mathcal{F}))} & = &
o(\|(\hat{v}, \hat{q}, \hat{k}', \hat{\varpi})\|_{\mathbb{H}}), \\
\| \overline{G} \|_{\H^{2,1}(Q_{\infty}^0)} & = &
o(\|(\hat{v}, \hat{q}, \hat{k}', \hat{\varpi})\|_{\mathbb{H}}), \\
\| \overline{W} \|_{\L^2(\mathbf{H}^{3/2}(\p \mathcal{S}))\cap \H^1(\mathbf{H}^{1/2}(\p \mathcal{S}))} & = &
o(\|(\hat{v}, \hat{q}, \hat{k}', \hat{\varpi})\|_{\mathbb{H}}), \\
\| \overline{F}_M \|_{\L^2(0,\infty;\R^3)} & = &
o(\|(\hat{v}, \hat{q}, \hat{k}', \hat{\varpi})\|_{\mathbb{H}}), \\
\| \overline{F}_I \|_{\L^2(0,\infty;\R^3)} & = &
o(\|(\hat{v}, \hat{q}, \hat{k}', \hat{\varpi})\|_{\mathbb{H}}), \\.
\end{eqnarray*}
Then the estimate \eqref{estNNN} can be applied for the quadruplet $(\hat{u}, \hat{p}, \hat{h}', \hat{\omega})$:
\begin{eqnarray*}
\| (\hat{u}, \hat{p}, \hat{h}', \hat{\omega}) \|_{\mathbb{H}} & \leq &
 C_0\left( \| \overline{F} \|_{\L^2(0,\infty;\mathbf{L}^2(\mathcal{F}))} +  \| \overline{G} \|_{\H^{2,1}(Q_{\infty}^0)} \right. \\
& & + \left. \| \overline{W} \|_{\L^2(0,\infty;\mathbf{H}^{3/2}(\p \mathcal{S}))\cap \H^1(0,\infty;\mathbf{H}^{1/2}(\p \mathcal{S}))} \right.
\nonumber \\
& & \left.   + \| \overline{F}_M \|_{\L^2(0,\infty;\R^3)}
+ \| \overline{F}_I \|_{\L^2(0,\infty;\R^3)} \right).
\end{eqnarray*}
Then we have for $R$ small enough
\begin{eqnarray*}
\| (\hat{u}, \hat{p}, \hat{h}', \hat{\omega}) \|_{\mathbb{H}}  \leq  C_0\times o(\|(\hat{v}, \hat{q}, \hat{k}', \hat{\varpi})\|_{\mathbb{H}})
 \leq  \frac{1}{2} \|(\hat{v}, \hat{q}, \hat{k}', \hat{\varpi})\|_{\mathbb{H}},
\end{eqnarray*}
and thus the mapping $\mathcal{N}$ is a contraction in $B_R$. It admits a unique fixed point, that is to say there exists a unique solution of system \eqref{hhpremsfix}--\eqref{hhdersfix}. The announced estimate is easily deduced from what precedes.

\section{Conclusion} \label{secconclusion}
After having proven Theorem \ref{thstabnonlinX} for the unknowns $(\hat{u},\hat{p},\hat{h}',\hat{\omega})$ of system \eqref{hhpremsfix}--\eqref{hhdersfix}, let us remind the relations given by \eqref{tildeu}-\eqref{tildeh} and \eqref{chhat}
\begin{eqnarray*}
\begin{array} {ll}
e^{\lambda t}u(x,t) =  \mathbf{R}(t)\hat{u}(Y(x,t),t), & e^{\lambda t}p(x,t) =  \hat{p}(Y(x,t),t), \\
e^{\lambda t}h'(t) =  \mathbf{R}(t)\hat{h}'(t), & e^{\lambda t}\omega(t) =  \mathbf{R}(t)\hat{\omega}(t),
\end{array}
\end{eqnarray*}
and the relations given by \eqref{definitilde}
\begin{eqnarray*}
\begin{array} {ll}
X(y,t) = h(t) + \mathbf{R}(t)\tilde{X}(y,t), & \quad (y,t) \in \mathcal{F}\times (0,\infty), \\
Y(x,t) = \tilde{Y}(\mathbf{R}(t)^T(x-h(t)),t), & \quad (x,t) \in Q_{\infty} = \displaystyle \bigcup_{t \geq 0} \mathcal{F}(t)\times \{ t \}.
\end{array}
\end{eqnarray*}
Let us deduce the stabilization of the main system in the sense of Definition \ref{defstab}.

\begin{theorem} \label{maintheorem}
For $(u_0,h_1,\omega_0)$ small enough in $\mathbf{H}^1(\mathcal{F}) \times \R^3 \times \R^3$, system \eqref{prems}--\eqref{ders} is stabilizable with an arbitrary exponential decay rate, in the sense of Definition \ref{defstab}.
\end{theorem}

\begin{proof}
First, from Remark \ref{remarkc}, we can see that for $(\hat{u},\hat{p},\hat{h}',\hat{\omega})$ given the quadruplet is determined in a unique way; In particular we can define the rotation $\mathbf{R}$ associated with the angular velocity $\omega$ as the solution of the following problem
\begin{eqnarray*}
\begin{array} {ccccc}
\displaystyle \frac{\d}{\d t}(\mathbf{R}) & = & \mathbb{S}\left(
e^{-\lambda t}\mathbf{R}\hat{\omega}\right) \mathbf{R} & = & e^{-\lambda t}\mathbf{R} \mathbb{S}\left(\hat{\omega}\right) \\
\mathbf{R}(t=0) & = & \I_{\R^3}, & &
\end{array}
\quad \text{with }
\mathbb{S}(\tilde{\omega}) = \left(
\begin{matrix}
0 & -\tilde{\omega}_3 & \tilde{\omega}_2 \\
\tilde{\omega}_3 & 0 & -\tilde{\omega}_1 \\
-\tilde{\omega}_2 & \tilde{\omega}_1 & 0
\end{matrix} \right).
\end{eqnarray*}
Now by considering
\begin{eqnarray*}
\frac{\d}{\d t}\left(e^{\lambda t} h'\right)  =
e^{-\lambda t}\mathbf{R}\left( \hat{\omega}\wedge\hat{h}'\right) + \mathbf{R}\hat{h}'', & \quad &
\frac{\d}{\d t}\left(e^{\lambda t} \omega\right)  =  \mathbf{R} \hat{\omega}',
\end{eqnarray*}
we can deduce from the estimate of Theorem \ref{thstabnonlinX} that
\begin{eqnarray*}
\| e^{\lambda t} h' \|_{\H^1(0,\infty;\R^3)} \leq C, & \quad & \| e^{\lambda t} \omega \|_{\H^1(0,\infty;\R^3)} \leq C,
\end{eqnarray*}
and by considering
\begin{eqnarray*}
\frac{\d}{\d t}\left(e^{\lambda t} u(x,t)\right)_{\left| x=X(y,t)\right.} & = & e^{-\lambda t}\mathbf{R}\left(\hat{\omega}(t)\wedge \hat{u}(y,t)\right)
+ \mathbf{R} \frac{\p \hat{u}}{\p t}(y,t)\\
& & - \mathbf{R}\nabla\hat{u}(y,t)\nabla \tilde{Y}(\tilde{X}(y,t),t) \times \\
& & \left(e^{-\lambda t}\hat{h}'+e^{-\lambda t}\hat{\omega}\wedge \tilde{X}(y,t)
+ e^{\lambda t} \frac{\p \tilde{X}}{\p t}\right)
\end{eqnarray*}
we can deduce from the estimate of Theorem \ref{thstabnonlinX} and the estimates of Proposition \ref{lemmaKtilde} that, since we have $\det \nabla \tilde{X} = 1$, the following estimates hold
\begin{eqnarray*}
\| e^{\lambda t}u\|_{\L^2(0,\infty;\mathbf{L}^2(\mathcal{F}(t)))} & = & \| \hat{u} \|_{\L^2(0,\infty;\mathbf{L}^2(\mathcal{F}(t)))} \leq C, \\
\left\| \frac{\d}{\d t}\left(e^{\lambda t} u\right) \right\|_{\L^2(0,\infty;\mathbf{L}^2(\mathcal{F}(t)))} & \leq & C.
\end{eqnarray*}
Finally by considering the equalities
\begin{eqnarray*}
e^{\lambda t}\nabla p(X(y,t),t) & = & \mathbf{R}(t)\nabla \tilde{Y}(\tilde{X}(y,t),t)^T\nabla \hat{p}(y,t), \\
e^{\lambda t}\nabla u(X(y,t),t) & = & \mathbf{R}(t)\nabla \hat{u}(u,t) \nabla \tilde{Y}(\tilde{X}(y,t),t) \mathbf{R}(t)^T,
\end{eqnarray*}
and for the $i$-st component of vector $u$ the following equality
\begin{eqnarray*}
e^{\lambda t} \nabla^2 u_i(X(y,t),t) & = & \nabla \left( \nabla u_i(X(y,t),t)\right) \nabla \tilde{Y}(\tilde{X}(y,t),t)\mathbf{R}(t)^T
\end{eqnarray*}
we can use the estimates of Lemma \ref{lemmaKtilde} to conclude the proof with the estimate
\begin{eqnarray*}
\| e^{\lambda t } p \|_{\L^2(0,\infty;\H^1(\mathcal{F}(t)))} \leq C, & \quad &
\| e^{\lambda t } u \|_{\L^2(0,\infty;\mathbf{H}^2(\mathcal{F}(t)))} \leq C.
\end{eqnarray*}
\end{proof}




\section{Appendix~A: The change of variables} \label{subsecchvar}
Let us consider an {\it admissible} deformation $X^{\ast} \in \mathcal{W}_{\lambda}(S_{\infty}^0)$ - in the sense of Definition \ref{defcontrol} - which satisfies, in particular, for all $t>0$ the following condition
\begin{eqnarray}
\int_{\p \mathcal{S}} \com \nabla X^{\ast}(y,t)^T\frac{\p X^{\ast}}{\p t}(y,t) \cdot  n  \d \Gamma(y) & = & 0. \label{eqreccomp}
\end{eqnarray}
The regularity considered for the datum $\displaystyle e^{\lambda t}\frac{\p X^{\ast}}{\p t}$ in this section is
\begin{eqnarray*}
\mathcal{H}(S_{\infty}^0) & = & \L^2(0,\infty;\mathbf{H}^{3}(\mathcal{S}))\cap \H^1(0,\infty;\mathbf{H}^{1}(\mathcal{S})).
\end{eqnarray*}

The goal of this subsection is to extend to the whole domain $\overline{\mathcal{O}}$ the mappings $X_{\mathcal{S}}(\cdot,t)$ and $Y_{\mathcal{S}}(\cdot , t)$, initially defined respectively on $\mathcal{S}$ and $\mathcal{S}(t)$. The process we use is not the same as the one given in \cite{SMSTT}. Instead of extending the Eulerian flow given by the deformation of the solid, we directly extend the deformation of the solid, because the difference in our case lies in the fact that the regularity of the Dirichlet data - written in Eulerian formulation on the time-dependent boundary $\p \mathcal{S}(t)$ - is limited.\\
The goal is to construct a mapping $X$ such that
\begin{eqnarray*}
\left\{ \begin{array} {lll}
\det \nabla X = 1, &  & \text{in } \mathcal{F} \times (0,\infty), \\
X = X_{\mathcal{S}} , &  & \text{on } \p \mathcal{S} \times (0,\infty), \\
X = \Id_{\p \mathcal{O}}, &  & \text{on } \p \mathcal{O} \times (0,\infty).
\end{array} \right.
\end{eqnarray*}

\subsection{Preliminary results}
Let us remind a result stated in the Appendix~B of \cite{Grubb} (Proposition B.1), which treats of Sobolev regularities for products of functions, and that we state as:
\begin{lemma} \label{lemmaGrubb}
Let $s$, $\mu$, and $\kappa$ in $\R$. If $f\in \H^{s+\mu}(\mathcal{F})$ and $g\in \H^{s+\kappa}(\mathcal{F})$, then there exists a positive constant $C$ such that
\begin{eqnarray*}
\| f g \|_{\H^s(\mathcal{F})} & \leq & C \| f  \|_{\H^{s+\mu}(\mathcal{F})} \| g \|_{\H^{s+\kappa}(\mathcal{F})},
\end{eqnarray*}
(i) when $s + \mu + \kappa \geq d/2$, \\
(ii) with $\mu \geq 0$, $\kappa \geq 0$, $2s + \mu + \kappa \geq 0$, \\
(iii) except that $s + \mu + \kappa > d/2$ if equality holds somewhere in (ii).
\end{lemma}

A consequence of this Lemma is the following result.

\begin{lemma} \label{lemmecomatrice}
Let $\overline{X}^{\ast}$ be in $\tilde{\mathcal{W}}(Q_{\infty}^0)$. Then
\begin{eqnarray}
\com \nabla \overline{X}^{\ast} & \in & \L^{\infty}(0,\infty;\mathbf{H}^{2}(\mathcal{F})) \cap \W^{1,\infty}(0,\infty;\mathbf{L}^2(\mathcal{F})), \label{resfond1}
\end{eqnarray}
and, if $\overline{X}^{\ast}-\Id_{\mathcal{F}}$ is small enough in $\tilde{\mathcal{W}}(Q_{\infty}^0)$, there exists a positive constant $C$ such that
\begin{eqnarray}
 \| \com \nabla \overline{X}^{\ast} - \I_{\R^3} \|_{\L^{\infty}(\mathbf{H}^{2}) \cap \W^{1,\infty}(\mathbf{L}^2)}
& \leq & C \| \nabla \overline{X}^{\ast} - \I_{\R^3} \|_{\L^{\infty}(\mathbf{H}^{2}) \cap \W^{1,\infty}(\mathbf{L}^2)}.
 \label{estcof}
\end{eqnarray}
Besides, if $\overline{X}_1^{\ast}-\Id_{\mathcal{F}}$ and $\overline{X}_2^{\ast}-\Id_{\mathcal{F}}$ are small enough in $\tilde{\mathcal{W}}(Q_{\infty}^0)$, there exists a positive constant $C$ such that
\begin{eqnarray}
 \| \com \nabla \overline{X}_2^{\ast} - \com \nabla \overline{X}_1^{\ast} \|_{\L^{\infty}(\mathbf{H}^{2}) \cap \W^{1,\infty}(\mathbf{L}^2)}
& \leq & C \| \nabla \overline{X}_2^{\ast} - \nabla \overline{X}_1^{\ast} \|_{\L^{\infty}(\mathbf{H}^{2}) \cap \W^{1,\infty}(\mathbf{L}^2)}.
\nonumber \\ \label{estcof21}
\end{eqnarray}
\end{lemma}

\begin{proof}
For proving \eqref{resfond1}, the case $d=2$ is obvious. For the general case, let us show that the space $\L^{\infty}(0,\infty;\H^{2}(\mathcal{F})) \cap \W^{1,\infty}(0,\infty;\L^2(\mathcal{F}))$ is stable by product. For that, let us consider two functions $f$ and $g$ which lie in this space. Applying Lemma \ref{lemmaGrubb} with $s = 2$ and $\mu = \kappa = 0$, we get
\begin{eqnarray*}
\| f g \|_{\L^{\infty}(0,\infty;\mathbf{H}^2(\mathcal{F}))}
& \leq & C \| f \|_{\L^{\infty}(0,\infty;\mathbf{H}^2(\mathcal{F}))} \| g \|_{\L^{\infty}(0,\infty;\mathbf{H}^2(\mathcal{F}))}.
\end{eqnarray*}
For the regularity in $\W^{1,\infty}(0,\infty;\mathbf{L}^2(\mathcal{F})$, we write
\begin{eqnarray*}
\frac{\p (fg)}{\p t} & = & \frac{\p f}{\p t}g + f\frac{\p g}{\p t}. \label{eqtriv}
\end{eqnarray*}
Using the continuous embedding $\mathbf{H}^2(\mathcal{F}) \hookrightarrow \mathbf{L}^{\infty}(\mathcal{F})$, we get
\begin{eqnarray*}
\left\|\frac{\p (fg)}{\p t} \right\|_{\L^{\infty}(0,\infty;\L^{2}(\mathcal{F}))} & \leq & C
\left\| \frac{\p f}{\p t} \right\|_{\L^{\infty}(0,\infty;\L^{2}(\mathcal{F}))} \left\| g \right\|_{\L^{\infty}(0,\infty;\H^{2}(\mathcal{F}))} \\
 & & + \left\| \frac{\p g}{\p t} \right\|_{\L^{\infty}(0,\infty;\L^{2}(\mathcal{F}))} \left\| f \right\|_{\L^{\infty}(0,\infty;\H^{2}(\mathcal{F}))},
\end{eqnarray*}
and thus the desired regularity. Thus the space $\L^{\infty}(0,\infty;\mathbf{H}^2(\mathcal{F}))\cap \W^{1,\infty}(0,\infty;\mathbf{L}^2(\mathcal{F}))$ is an algebra. The estimate \eqref{estcof} is obtained by the differentiability of the mapping $\nabla X^{\ast} \mapsto \com \nabla X^{\ast}$ (see \cite{Allaire} for instance); More precisely, we have
\begin{eqnarray*}
\I_{\R^3} - \com \nabla \overline{X}^{\ast} & = &
\left(\nabla \overline{X}^{\ast} - \I_{\R^3}\right)^T - \div\left(\overline{X}^{\ast} - \Id_{\mathcal{F}}\right)
+ o\left(\| \overline{X}^{\ast} - \Id_{\mathcal{F}}\|_{\tilde{\mathcal{W}}(Q_{\infty}^0)} \right),
\end{eqnarray*}
so that we get
\begin{eqnarray*}
\| \com \nabla \overline{X}^{\ast} - \I_{\R^3} \|_{\L^{\infty}(\mathbf{H}^2)\cap \W^{1,\infty}(\mathbf{L}^2)}
& \leq  & C\| \nabla \overline{X}^{\ast} - \I_{\R^3} \|_{\L^{\infty}(\mathbf{H}^2)\cap \W^{1,\infty}(\mathbf{L}^2)}.
\end{eqnarray*}
The estimate \eqref{estcof21} can be obtained by the mean-value theorem, so its proof is left to the reader.
\end{proof}

\subsection{Extension of the Lagrangian mappings}
Let us first extend the deformation of the solid $X^{\ast}$ to the fluid domain $\mathcal{F}$, in a mapping that we have already denoted by $\overline{X}^{\ast}$.

\begin{proposition} \label{lemmaxtension}
Let $X^{\ast} - \Id_{\mathcal{S}} \in \mathcal{W}_{\lambda}(S_{\infty}^0)$ be an {\it admissible} deformation, in the sense of Definition \ref{defcontrol}. Let us assume that $X^{\ast} - \Id_{\mathcal{S}}$ is small enough in $\mathcal{W}_{\lambda}(S_{\infty}^0)$, that is to say that the function
\begin{eqnarray*}
(y,t) \mapsto e^{\lambda t} \frac{\p X^{\ast}}{\p t}
\end{eqnarray*}
is small enough in $\L^{2}(0,\infty;\mathbf{H}^3(\mathcal{S})) \cap \H^{1}(0,\infty;\mathbf{H}^1(\mathcal{S}))$. Then there exists a mapping $\overline{X}^{\ast} \in \tilde{\mathcal{W}}(Q_{\infty}^0)$ satisfying
\begin{eqnarray}
\left\{ \begin{array} {lcl}
\det \nabla \overline{X}^{\ast}  = 1 & \quad & \text{in } \mathcal{F} \times (0,\infty), \label{ext1}\\
\overline{X}^{\ast} = X^{\ast} & \quad & \text{on } \p \mathcal{S} \times (0,\infty), \\
\overline{X}^{\ast} = \Id_{\p \mathcal{O}} & \quad & \text{on } \p \mathcal{O} \times (0,\infty), \label{ext3}
\end{array} \right.
\end{eqnarray}
and such that
\begin{eqnarray}
\| \overline{X}^{\ast}-\Id_{\mathcal{F}} \|_{\tilde{\mathcal{W}}(Q_{\infty}^0)} & \leq & C
\left\| e^{\lambda t}\frac{\p X^{\ast}}{\p t} \right\|_{\L^2(\mathbf{H}^{5/2}(\p \mathcal{S}))\cap \H^1(\mathbf{H}^{1/2}(\p \mathcal{S}))} \label{contSF}
\end{eqnarray}
for some positive constant $C$ independent of $X^{\ast}$. Besides, if $X^{\ast}_1 - \Id_{\mathcal{S}}$ and $X^{\ast}_2 - \Id_{\mathcal{S}}$ are two displacements small enough in $\mathcal{W}_{\lambda}(S_{\infty}^0)$, then the solutions $\overline{X}^{\ast}_1$ and $\overline{X}^{\ast}_2$ of problem \eqref{ext1}, corresponding to $X^{\ast}_1$ and $X^{\ast}_2$ as data respectively, satisfy
\begin{eqnarray}
\| \overline{X}^{\ast}_2-\overline{X}^{\ast}_1  \|_{\tilde{\mathcal{W}}(Q_{\infty}^0)} & \leq & C
\left\| e^{\lambda t}\frac{\p X^{\ast}_2}{\p t} - e^{\lambda t}\frac{\p X^{\ast}_1}{\p t} \right\|_{\L^2(\mathbf{H}^{5/2}(\p \mathcal{S}))\cap \H^1(\mathbf{H}^{1/2}(\p \mathcal{S}))}. \label{contSF21}
\end{eqnarray}
\end{proposition}


\begin{proof}
Given the initial datum $X^{\ast}(y,0) = y$ for $y\in \overline{\mathcal{S}}$, let us consider the system \eqref{ext1} derived in time, as follows
\begin{eqnarray*}
\left\{ \begin{array} {lcl}
\displaystyle \left(\com \nabla \overline{X}^{\ast} \right) : \frac{\p \nabla \overline{X}^{\ast} }{\p t} = 0 &  & \text{in } \mathcal{F} \times (0,\infty), \\
\displaystyle \frac{\p \overline{X}^{\ast}}{\p t} = \frac{\p X^{\ast}}{\p t} &  & \text{on } \p \mathcal{S} \times (0,\infty),
 \\
\displaystyle \frac{\p \overline{X}^{\ast}}{\p t} = 0 &  & \text{on } \p \mathcal{O} \times (0,\infty).
\end{array} \right.
\end{eqnarray*}
This system can be viewed as a modified nonlinear divergence problem, that we state as
\begin{eqnarray}
\left\{ \begin{array} {lcl}
\displaystyle \div \ \frac{\p \overline{X}^{\ast}}{\p t} = f(\overline{X}^{\ast})  &  & \text{in } \mathcal{F} \times (0,\infty), \\
\displaystyle \frac{\p \overline{X}^{\ast}}{\p t} = \frac{\p X^{\ast}}{\p t} &  & \text{on } \p \mathcal{S} \times (0,\infty), \\
\displaystyle \frac{\p \overline{X}^{\ast}}{\p t} = 0 & \quad & \text{on } \p \mathcal{O} \times (0,\infty), \\
\displaystyle \overline{X}^{\ast}(\cdot,0) = \Id_{\mathcal{F}}, \quad \frac{\p \overline{X}^{\ast}}{\p t}(\cdot,0) = 0, & &
\end{array} \right. \label{pbdiv0}
\end{eqnarray}
with
\begin{eqnarray*}
f(\overline{X}^{\ast}) & = & \left(\I_{\R^3} - \com \nabla \overline{X}^{\ast} \right) : \frac{\p \nabla \overline{X}^{\ast} }{\p t}.
\end{eqnarray*}
\textcolor{black}{Let us notice that if we assume in addition that the condition below is satisfied
\begin{eqnarray*}
\int_{\p \mathcal{S}}  \frac{\p \overline{X}^{\ast}}{\p t}\cdot \left(\com  \nabla \overline{X}^{\ast} \right)n \d \Gamma &  = & 0,
\end{eqnarray*}
then from the Piola identity we have
\begin{eqnarray*}
\int_{\mathcal{F}}f(\overline{X}^{\ast}) = \int_{\mathcal{F}} \div\left(
\left(\I_{\R^3} - \com \nabla {\overline{X}^{\ast}}^T\right)\frac{\p \overline{X}^{\ast}}{\p t}
\right)
 = \int_{\p \mathcal{S}} \frac{\p \overline{X}^{\ast}}{\p t}\cdot n \d \Gamma = 
 \int_{\p \mathcal{S}} \frac{\p X^{\ast}}{\p t}\cdot n \d \Gamma
\end{eqnarray*}
and thus the compatibility condition for this divergence system is satisfied.}\\
A solution of this system can be viewed as a fixed point of the mapping
\begin{eqnarray}
\begin{array} {cccc} 
\mathfrak{T} : & \text{\textcolor{black}{$\mathfrak{W}_\lambda$}} & \rightarrow & \text{\textcolor{black}{$\mathfrak{W}_\lambda$}} \\
& \overline{X}^{\ast}_1-\Id_{\mathcal{F}} & \mapsto & \overline{X}^{\ast}_2 - \Id_{\mathcal{F}},
\end{array}
\end{eqnarray}
\textcolor{black}{where\footnote{\textcolor{black}{The set $\mathfrak{W}_{\lambda}$ is non-trivial; Indeed, it contains extensions of $X^{\ast}-\Id_{\mathcal{S}}$ obtained by the use of {\it plateau fonctions} - see \cite{SMSTT} for instance. The difficulty here is to obtain Lipschitz estimates on the Lagrangian mappings.}}
\begin{eqnarray*}
\text{\textcolor{black}{$\mathfrak{W}_{\lambda}$}}  =  \text{\textcolor{black}{
$\displaystyle \left\{\overline{Z}^{\ast} \in\mathcal{W}_{\lambda}(Q_{\infty}^0) \mid  \overline{Z}^{\ast} = X^{\ast} - \Id_{\mathcal{S}} \text{ on } \p \mathcal{S}, 
\ \overline{Z}^{\ast}(\cdot,0) = 0,
\right.$}} & & \\
\text{\textcolor{black}{$ \displaystyle \left.
\frac{\p \overline{Z}^{\ast}}{\p t}(\cdot,0) = 0, \
\int_{\p \mathcal{S}}  \frac{\p \overline{Z}^{\ast}}{\p t}\cdot \left(\com  (\nabla \overline{Z}^{\ast} + \I_{\R^3})  \right)n \d \Gamma  = 0 
\right\}$}} & & 
\end{eqnarray*}
and} where $\overline{X}^{\ast}_2$ satisfies the classical divergence problem
\begin{eqnarray*}
\left\{ \begin{array} {lcl}
\displaystyle \div \ \frac{\p \overline{X}^{\ast}_2}{\p t} = f(\overline{X}^{\ast}_1)  & \quad & \text{in } \mathcal{F} \times (0,\infty), \\
\displaystyle \frac{\p \overline{X}^{\ast}_2}{\p t} = \frac{\p X^{\ast}}{\p t} & \quad & \text{on } \p \mathcal{S} \times (0,\infty), \\
\displaystyle \frac{\p \overline{X}^{\ast}_2}{\p t} = 0 & \quad & \text{on } \p \mathcal{O} \times (0,\infty), \\
\displaystyle \overline{X}^{\ast}_2(\cdot,0) = \Id_{\mathcal{F}}, \quad \frac{\p \overline{X}^{\ast}_2}{\p t}(\cdot,0) = 0. & &
\end{array} \right. \label{pbdiv}
\end{eqnarray*}
Indeed, let us first verify that for $\overline{X}^{\ast}-\Id_{\mathcal{F}} \in \mathcal{W}_{\lambda}(Q_{\infty}^0)$ we have $e^{\lambda t}f(\overline{X}^{\ast}) \in \L^{2}(0,\infty;\mathbf{H}^{2}(\mathcal{F})) \cap \H^{1}(0,\infty;\mathbf{L}^2(\mathcal{F}))$. For that, we remind from the previous lemma that $\com \nabla \overline{X}^{\ast} \in \L^{\infty}(0,\infty;\mathbf{H}^{2}(\mathcal{F})) \cap \W^{1,\infty}(0,\infty;\mathbf{L}^2(\mathcal{F}))$, and we first use the result of Lemma \ref{lemmaGrubb} with $s = 2$ and $\mu = \kappa = 0$ to get
\begin{eqnarray*}
\|e^{\lambda t} f(\overline{X}^{\ast}) \|_{\L^2(\mathbf{H}^{2}(\mathcal{F}))} & \leq & C \| \I_{\R^3} - \com \nabla \overline{X}^{\ast}  \|_{\L^{\infty}(\mathbf{H}^{2}(\mathcal{F}))}
\left\|e^{\lambda t}\frac{\p \nabla \overline{X}^{\ast} }{\p t}\right\|_{\L^2(\mathbf{H}^{2}(\mathcal{F}))}.
\end{eqnarray*}
For the regularity in $\H^{1}(0,\infty;\mathbf{L}^2(\mathcal{F}))$, we write
\begin{eqnarray*}
e^{\lambda t} \frac{\p f(\overline{X}^{\ast})}{\p t} = (\I_{\R^3} - \com \nabla \overline{X}^{\ast}):
\left( e^{\lambda t} \frac{\p^2 \nabla \overline{X}^{\ast} }{\p t^2}\right) - \frac{\p \com \nabla \overline{X}^{\ast}}{\p t}:
\left( e^{\lambda t} \frac{\p \nabla \overline{X}^{\ast} }{\p t} \right),
\end{eqnarray*}
\begin{eqnarray*}
\left\| e^{\lambda t}\frac{\p f(\overline{X}^{\ast})}{\p t} \right\|_{\L^2(\mathbf{L}^2(\mathcal{F}))} & \leq &
C\| \I_{\R^3} - \com \nabla \overline{X}^{\ast} \|_{\L^{\infty}(\mathbf{H}^{2}(\mathcal{F}))}
\left\| e^{\lambda t} \frac{\p \nabla \overline{X}^{\ast} }{\p t} \right\|_{\H^1(\mathbf{L}^2(\mathcal{F}))} \\
& & + C\left\| \frac{\p  \com \nabla \overline{X}^{\ast}}{\p t} \right\|_{\L^{\infty}(\mathbf{L}^{2}(\mathcal{F}))}
\left\| e^{\lambda t} \frac{\p \nabla \overline{X}^{\ast} }{\p t} \right\|_{\L^2(\mathbf{H}^{2}(\mathcal{F}))}
\end{eqnarray*}
where we have used the continuous embedding $\mathbf{H}^2(\mathcal{F}) \hookrightarrow \mathbf{L}^{\infty}(\mathcal{F})$. Thus there exists a positive constant $C$ such that
\begin{eqnarray}
\left\|e^{\lambda t}f(\overline{X}^{\ast})\right\|_{\H^{2,1}(Q_{\infty}^0)} & \leq &
C\| \I_{\R^3} - \com \nabla \overline{X}^{\ast} \|_{\L^{\infty}(\mathbf{H}^{2}(\mathcal{F}))\cap \W^{1,\infty}(\mathbf{L}^{2}(\mathcal{F}))}
\left\|e^{\lambda t}\frac{\p \nabla \overline{X}^{\ast} }{\p t}\right\|_{\H^{2,1}(Q_{\infty}^0)}. \nonumber \\ \label{ineq001}
\end{eqnarray}
The estimate \eqref{ineq001} shows in particular that the mapping $\mathfrak{T}$ is well-defined. Moreover, for the divergence problem \eqref{pbdiv} there exists a positive constant $C$ (see \cite{Galdi1} for instance\footnote{Using results of \cite{Galdi1}, the nonhomogeneous Dirichlet condition can be lifted (see Theorem 3.4, Chapter II) and the resolution made by using Exercise 3.4 and Theorem 3.2 of Chapter III.}) such that
\begin{eqnarray*}
 \left\|e^{\lambda t} \frac{\p \overline{X}^{\ast}_2}{\p t} \right\|_{\L^{2}(\mathbf{H}^{3}(\mathcal{F}))} &  \leq &
C \left(\| e^{\lambda t}f(\overline{X}^{\ast}_1) \|_{\L^2(\mathbf{H}^2(\mathcal{F}))}
+ \left\|e^{\lambda t} \frac{\p X^{\ast}}{\p t} \right\|_{\L^2(\mathbf{H}^{5/2}(\p \mathcal{S}))} \right), \label{newineq}
\end{eqnarray*}
and also
\begin{eqnarray*}
  \left\|e^{\lambda t} \frac{\p \overline{X}^{\ast}_2}{\p t} \right\|_{\H^{1}(\mathbf{H}^{1}(\mathcal{F}))}  & \leq &
C \left((1+\lambda)\|e^{\lambda t} f(\overline{X}^{\ast}_1) \|_{\H^1(\mathbf{L}^2(\mathcal{F}))}
+ \left\|e^{\lambda t} \frac{\p X^{\ast}}{\p t} \right\|_{\H^1(\mathbf{H}^{1/2}(\p \mathcal{S}))} \right). \label{newineqbis}
\end{eqnarray*}
Thus there exists a positive constant $C_0$ such that
\begin{eqnarray}
\left\|e^{\lambda t} \frac{\p \overline{X}^{\ast}_2}{\p t} \right\|_{\mathcal{H}(Q_{\infty}^0) }  & \leq &
 C_0\left(\| e^{\lambda t}f(\overline{X}^{\ast}_1) \|_{\H^{2,1}(Q^0_{\infty})}
+ \left\|e^{\lambda t} \frac{\p X^{\ast}}{\p t} \right\|_{\mathcal{H}(S_{\infty}^0)}\right). \label{newineqter}
\end{eqnarray}
Let us consider the set
\begin{eqnarray*}
\mathfrak{B}_R & = & \left\{   \overline{Z}^{\ast}  \in \text{\textcolor{black}{$\mathfrak{W}_{\lambda}$}}, \
\|  \overline{Z}^{\ast} \|_{\mathcal{W}_{\lambda}(Q_{\infty}^0)} \leq R \right\}
\end{eqnarray*}
with
\begin{eqnarray*}
R & = & 2C_0\left\|e^{\lambda t} \frac{\p X^{\ast}}{\p t} \right\|_{\mathcal{H}(S_{\infty}^0)}.
\end{eqnarray*}
Notice that a mapping $\overline{X}^{\ast}_1 -\Id_{\mathcal{F}} \in \mathfrak{B}_R$ satisfies in particular the following inequality, obtained in the same way we have proceeded to get the embedding \eqref{ineqstar}:
\begin{eqnarray*}
\| \overline{X}^{\ast}_1 - \Id_{\mathcal{F}} \|_{\tilde{\mathcal{W}}(Q_{\infty}^0)} & \leq &
C \left\|e^{\lambda t} \frac{\p \overline{X}^{\ast}_1}{\p t} \right\|_{\mathcal{H}(Q_{\infty}^0)}
= C \| \overline{X}^{\ast}_1 - \Id_{\mathcal{F}} \|_{\mathcal{W}_{\lambda}(Q_{\infty}^0)} \\
& \leq & C R .
\end{eqnarray*}
Then the inequality \eqref{newineqter} combined to the estimates \eqref{ineq001} and \eqref{estcof} show that for $\overline{X}^{\ast}_1 - \Id_{\mathcal{F}} \in \mathfrak{B}_R$ we have
\begin{eqnarray*}
\| \overline{X}^{\ast}_2 - \Id_{\mathcal{F}} \|_{\mathcal{W}_{\lambda}(Q_{\infty}^0)} & \leq &
 C_0\left( CR^2(CR+1)+ \frac{R}{2C_0} \right),
\end{eqnarray*}
and thus for $R$ small enough, $\mathfrak{B}_R$ is stable by $\mathfrak{T}$. Notice that $\mathfrak{B}_R$ is a closed subset of $\tilde{\mathcal{W}}(Q_{\infty}^0)$. Let us verify that $\mathfrak{T}$ is a contraction in $\mathfrak{B}_R$.\\
For $\overline{X}^{\ast}_1-\Id_{\mathcal{F}}$ and $\overline{X}^{\ast}_2-\Id_{\mathcal{F}}$ in $\mathfrak{B}_R$, we denote $\overline{Z}^{\ast} = \mathfrak{T}(\overline{X}^{\ast}_2-\Id_{\mathcal{F}}) - \mathfrak{T}(\overline{X}^{\ast}_1-\Id_{\mathcal{F}})$ which satisfies the divergence system
\begin{eqnarray*}
\left\{ \begin{array} {lcl}
\displaystyle \div \ \frac{\p \overline{Z}^{\ast}}{\p t} = f(\overline{X}^{\ast}_2) - f(\overline{X}^{\ast}_1)  & \quad & \text{in } \mathcal{F} \times (0,\infty), \\
\displaystyle \frac{\p \overline{Z}^{\ast}}{\p t} = 0 & \quad & \text{on } \p \mathcal{S} \times (0,\infty), \\
\displaystyle \frac{\p \overline{Z}^{\ast}}{\p t} = 0 & \quad & \text{on } \p \mathcal{O} \times (0,\infty).
\end{array} \right.
\end{eqnarray*}
and thus the estimate
\begin{eqnarray*}
 \left\|e^{\lambda t} \frac{\p \overline{Z}^{\ast}}{\p t} \right\|_{\mathcal{H}(Q_{\infty}^0)} \leq
  C_{\mathcal{F}} \left\|e^{\lambda t}\left( f(\overline{X}^{\ast}_2) - f(\tilde{X}_1)\right) \right\|_{\H^{2,1}(Q_{\infty}^0)}.
\end{eqnarray*}
For tackling the Lipschitz property of the nonlinearity, we write
\begin{eqnarray*}
f(\overline{X}^{\ast}_2) - f(\overline{X}^{\ast}_1) & = & \left(\com \nabla \overline{X}^{\ast}_2 - \com \nabla \overline{X}^{\ast}_1 \right) : \frac{\p \nabla \overline{X}^{\ast}_2}{\p t} \\
& &  + \left(\I_{\R^3} - \com \nabla \overline{X}^{\ast}_1 \right) : \frac{\p \nabla (\overline{X}^{\ast}_2 - \overline{X}^{\ast}_1)}{\p t}.
\end{eqnarray*}
By reconsidering the steps of the proof of the estimate \eqref{ineq001} and by using \eqref{estcof21}, we can verify that for $R$ small enough the mapping $\mathfrak{T}$ is a contraction in $\mathfrak{B}_R$. Thus $\mathfrak{T}$ admits a unique fixed point in $\mathfrak{B}_R$.\\
For the estimate \eqref{contSF21}, if $\overline{X}^{\ast}_1$ and $\overline{X}^{\ast}_2$ are two solutions corresponding to $X^{\ast}_1$ and $X^{\ast}_2$ respectively, let us just write the system satisfied by the difference $\overline{Z}^{\ast} = \overline{X}^{\ast}_2- \overline{X}^{\ast}_1$:
\begin{eqnarray*}
\left\{ \begin{array} {lcl}
\displaystyle \div \ \frac{\p \overline{Z}^{\ast}}{\p t} = f(\overline{X}^{\ast}_2) - f(\overline{X}^{\ast}_1)  &  & \text{in } \mathcal{F} \times (0,\infty), \\
\displaystyle \frac{\p \overline{Z}^{\ast}}{\p t} = \frac{\p X^{\ast}_2}{\p t} - \frac{\p X^{\ast}_1}{\p t} &  & \text{on } \p \mathcal{S} \times (0,\infty), \\
\displaystyle \frac{\p \overline{Z}^{\ast}}{\p t} = 0 & & \text{on } \p \mathcal{O} \times (0,\infty).
\end{array} \right.
\end{eqnarray*}
Then the methods used above can be similarly applied to this system in order to deduce from it the announced result.
\end{proof}
\hfill \\

Let us now consider $h \in \H^2(0,\infty;\R^3)$, and $\mathbf{R} \in \H^2(0,\infty;\R^9)$ which provides $\omega \in \H^1(0,\infty;\R^3)$. Let us construct a mapping $X$ such that $X(\cdot,0) = \Id_{\mathcal{F}}$ and
\begin{eqnarray*}
\left\{ \begin{array} {lll}
\det \nabla X = 1 &  & \text{in } \mathcal{F} \times (0,\infty), \\
X = h(t) + \mathbf{R}(t)X^{\ast}  &  & \text{on } \p \mathcal{S} \times (0,\infty), \\
X = \Id_{\p \mathcal{O}} & & \text{on } \p \mathcal{S} \times (0,\infty).
\end{array} \right.
\end{eqnarray*}
We cannot solve this problem as we have done for problem \eqref{ext1}, because the proof would require the unknowns $h$ and $\omega$ arbitrarily small enough, a thing that we cannot assume, even {\it a posteriori}. Instead of that, we utilize the mapping $\overline{X}^{\ast}$ provided by Proposition \ref{lemmaxtension}, and we search for a mapping $\overline{X}^R$ such that
\begin{eqnarray*}
X(\cdot,t) &  =  & \overline{X}^R(\cdot,t)  \circ  \overline{X}^{\ast}(\cdot,t).
\end{eqnarray*}
Such a mapping $\overline{X}^R$ has to satisfy
\begin{eqnarray*}
\left\{ \begin{array} {lll}
\det \nabla \overline{X}^{R} = 1 &  & \text{in } \mathcal{F} \times (0,\infty), \\
\overline{X}^{R} = h + \mathbf{R} \Id_{\p \mathcal{S}} &  & \text{on } \p \mathcal{S} \times (0,\infty), \\
\overline{X}^{R} = \Id_{\p \mathcal{O}} &  & \text{on } \p \mathcal{O}\times (0,\infty).
\end{array} \right.
\end{eqnarray*}
For that, let us proceed as in \cite{TT}: We consider a cut-off function $\xi \in C^{\infty}(\mathcal{F})$, such that $\xi \equiv 1$ in a vicinity of $\p \mathcal{S}$ and $\xi \equiv 0$ in a vicinity of $\p \mathcal{O}$. We define the function
\begin{eqnarray*}
 \mathfrak{F}_R(x,t)  = \frac{1}{2}h'(t) \wedge (x-h(t))-\frac{1}{2}|x-h(t)|^2\omega(t),
\end{eqnarray*}
so that $\rot (\mathfrak{F}_R)(x,t) = h'(t) + \omega(t) \wedge (x-h(t))$, and we construct $\overline{X}^R$ as the solution of the following Cauchy problem
\begin{eqnarray}
\frac{\p \overline{X}^R}{\p t}(\overline{x}^{\ast},t) = \rot(\xi \mathfrak{F}_R)(\overline{X}^R(\overline{x}^{\ast},t),t),
\quad \overline{X}^R(\overline{x}^{\ast},0) = \overline{x}^{\ast},
\quad \overline{x}^{\ast} \in \overline{X}^{\ast}(\mathcal{F},t) = \mathcal{F}^{\ast}. \nonumber \\ \label{pbextR}
\end{eqnarray}
We can verify (see \cite{TT} for instance) that the mapping $\overline{X}^R$ so obtained has the desired properties, and thus we can set
\begin{eqnarray}
X(y,t) &  =  & \overline{X}^R(\overline{X}^{\ast}(y,t),t), \quad (y,t) \in \mathcal{F} \times (0,\infty). \label{eqextension}
\end{eqnarray}
Since $\overline{X}^{\ast}(\cdot,t)$ and $\overline{X}^R(\cdot,t)$ are invertible, the mapping $X(\cdot,t)$ is invertible, and we denote by $Y(\cdot,t)$ its inverse. The mapping $X$ presents the same type of regularity as the mapping $\overline{X}^{\ast}$ . We sum its properties in the following proposition.\\

\begin{proposition} \label{restildeX} \label{resX}
Let $X^{\ast}$ be an {\it admissible} control - in the sense of Definition \ref{defcontrol} - and $\overline{X}^{\ast}$ the extension of $X^{\ast}$ provided by Proposition \ref{lemmaxtension} (for $ e^{\lambda t}\frac{\p X^{\ast}}{\p t}$ small enough in $\mathcal{H}_{3}(S_{\infty}^0)$). Let $X$ be the mapping given by \eqref{eqextension}. For all $t \geq 0$, the mapping $X(\cdot,t)$ is a $C^1$-diffeomorphism from $\mathcal{O}$ onto $\mathcal{O}$, from $\p \mathcal{S}$ onto $\p \mathcal{S}(t)$, and from $\mathcal{F}$ onto $\mathcal{F}(t)$. We denote by $Y(\cdot,t)$ its inverse at some time $t$. We have
\begin{eqnarray*}
& & (y,t) \mapsto X(y,t)  \in  \tilde{\mathcal{W}}(Q_{\infty}^0), \\
& & \det \nabla X(y,t)  =   1, \ \text{for all } (y,t) \in \mathcal{F} \times (0,\infty).
\end{eqnarray*}
\end{proposition}

The proof for the regularity of $X$ can be straightforwardly deduced from Lemma \ref{lemmaX21} in the Appendix~B of this chapter. We do not give more detail in this section, because here the aim is only to get a change of variables which enables us rewrite the main system as an equivalent one written in fixed domains (see section \ref{secchange}).

\section{Appendix~B: Proofs of estimates for the changes of variables}
Let us remind that for $X^{\ast}-\Id_{\mathcal{S}} \in \mathcal{W}_{\lambda}(S_{\infty}^0)$ Proposition \ref{lemmaxtension} enables us to define the extension $\overline{X}^{\ast} \in \mathcal{W}_{\lambda}(Q_{\infty}^0) + \Id_{\mathcal{F}}$ satisfying
\begin{eqnarray*}
\left\{ \begin{array} {lcl}
\det \nabla \overline{X}^{\ast}  = 1 & & \text{in } \mathcal{F} \times (0,\infty), \\
\overline{X}^{\ast} = X^{\ast} & & \text{on } \p \mathcal{S} \times (0,\infty), \\
\overline{X}^{\ast} = \Id_{\p \mathcal{O}} & & \text{on } \p \mathcal{O} \times (0,\infty),
\end{array} \right.
\end{eqnarray*}
and
\begin{eqnarray*}
\| \overline{X}^{\ast}-\Id_{\mathcal{F}} \|_{\tilde{\mathcal{W}}(Q_{\infty}^0)} & \leq & C
\left\| e^{\lambda t}\frac{\p X^{\ast}}{\p t} \right\|_{\mathcal{H}(Q_{\infty}^0)}.
\end{eqnarray*}
For $h\in \H^2(0,\infty;\R^3)$ and $\mathbf{R} \in \H^2(0,\infty;\R^9)$ which provides $\omega \in \H^1(0,\infty;\R^3)$ such that
\begin{eqnarray*}
\left\{\begin{array} {ccccc}
\displaystyle \frac{\d \mathbf{R}}{\d t} & = & \mathbb{S}\left( \omega\right) \mathbf{R} \\
\mathbf{R}(0) & = & \I_{\R^3},
\end{array}
\right.
\quad \text{with }
\mathbb{S}(\omega) = \left(
\begin{matrix}
0 & -\omega_3 & \omega_2 \\
\omega_3 & 0 & -\omega_1 \\
-\omega_2 & \omega_1 & 0
\end{matrix} \right),
\end{eqnarray*}
we can define $\overline{X}^R$ through the problem
\begin{eqnarray*}
\frac{\p \overline{X}^R}{\p t}(\overline{x}^{\ast},t) = \rot(\xi \mathfrak{F}_R)(\overline{X}^R(\overline{x}^{\ast},t),t),
\quad \overline{X}^R(\overline{x}^{\ast},0) = \overline{x}^{\ast},
\quad \overline{x}^{\ast} \in \overline{X}^{\ast}(\mathcal{F},t) = \mathcal{F}^{\ast}(t),
\end{eqnarray*}
where $\xi$ is a regular cut-off function, and
\begin{eqnarray*}
 \mathfrak{F}_R(x,t)  = \frac{1}{2}h'(t) \wedge (x-h(t))-\frac{1}{2}|x-h(t)|^2\omega(t).
\end{eqnarray*}
Then we define
\begin{eqnarray*}
X & = & \overline{X}^R \circ \overline{X}^{\ast},
\end{eqnarray*}
and
\begin{eqnarray*}
\tilde{X} & = & \mathbf{R}^T(X-h).
\end{eqnarray*}

\begin{lemma} \label{lemmaKstar}
For $X^{\ast}_1 - \Id_{\mathcal{S}}$ and $X^{\ast}_2 - \Id_{\mathcal{S}}$ small enough in $\mathcal{W}_{\lambda}(S_{\infty}^0)$, let $\overline{X}^{\ast}_1$ and $\overline{X}^{\ast}_2$ be the solutions of Problem \eqref{ext1} (see Proposition \ref{lemmaxtension}) corresponding to the data $X^{\ast}_1$ and $X^{\ast}_2$ respectively. If we denote by $\overline{Y}^{\ast}_1(\cdot,t)$ and $\overline{Y}^{\ast}_2(\cdot,t)$ the inverses of $\overline{X}^{\ast}_1$ and $\overline{X}^{\ast}_2$ respectively, we have
\begin{eqnarray}
& & \| \nabla \overline{Y}^{\ast}_2(\overline{X}^{\ast}_2) - \nabla \overline{Y}^{\ast}_1(\overline{X}^{\ast}_1) \|_{\L^{\infty}(0,\infty;\mathbf{H}^2(\mathcal{F}))\cap \W^{1,\infty}(0,\infty;\mathbf{L}^2(\mathcal{F}))}  \nonumber \\
& & \leq C \left\| e^{\lambda t}\frac{\p X^{\ast}_2}{\p t}-e^{\lambda t}\frac{\p X^{\ast}_1}{\p t} \right\|_{\mathcal{H}(S_{\infty}^0)}.
 \label{estnablatilde}
\end{eqnarray}
\end{lemma}

\begin{proof}
Let us remind the estimate \eqref{contSF21}, obtained for $X^{\ast}_1 - \Id_{\mathcal{S}}$ and $X^{\ast}_2 - \Id_{\mathcal{S}}$ small enough in $\mathcal{W}_{\lambda}(S_{\infty}^0)$:
\begin{eqnarray*}
\| \overline{X}^{\ast}_2-\overline{X}^{\ast}_1  \|_{\tilde{\mathcal{W}}(Q_{\infty}^0)} & \leq & C
\left\| e^{\lambda t}\frac{\p X^{\ast}_2}{\p t} - e^{\lambda t}\frac{\p X^{\ast}_1}{\p t} \right\|_{\mathcal{H}(S_{\infty}^0)}.
\end{eqnarray*}
First, let us give two intermediate estimates:
\begin{eqnarray*}
\| \nabla \overline{X}^{\ast}_1 - \I_{\R^3} \|_{\L^{\infty}(\mathbf{H}^2)\cap\W^{1,\infty}(\mathbf{L}^2)} & \leq & \left\| e^{\lambda t}\frac{\p X^{\ast}_1}{\p t} \right\|_{\mathcal{H}(S_{\infty}^0)}
\end{eqnarray*}
and
\begin{eqnarray*}
\| \nabla \overline{Y}^{\ast}_2(\overline{X}^{\ast}_2) - \I_{\R^3} \|_{\L^{\infty}(\mathbf{H}^2)\cap \W^{1,\infty}(\mathbf{L}^2)} & \leq &
\frac{ \| \nabla \overline{X}^{\ast}_2 - \I_{\R^3} \|_{\L^{\infty}(\mathbf{H}^2)\cap \W^{1,\infty}(\mathbf{L}^2)} }{1 - C \| \nabla \overline{X}^{\ast}_2 - \I_{\R^3} \|_{\L^{\infty}(\mathbf{H}^2)\cap \W^{1,\infty}(\mathbf{L}^2)} } \\
& \leq & 2  \| \nabla \overline{X}^{\ast}_2 - \I_{\R^3} \|_{\L^{\infty}(\mathbf{H}^2)\cap \W^{1,\infty}(\mathbf{L}^2)},
\end{eqnarray*}
provided by the equality
\begin{eqnarray*}
\nabla \overline{Y}^{\ast}_2(\overline{X}^{\ast}_2) - \I_{\R^3}  =
\left(\I_{\R^3} - \nabla \overline{X}^{\ast}_2 \right)\left(\nabla \overline{Y}^{\ast}_2(\overline{X}^{\ast}_2) - \I_{\R^3} \right) + \left(\I_{\R^3} - \nabla \overline{X}^{\ast}_2 \right)
\end{eqnarray*}
and by the fact that $\L^{\infty}(0,\infty;\mathbf{H}^2(\mathcal{F}))\cap \W^{1,\infty}(0,\infty;\mathbf{L}^2(\mathcal{F}))$ is an algebra. Then the estimate \eqref{estnablatilde} is obtained by writing
\begin{eqnarray*}
\nabla \overline{Y}^{\ast}_2(\overline{X}^{\ast}_2) - \nabla \overline{Y}^{\ast}_1(\overline{X}^{\ast}_1) & = &
\left( \nabla \overline{Y}^{\ast}_2(\overline{X}^{\ast}_2) - \nabla \overline{Y}^{\ast}_1(\overline{X}^{\ast}_1) \right)
\left(\I_{\R^3} - \nabla \overline{X}^{\ast}_1\right) \\
& &  - \left(\nabla \overline{X}^{\ast}_2 - \nabla \overline{X}^{\ast}_1 \right) \nabla \overline{Y}^{\ast}_2(\overline{X}^{\ast}_2), \\
\| \nabla \overline{Y}^{\ast}_2(\overline{X}^{\ast}_2) - \nabla \overline{Y}^{\ast}_1(\overline{X}^{\ast}_1)\| & \leq &
\frac{ \left(1+ \| \nabla \overline{Y}^{\ast}_2(\overline{X}^{\ast}_2) - \I_{\R^3} \|\right)\|
\nabla \overline{X}^{\ast}_2 - \nabla \overline{X}^{\ast}_1 \|}
{1-C\| \nabla \overline{X}^{\ast}_1 - \I_{\R^3} \|}.
\end{eqnarray*}
\end{proof}

\begin{lemma}
Let $\overline{X}^R_1$ and $\overline{X}^R_2$ be the extensions defined by problem \eqref{pbextR}, with data $(h_1,\mathbf{R}_1) \in \H^2(0,\infty;\R^3)\times \H^2(0,\infty;\R^9)$ and $(h_2,\mathbf{R}_2) \in \H^2(0,\infty;\R^3)\times \H^2(0,\infty;\R^9)$ respectively. Then we have
\begin{eqnarray*}
\| \overline{X}^R_2 - \overline{X}^R_1 \|_{\H^2(0,\infty;W^{4,\infty}(\R^3))} & \leq &
O\left(\|\hat{h}'_2-\hat{h}'_1\|_{\H^1(0,\infty;\R^3)}+ \|\hat{\omega}_2 - \hat{\omega}_1\|_{\H^1(0,\infty;\R^3)}  \right),
\end{eqnarray*}
where we remind that
\begin{eqnarray*}
\begin{array} {lll}
\hat{h}'_1 = e^{\lambda t}\mathbf{R}_1^T h'_1, & \quad & \hat{h}'_2 = e^{\lambda t}\mathbf{R}_2^T h'_2, \\
\hat{\omega}_1 = e^{\lambda t}\mathbf{R}_1^T \omega_1, & \quad & \hat{\omega}_2 = e^{\lambda t}\mathbf{R}_2^T \omega_2.
\end{array}
\end{eqnarray*}
\end{lemma}

\begin{proof}
The change of variables given by a mapping $\overline{X}^R$ is slightly the same as the one utilized in \cite{TT}; In considering the writing
\begin{eqnarray*}
\mathbf{R}(t)^T \mathfrak{F}_R(x,t)  = \frac{1}{2}\tilde{h}'(t) \wedge \mathbf{R}(t)^T(x-h(t))
-\frac{1}{2}|x-h(t)|^2\tilde{\omega}(t),
\end{eqnarray*}
the steps of the proofs of Lemmas 6.11 and 6.12 of \cite{TT} can be then repeated, with the difference that in infinite time horizon we rather have
\begin{eqnarray*}
\| \overline{X}^R_2 - \overline{X}^R_1 \|_{\H^2(0,\infty;W^{4,\infty}(\R^3))} & \leq &
K^R\left(\|\tilde{h}'_2-\tilde{h}'_1\|_{\H^1(0,\infty;\R^3)}+ \|\tilde{\omega}_2 - \tilde{\omega}_1\|_{\H^1(0,\infty;\R^3)} \right. \\
& & \left. + \|h_2-h_1\|_{\L^{\infty}(0,\infty;\R^3)} + \| \mathbf{R}_2 - \mathbf{R}_1 \|_{\L^{\infty}(0,\infty;\R^9)} \right),
\end{eqnarray*}
where $K^R$ is bounded when $h_1, \ h_2$ are close to $0$ and $\mathbf{R}_1, \ \mathbf{R}_2$ are close to $\I_{\R^3}$. In order to estimate $\|h_2-h_1\|_{\L^{\infty}(0,\infty;\R^3)}$ and $\| \mathbf{R}_2 - \mathbf{R}_1 \|_{\L^{\infty}(0,\infty;\R^9)}$, we first apply the Gr\"onwall's lemma on
\begin{eqnarray*}
\left\{\begin{array} {rcl}
\displaystyle \frac{\p }{\p t}\left( \mathbf{R}_2 - \mathbf{R}_1 \right) & = &
(\mathbf{R}_2 - \mathbf{R}_1)\mathbb{S}(\tilde{\omega}_2) + \mathbf{R}_1 \mathbb{S}(\tilde{\omega}_2 - \tilde{\omega}_1) \\
\displaystyle (\mathbf{R}_2 - \mathbf{R}_1)(0) & = & 0
\end{array} \right.
\end{eqnarray*}
in order to get
\begin{eqnarray*}
\| \mathbf{R}_2 - \mathbf{R}_1 \|_{\L^{\infty}(0,\infty;\R^9)} & \leq &
C \|\tilde{\omega}_2 - \tilde{\omega}_1\|_{\L^1(0,\infty;\R^3)} \exp \left( C \|\tilde{\omega}_2 \|_{\L^1(0,\infty;\R^3)} \right).
\end{eqnarray*}
Besides, it is easy to see that
\begin{eqnarray*}
\|\tilde{\omega}_2 - \tilde{\omega}_1\|_{\L^1(0,\infty;\R^3)} & \leq &
\frac{1}{\sqrt{2\lambda}} \|\hat{\omega}_2 - \hat{\omega}_1\|_{\L^2(0,\infty;\R^3)},
\end{eqnarray*}
so that $\| \mathbf{R}_2 - \mathbf{R}_1 \|_{\L^{\infty}(0,\infty;\R^9)}$ is controlled by $\|\hat{\omega}_2 - \hat{\omega}_1\|_{\L^2(0,\infty;\R^3)}$. Then the term $\|h_2-h_1\|_{\L^{\infty}(0,\infty;\R^3)}$ can be treated by writing
\begin{eqnarray*}
\|h_2 - h_1 \|_{\L^{\infty}(0,\infty;\R^3)} & \leq & \| h'_2 - h'_1 \|_{\L^{1}(0,\infty;\R^3)}, \\
h'_2 - h'_1 & = & \left(\mathbf{R}_2 - \mathbf{R}_1 \right) \tilde{h}'_2 + \mathbf{R}_1\left(\tilde{h}'_2 - \tilde{h}'_1 \right), \\
\| h'_2 - h'_1 \|_{\L^{1}(0,\infty;\R^3)} & \leq &
\| \mathbf{R}_2 - \mathbf{R}_1 \|_{\L^{\infty}(0,\infty;\R^9)}\| \tilde{h}'_2 \|_{\L^{1}(0,\infty;\R^3)}
 + \| \tilde{h}'_2 - \tilde{h}'_1 \|_{\L^{1}(0,\infty;\R^3)} \\
\| \tilde{h}'_2  \|_{\L^{1}(0,\infty;\R^3)} & \leq & \frac{1}{\sqrt{2\lambda}} \| \hat{h}'_2 \|_{\L^{2}(0,\infty;\R^3)}, \\
\| \tilde{h}'_2 - \tilde{h}'_1 \|_{\L^{1}(0,\infty;\R^3)} & \leq &
\frac{1}{\sqrt{2\lambda}} \| \hat{h}'_2 - \hat{h}'_1 \|_{\L^{2}(0,\infty;\R^3)}.
\end{eqnarray*}
Finally, it is easy to verify that
\begin{eqnarray*}
\|\tilde{h}'_2-\tilde{h}'_1\|_{\H^1(0,\infty;\R^3)} & \leq & (1+ \lambda)\|\hat{h}'_2-\hat{h}'_1\|_{\H^1(0,\infty;\R^3)}, \\
\|\tilde{\omega}_2-\tilde{\omega}_1\|_{\H^1(0,\infty;\R^3)} & \leq & (1+ \lambda)\|\hat{\omega}_2-\hat{\omega}_1\|_{\H^1(0,\infty;\R^3)}.
\end{eqnarray*}
\end{proof}

\begin{lemma} \label{lemmaX21}
Let $X_1$ and $X_2$ be defined by
\begin{eqnarray*}
X_1 = \overline{X}^R_1 \circ \overline{X}^{\ast}_1, & \quad & X_2 = \overline{X}^R_2 \circ \overline{X}^{\ast}_2,
\end{eqnarray*}
where $\overline{X}^{\ast}_1$, $\overline{X}^{\ast}_2$, $\overline{X}^R_1$ and $\overline{X}^R_2$ are given in the assumptions of the previous lemmas. Then we have
\begin{eqnarray*}
\|X_2 - X_1 \|_{\tilde{\mathcal{W}}(Q_{\infty}^0)} & = & r K(r),
\end{eqnarray*}
where
\begin{eqnarray*}
r & = & \|\hat{h}'_2-\hat{h}'_1\|_{\H^1(0,\infty;\R^3)}
+ \| \hat{\omega}_2 -\hat{\omega}_1\|_{\H^1(0,\infty;\R^3)}
+ \left\| e^{\lambda t}\frac{\p X^{\ast}_2}{\p t}-e^{\lambda t}\frac{\p X^{\ast}_1}{\p t} \right\|_{\mathcal{H}(S_{\infty}^0)},
\end{eqnarray*}
and $K(r)$ is bounded when $r$ goes to $0$.
\end{lemma}

\begin{proof}
Let us write
\begin{eqnarray*}
X_2-X_1 & = & \overline{X}^R_2 \circ \overline{X}^{\ast}_2 - \overline{X}^R_2 \circ \overline{X}^{\ast}_1 +
\left( \overline{X}^R_2 - \overline{X}^R_1\right) \circ \overline{X}^{\ast}_1.
\end{eqnarray*}
For tackling the difference $\overline{X}^R_2 \circ \overline{X}^{\ast}_2 - \overline{X}^R_2 \circ \overline{X}^{\ast}_1$, let us apply Lemma~A.3 of the Appendix of \cite{Bourguignon}; We get the estimate
\begin{eqnarray*}
\| \overline{X}^R_2 \circ \overline{X}^{\ast}_2 - \overline{X}^R_2 \circ \overline{X}^{\ast}_1 \|_{\mathbf{H}^3(\mathcal{F})} & \leq &
C \| \overline{X}^R_2 \|_{C^4(\overline{\mathcal{F}})} \| \overline{X}^{\ast}_2 - \overline{X}^{\ast}_1 \|_{\mathbf{H}^3(\mathcal{F})} \times \\
& & \left(\|\overline{X}^{\ast}_1\|^3_{\mathbf{H}^3(\mathcal{F})} + \|\overline{X}^{\ast}_2\|^3_{\mathbf{H}^3(\mathcal{F})} + 1 \right),
\end{eqnarray*}
and thus the regularity in $\L^{\infty}(0,\infty;\mathbf{H}^3(\mathcal{F}))$. The regularity in $\W^{1,\infty}(0,\infty;\mathbf{H}^1(\mathcal{F}))$ can be also obtained by applying Lemma~A.3 of \cite{Bourguignon} for the time derivative of $\overline{X}^R_2 \circ \overline{X}^{\ast}_2 - \overline{X}^R_2 \circ \overline{X}^{\ast}_1$.\\
For the term $\left( \overline{X}^R_2 - \overline{X}^R_1\right) \circ \overline{X}^{\ast}_1$, we apply Lemma A.2 of the Appendix of \cite{Bourguignon}; We get the estimate
\begin{eqnarray*}
\left\| \left( \overline{X}^R_2 - \overline{X}^R_1\right) \circ \overline{X}^{\ast}_1 \right\|_{\mathbf{H}^3(\mathcal{F})} & \leq &
\| \overline{X}^R_2 - \overline{X}^R_1 \|_{C^3(\overline{\mathcal{F}})}\left( \|\overline{X}^{\ast}_1\|^3_{\mathbf{H}^3(\mathcal{F})}+1 \right),
\end{eqnarray*}
and thus the regularity in $\L^{\infty}(0,\infty;\mathbf{H}^3(\mathcal{F}))$. Here again the regularity in the space $\W^{1,\infty}(0,\infty;\mathbf{H}^1(\mathcal{F}))$ is obtained by applying the same lemma on the time derivative.
\end{proof}

\begin{proposition} \label{lemmaKtilde}
Let $\tilde{X}_1$ and $\tilde{X}_2$ be defined by
\begin{eqnarray*}
\tilde{X}_1 = \mathbf{R}^T_1(X_1-h_1), & \quad & \tilde{X}_2 = \mathbf{R}^T_2(X_1-h_2),
\end{eqnarray*}
where $X_1$ and $X_2$ are given in the assumptions of the previous lemma. Then
\begin{eqnarray}
\|\tilde{X}_2 - \tilde{X}_1 \|_{\tilde{\mathcal{W}}(Q_{\infty}^0)} & \leq & r\tilde{K}(r), \label{app11} \\
\| \nabla \tilde{Y}_2(\tilde{X}_2) - \nabla \tilde{Y}_1(\tilde{X}_1) \|_{\L^{\infty}(\mathbf{H}^2)\cap \W^{1,\infty}(\mathbf{L}^2)} & \leq & r\tilde{K}(r), \label{app12} \\
\| \Delta \tilde{Y}_2(\tilde{X}_2) - \Delta \tilde{Y}_1(\tilde{X}_1) \|_{\L^{\infty}(\mathbf{H}^1)} & \leq & r\tilde{K}(r), \label{app13}
\end{eqnarray}
where
\begin{eqnarray*}
r & = & \|\hat{h}'_2-\hat{h}'_1\|_{\H^1(0,\infty;\R^3)}
+ \| \hat{\omega}_2 -\hat{\omega}_1\|_{\H^1(0,\infty;\R^3)}
+ \left\| e^{\lambda t}\frac{\p X^{\ast}_2}{\p t}-e^{\lambda t}\frac{\p X^{\ast}_1}{\p t} \right\|_{\mathcal{H}(S_{\infty}^0)},
\end{eqnarray*}
and $\tilde{K}(r)$ is bounded when $r$ goes to $0$.
\end{proposition}

\begin{proof}
For proving \eqref{app11}, it is sufficient to write
\begin{eqnarray*}
\tilde{X}_2 - \tilde{X}_1 & = & \mathbf{R}^T_2(X_2-X_1) + \left(\mathbf{R}^T_2 - \mathbf{R}^T_1 \right)(X_1-h_1) - \mathbf{R}^T_2(h_2-h_1)
\end{eqnarray*}
and to apply the previous lemma. The estimate \eqref{app12} can be proven exactly like the estimate \eqref{estnablatilde}. Finally, for the estimate \eqref{app13} we denote by $\tilde{Y}_{i,1}$ and $\tilde{Y}_{i,2}$ the i-th component of $\tilde{Y}_1$ and $\tilde{Y}_2$ respectively, and we write the equality
\begin{eqnarray*}
\nabla^2 \tilde{Y}_{i,2}(\tilde{X}_2) - \nabla^2 \tilde{Y}_{i,1}(\tilde{X}_1) & = &
\left( \nabla \left( \nabla \tilde{Y}_{i,2}(\tilde{X}_2) - \nabla \tilde{Y}_{i,1}(\tilde{X}_1) \right) \right)\nabla \tilde{Y}_{2}(\tilde{X}_2) \\
& & + \left( \nabla \left( \tilde{Y}_{i,1}(\tilde{X}_1) \right) \right)\left( \nabla \tilde{Y}_2(\tilde{X}_2) - \nabla \tilde{Y}_1(\tilde{X}_1) \right),
\end{eqnarray*}
and apply Lemma \ref{lemmaGrubb}.
\end{proof}

\medskip
Received xxxx 20xx; revised xxxx 20xx.
\medskip

\end{document}